\documentclass[11pt]{article}
\usepackage{amsfonts}
\usepackage{color}
\newtheorem{theo}{Theorem}[section]
\newtheorem{lem}[theo]{Lemma}

\newtheorem{defi}[theo]{Definition}

\newcommand{\mysection}[1]{\section{#1} \setcounter{equation}{0}}
\newcommand{\proof}{{\sc Proof.} \quad}
\newcommand{\proofc}{{\sc Proof} \ }
\newcommand{\be}{\begin{equation} \label}
\newcommand{\ee}{\end{equation}}
\newcommand{\bea}{\begin{eqnarray}\label}
\newcommand{\eea}{\end{eqnarray}}
\newcommand{\bas}{\begin{eqnarray*}}
\newcommand{\eas}{\end{eqnarray*}}
\newcommand{\bit}{\begin{itemize}}
\newcommand{\eit}{\end{itemize}}
\newcommand{\qed}{\hfill$\Box$ \vskip.2cm}
\newcommand{\nn}{\nonumber}
\newcommand{\R}{\mathbb{R}}
\newcommand{\N}{\mathbb{N}}
\newcommand{\pO}{\partial\Omega}

\newcommand{\eps}{\varepsilon}

\newcommand{\supp}{{\rm supp} \, }

\newcommand{\wto}{\rightharpoonup}

\newcommand{\hra}{\hookrightarrow}
\newcommand{\io}{\int_\Omega}
\newcommand{\bom}{\overline{\Omega}}

\newcommand{\cb}{\normalcolor}
\newcommand{\cred}{\normalcolor}
\newcommand{\abs}{\\[5pt]}
\newcommand{\ueps}{u_\eps}
\newcommand{\weps}{w_\eps}
\newcommand{\deps}{d_\eps}
\newcommand{\uepsx}{u_{\eps x}}
\newcommand{\wepsx}{w_{\eps x}}
\newcommand{\depsx}{d_{\eps x}}
\newcommand{\wepsxx}{w_{\eps xx}}
\newcommand{\yeps}{y_\eps}
\newcommand{\heps}{h_\eps}
\newcommand{\tme}{T_{max,\eps}}

\newcommand{\ug}{\underline{\gamma}}
\newcommand{\og}{\overline{\gamma}}

\newcommand{\tphi}{\widetilde{\varphi}}
\newcommand{\hphi}{\widehat{\varphi}}

\newcommand{\epsin}{\eps\in (\eps_j)_{j\in\N}}

\newcommand{\zeps}{z_\eps}
\newcommand{\od}{\omega_d}
\begin{document}
\enlargethispage{10mm}
\title{Singular structure formation in a degenerate haptotaxis model involving myopic diffusion}
\author{
Michael Winkler\footnote{michael.winkler@math.uni-paderborn.de}\\
{\small Institut f\"ur Mathematik, Universit\"at Paderborn,}\\
{\small 33098 Paderborn, Germany} }
\date{}
\maketitle
\begin{abstract}
\noindent 
  We consider the system			%haptotaxis system
  \be{a1}
	\left\{ \begin{array}{l}
	u_t=\big(d(x)u\big)_{xx} - \big(d(x)uw_x\big)_x, \\[1mm]
	w_t=-ug(w),
	\end{array} \right.
  \ee
  which arises as a simple model for haptotactic migration in heterogeneous environments, 
  such as typically occurring in the invasive dynamics of glioma.
  A particular focus is on situations when the diffusion herein is degenerate in the sense that
  the zero set of $d$ is not empty.\abs 
  It is shown that if such possibly present degeneracies are sufficiently mild in the sense that
  \be{a2}
	\io \frac{1}{d}<\infty,
  \ee
  then under appropriate assumptions on the initial data a corresponding initial-boundary value problem for
  (\ref{a1}), posed under no-flux boundary conditions in a bounded open interval $\Omega\subset\R$,
  possesses at least one globally defined generalized solution. \\
  Moreover, despite such degeneracies the myopic diffusion mechanism in (\ref{a1})
  is seen to asymptotically determine the solution behavior
  in the sense that for some constant $\mu_\infty>0$, the obtained solution satisfies
  \be{a3}
	u(\cdot,t)\wto \frac{\mu_\infty}{d}
	\quad \mbox{in } L^1(\Omega)
	\qquad \mbox{and} \qquad
	w(\cdot,t) \to 0
	\quad \mbox{in } L^\infty(\Omega)
	\qquad \mbox{as } t\to\infty,
  \ee
  and that hence in the degenerate case the solution component $u$ stabilizes toward a state involving infinite densities,
  which is in good accordance with experimentally observed phenomena of cell aggregation.\\
  Finally, under slightly stronger hypotheses inter alia requiring that $\frac{1}{d}$ belong
  to $L\log L(\Omega)$, a substantial effect of diffusion is shown to appear already immediately by proving that
  for a.e.~$t>0$, the quantity $\ln (du(\cdot,t))$ is bounded in $\Omega$. 
  In degenerate situations, this particularly implies that the blow-up phenomena expressed in (\ref{a3})
  in fact occur instantaneously.\abs
\end{abstract}
\newpage
\section{Introduction}\label{intro}
In the theoretical description of collective cell behavior at macroscopic scales, 
taxis mechanisms have been playing an increasingly substantial role (\cite{hillen_painter}). 
In the past two decades, an accordingly growing literature on mathematical analysis of such processes 
has brought about quite a thorough knowledge of various classes of corresponding PDE models, 
containing cross-diffusive parabolic equations as their most characteristic ingredient,
especially in situations when 	
the attractive signal is a chemical and hence diffusible (see \cite{BBTW} for a recent survey).
{\cred
Unlike such chemotaxis systems, considerably less understood seem so-called haptotaxis systems 
which substantially differ from the former in that they 
address cases of non-diffusible cues, as naturally involved when tumors invade healthy tissue.
}
Moreover, virtually all analytical studies {\cred on taxis systems} 
assume that random movement of cells is of Fickian diffusion type,
either linear or nonlinear, with few exceptions considering fractional diffusion chemotaxis models
(\cite{bournaveas_calvez2010}, \cite{burczak_et_al_M3AS2016}).
Recent modeling approaches, however, indicate that in situations of significantly heterogeneous environments,
adequate macroscopic limits of random walks based on individually local sensing  
rather lead to certain non-Fickian diffusion operators (\cite{maini_et_al_JTB2013}, 
{\cred \cite{engwer_hunt_surulescu}}).\abs
The main focus of the present work is on the question how far the latter concept, {\cred in the literature
also referred to as {\em myopic diffusion} (\cite{maini_et_al_JTB2013}),}
can rigorously be proved appropriate for the description of spontaneous structure generation in the 
context of simple haptotaxis systems in heterogeneous environments.
We thereby intend to provide some analytical evidence for
heuristic reasonings (\cite{maini_et_al_JTB2013}) suggesting that in contrast to those based on Fickian diffusion, 
this modeling framework may indeed much more accurately
describe the emergence of neighborhood-adapted structures in such populations of myopic individuals,
with aggregation phenomena of {\cb glioma} near {\cred thin} interfaces {\cb between} white and grey matter in mouse brains
forming a corresponding experimental observation of particular importance (\cite{burden_gulley2011}).\abs
{\cred
To this end, we will consider a particular version of an evolution system recently proposed as a model 
for the description of glioma spread in heterogeneous tissue (\cite{engwer_hunt_surulescu}), for mathematical
purposes simplified in that any proliferation effects are neglected and that 
the spatial setting is assumed to be one-dimensional.
Specifically, we shall be concerned with the initial-boundary value problem
}
\be{0}
	\left\{ \begin{array}{ll}
	u_t=\big(d(x)u\big)_{xx} - \big(d(x)uw_x\big)_x,
	\qquad & x\in\Omega, \ t>0, \\[1mm]
	w_t=-ug(w),
	\qquad & x\in\Omega, \ t>0, \\[1mm]
	\big(d(x)u\big)_x - d(x)uw_x=0, 
	\qquad & x\in\pO, \ t>0, \\[1mm]
	u(x,0)=u_0(x), \quad w(x,0)=w_0(x),
	\quad & x\in\Omega,
	\end{array} \right.
\ee
{\cred for the unknown cell density $u=u(x,t)$ and} the {\cb density} $w=w(x,t)$ {\cb of tissue fibers acting 
as a haptotactic cue} {\cred in a bounded open interval $\Omega\subset \R$,} 
with given nonnegative functions $d$, $u_0$ and $w_0$ on $\bom$ and 
$g$ generalizing the {\cred prototypical} 
choice $g(s)=s$, $s\ge 0$, in a sense to be specified in (\ref{g1}) and (\ref{g2}) below.
{\cb
The } {\cred formal parabolic limit procedure performed 
in \cite[Section 3.1]{engwer_hunt_surulescu},
adequately accounting for the} {\cb influence of the 
underlying tissue structure} 
{\cred on} {\cb tumor cell movement, led to the above concrete form of the macroscopic equations featuring myopic 
diffusion and haptotaxis. Both these types of terms in their respective coefficient functions, involve 
the so-called tumor diffusion tensor explicitly deduced e.g.~in \cite[Formula (3.11)]{engwer_hunt_surulescu}. 
In the latter reference, the distribution of the tissue density is assessed from medical data and plays 
the role of an input to the equation for the space-time evolution of the tumor cell population. 
When the tissue dynamics is taken into account, as done through the second equation in (\ref{0}),
then the mathematical analysis 
of the resulting system becomes challenging, the more so in situations with possibly 
degenerate diffusion, which can indeed occur during the migration of glioma through the tissue, 
either when the latter is 
locally too dense and isotropic, thus impairing the spread of cells which have to overcome it, or 
too sparse, which in turn is hindering the spread of cells, as they have to rely on 
it both for migration and proliferation. 
In this work we therefore concentrate on the one-dimensional version of 
the system obtained in \cite{engwer_hunt_surulescu}, which 
correspondingly uses the same motility coefficient function $d=d(x)$ in both 
the diffusive and the advective terms in (\ref{0}) and allow this function to degenerate.}\abs
We {\cred moreover note} that as can readily be verified on substituting $w=\Psi(v):=\int_0^v \psi(\xi)d\xi$ and 
$g(w)=\Psi^{-1}(w) \cdot \psi(\Psi^{-1}(w))$, for arbitrary
smooth positive $\psi: \ [0,\infty)\to\R$
this thereby implicitly includes solutions with sufficiently small component $v$
of the respective initial-boundary value problem for
\bas
	\left\{ \begin{array}{ll}
	u_t=\big(d(x)u\big)_{xx} - \big(d(x)u\psi(v)v_x\big)_x,
	\qquad & x\in\Omega, \ t>0, \\[1mm]
	v_t=-uv,
	\qquad & x\in\Omega, \ t>0, 
	\end{array} \right.
\eas
where the choice $\psi(v)=\frac{1}{(1+v)^2}$ corresponds to the particular tumor invasion model
recently {\cred analyzed} in \cite{zhigun_surulescu_uatay}.\abs
In the case $d\equiv 1$ representing spatially homogeneous conditions for both diffusion and cross-diffusion,
(\ref{0}) reduces to the apparently simplest reasonable model for haptotactic interaction
(\cite{levine_sleeman_nilsenh_JOMB2001}), containing the essential
aspects of several more complex systems that have been discussed in the modeling literature (\cite{perumpanani_byrne},
\cite{chaplain_anderson2003}, \cite{CL2}; cf.~also \cite{bellomo_li_maini_M3AS2008}) and also analyzed analytically.
Beyond statements on global existence in various functional frameworks 
(see \cite{walker_webb}, \cite{corrias_perthame_zaag2004} and \cite{litcanu_moralesrodrigo_M3AS} for some classical and 
e.g.~\cite{tao_NA2011} as well as \cite{ssw} for more recent examples)
and scattered results on boundedness (\cite{marciniak_ptashnyk}, \cite{taowin_PROCA2014},
\cite{friedman_tello}), however,
even in this non-degenerate and homogeneous setting 
a detailed description of further qualitative facets such as the large time behavior could be established
only in very particular cases up to now;
moreover, apparently all available results in this direction are either restricted to solutions suitably close to
equilibria (\cite{fontelos_friedman_hu}, \cite{friedman_tello}), or to situations when a strongly dissipative
action of additional logistic-type cell kinetic
terms can be shown to dominate on large time scales (\cite{li_lin_mu_AML2015},
\cite{yifu_wang_ke_JDE2016}, \cite{taowin_SIMA2015}, \cite{taowin_PROCA2014}, \cite{hpw_M3AS}),
meaning that in the latter cases solutions exclusively stabilize toward spatially 
homogeneous and hence unstructured equilibria.
This lack of rigorous knowledge in situations of expectedly more colorful solution behavior 
may be viewed as reflecting the circumstance that unless suitably compensated by further mechanisms,
tactic cross-diffusion of the form in (\ref{0}) 
may substantially affect the regularity of solutions and hence 
obstruct mathematical analysis at various stages.
This strongly destabilizing potential is well-known from various findings detecting 
unboundedness phenomena especially in self-reinforced taxis models,
{\cred 
even in apparently more regular settings determined by cross-diffusive interaction with a diffusible quantity
such as in}
the classical Keller-Segel chemotaxis
system and derivates thereof (\cite{herrero_velazquez}, \cite{nagai}, \cite{win_JMPA},
\cite{lankeit_threshold}),
{\cred already} in some spatially one-dimensional scenarios (\cite{kang_stevens_velazquez_CPDE2010}, \cite{win_JNLS}),
but also in some models for tactic migration toward non-diffusible attractants (\cite{levine_sleeman_SIMA1997},
\cite{sleeman_levine_MMAS2001}).\abs
{\bf Main results.}\quad
In the presently considered context of the model (\ref{0}), our analysis will reveal 
that under appropriate assumptions inter alia requiring mildness of possible degeneracies in diffusion,
such types of taxis-driven collapse do not occur, but that the
solution behavior is rather essentially prearranged by the environmental conditions.
Indeed, our main results will show that for a large class of initial data, 
certain global generalized solutions can be constructed which in the large time limit
approach a positive multiple of the reciprocal myopic diffusion coefficient $\frac{1}{d}$ in their first
component, as predicted in \cite{maini_et_al_JTB2013};
in particular, this reflects asymptotic aggregation of cells in regions where $d$ is small, 
in presence of zeros of $d$ even in the mathematically extreme sense of stabilization toward a singular state.
Beyond this asymptotic statement, we will identify a solution property that indicates
a certain predominance of the diffusion process in (\ref{0}) already at intermediate and even small time scales:
Namely, we shall see that under slightly stronger assumptions,	
for a.e.~$t>0$ 
the quantity $du(\cdot,t)$ is bounded from above and below in $\Omega$ by positive constants only depending on $t$.
This firstly ensures local boundedness of $u(\cdot,t)$ inside the positivity set of $d$ and hence rules out
any significant taxis-forced aggregation; secondly, and more drastically, however, this implies that 
singularities near points of degenerate diffusion, according to the above arising at least
in the long-term limit, in fact emerge instantaneously.\abs
In order to formulate these results more precisely, let us specify the framework to be considered henceforth by
assuming $d\in C^0(\bom)$ to be nonnegative and such that
\be{d1}
	d\in C^1(\{d>0\}),
\ee
as well as
\be{d2}
	\io \frac{1}{d} < \infty,
\ee
where $\{d>0\} := \big\{ x\in \overline{\Omega} \ \big| \ d(x)>0 \big\}$,
with this and similar notation frequently being used throughout the sequel without further explicit definition.
{\cred
We observe that (\ref{d2}) in particular requires the set of all zeros of $d$ to be a null set of points, thus
inter alia excluding situations when diffusion may become degenerate throughout entire subintervals of $\Omega$.
In application contexts, this corresponds to limiting situations of small interfacial layers
of inhibited diffusion, such as typically occurring in the mentioned framework of glioma spread addressed in
\cite{maini_et_al_JTB2013}.}
{\cred Mathematically, it may be noted that at least formally, (\ref{0}) would predict temporal constancy of $u$ 
inside the interior of such degeneracy regions; a partial rigorous justification thereof has recently been achieved in
\cite{surulescu_win1}.} 
{\cb Thinking of 	
the particular problem setting of glioma invasion,
let us recall that the 
tumor diffusion tensor obtained during the macroscopic scaling process in \cite{engwer_hunt_surulescu} is 
proportional to the water diffusion tensor assessed by diffusion tensor imaging (cf.~also
\cite{clatz} and \cite{swanson} for independently obtained similar links);
accordingly, in the one-dimensional framework at hand the corresponding scalar coefficient function $d$
is also supposed to be tightly related to the diffusivity of water molecules. 
Thereby, 
sharp intersections 
of the one-dimensional diffusion direction of water molecules by tissue fibers which are very thin, single objects, 
at those sites lead
to essentially single-point degeneracies in the diffusion of water molecules and, the more so, of tumor cells.}
{\cb In addition, (\ref{d2}) implicitly requires that $d$ grows suitably fast near its zeros, 
in the prototypical case when $d(x) = |x-x_0|^\theta$ for all $x\in\bom$, some $x_0\in\bom$ and some 
$\theta>0$ reducing to the hypothesis 
that $\theta<1$.
Biologically, this corresponds to situations in which the diffusivity 
undergoes a rapid enhancement 
in the immediate proximity of the sites of degeneration, e.g., where it was 'blocked' by the fibers
(\cite{engwer_hunt_surulescu});
further indications for the occurrence of such sudden increases in diffusivity at interfaces
is provided by experimental evidence reporting that the diffusivity in white brain matter is much higher than in 
grey matter and leads 
to differences in cell motility 5-25 times higher in white than in grey matter (see e.g. 
\cite{maini_et_al_JTB2013} and the references therein).
Besides their biological plausibility, these assumptions will also serve technical purposes} 
{\cred that will become evident in the discussion below, e.g.~around
the formulation of Theorem \ref{theo7000}.}\abs
As for the signal absorption coefficient function in (\ref{0}),	
we shall suppose that $g\in C^2([0,\infty))$ is such that
$g(0)=0$ and that with some positive constants $\ug$ and $\og$ we have
\be{g1}
	\ug \le g'(s) \le \og
	\qquad \mbox{for all } s\ge 0
\ee
and hence also
\be{g2}
	\ug s \le g(s) \le \og s
	\qquad \mbox{for all } s\ge 0,
\ee
and the intial data are required to be such that
\be{init}
	\left\{ \begin{array}{l}
{\cred	0\le u_0\in C^0(\overline{\Omega})  \mbox{ satisfies } u_0\not\equiv 0, }
	\quad \mbox{and that} \\[1mm]
	0\le w_0 \in C^0(\overline{\Omega}) \ \mbox{ is such that } 
	\sqrt{w_0} \in W^{1,2}(\Omega).	
	\end{array} \right.
\ee
Within this setting, the first of our main results establishes global existence of a solution to
(\ref{0}) under an appropriate additional condition requiring a certain smallness property of $w_0$ near zeros of $d$.
We emphasize already here that due to our mild assumptions on $d$, 
in view of the statement on instantaneous blow-up formulated in Theorem \ref{theo87}
we can in general not expect boundedness of the first solution component with respect to the norm in $L^p(\Omega)$
for any $p>1$, not even locally in time,
so that our notion of solution needs to be adequately adapted to this circumstance.
After all, our analysis will reveal that it is not necessary to resort to concepts involving measure-valued solutions,
but that it is rather possible to construct solutions with their first component
belonging to the space $C^0_w([0,\infty);L^1(\Omega))$ of $L^1(\Omega)$-valued functions defined on $[0,\infty)$
which are continuous with respect to the weak topology in $L^1(\Omega)$.
\begin{theo}\label{theo700}
  Let $\Omega\subset \R$ be a bounded interval, and suppose that $d\in C^0(\overline{\Omega})$ is nonnegative
  and such that (\ref{d1}) and (\ref{d2}) hold.
  Moreover, let $g\in C^2([0,\infty))$ be such that $g(0)=0$ and that (\ref{g1}) is valid with some $\ug>0$ and $\og>0$.
  Then for all initial data $u_0$ and $w_0$ which satisfy (\ref{init}) and which are such that furthermore
  \be{w0}
	\io \frac{d_x^2}{d} w_0 < \infty,
  \ee
  there exists at least one pair $(u,w)$ of nonnegative functions 
  \be{700.01}
	\left\{ \begin{array}{l}
	u \in C^0_w([0,\infty);L^1(\Omega)) \cap L^\infty((0,\infty);L^1(\Omega)), \\[1mm]
	w \in C^0(\bom\times [0,\infty)) \cap L^\infty(\Omega\times (0,\infty)) \cap L^1_{loc}([0,\infty);W^{1,1}(\Omega)),
  	\end{array} \right.
  \ee
  which form a global weak solution of (\ref{0}) in the sense of Definition \ref{defi_weak}, and for which we have
  \be{700.9}
	\io u(\cdot,t)=\io u_0
	\qquad \mbox{for all } t>0.
  \ee
\end{theo}
Next, our main result concerning qualitative
behavior in (\ref{0}) asserts that in the large time limit, each of these solutions approaches a steady state of (\ref{0}).
Here since the nonnegative equilibria of (\ref{0}) are precisely the pairs $(\frac{\mu}{d},0)$ with $\mu\ge 0$,
in light of the mass conservation property (\ref{700.9}) this a posteriori underlines the crucial role
of our overall integrability assumption (\ref{d2}) for this central result.
\begin{theo}\label{theo7000}
  Suppose that the assumptions of Theorem \ref{theo700} are fulfilled. Then the global generalized solution
  $(u,w)$ of (\ref{0}) obtained in Theorem \ref{theo700} satisfies
  \be{7000.2}
	u(\cdot,t)\wto \frac{\mu_\infty}{d}
	\quad \mbox{in } L^1(\Omega)
	\qquad \mbox{as } t\to\infty
  \ee
  and
  \be{7000.3}
	w(\cdot,t) \to 0
	\quad \mbox{in } L^\infty(\Omega)
	\qquad \mbox{as } t\to\infty
  \ee
  with the positive number
  \be{mu_infty}
	\mu_\infty:=\frac{\io u_0}{\io \frac{1}{d}}.
  \ee
\end{theo}
{\cred
We note that in presence of zeros of $d$, (\ref{7000.2}) actually asserts that the quantity $u$ undergoes a
certain blow-up phenomenon at least in the large time limit.
We finally make sure that this explosion
}
actually occurs immediately and persistently, provided that diffusion is slightly less degenerate than admitted
in Theorem \ref{theo700}, and that $\frac{w_0}{d}$ is bounded.
In fact, the following states that under these hypotheses, the regularizing action of diffusion is strong enough,
both relatively to haptotaxis and absolutely, so as to allow for the conclusion that, at least in 
an appropriate weakened form, the quantity $du$ enjoys
properties of instantaneous positivity and boundedness well-known for solutions of the heat equation.
\begin{theo}\label{theo87}
  Assume that in addition to the hypotheses of Theorem \ref{theo700},
  \be{87.1}
	\io \frac{1}{d} \ln \frac{1}{d} < \infty
  \ee
  and
  \be{87.3}
	\frac{w_0}{d} \in L^\infty(\Omega).
  \ee
  Then the global generalized solution $(u,w)$ of (\ref{0}) from Theorem \ref{theo700} has the property that
  \be{87.4}
	\int_\tau^T \Big\| \ln \Big(du(\cdot,t)\Big)\Big\|_{L^\infty(\Omega)}^3 dt < \infty
	\qquad \mbox{for all $T>0$ and } \tau \in (0,T).
  \ee
  In particular, for a.e.~$t>0$ there exist $C_1(t)>0$ and $C_2(t)>0$ such that
  \be{87.5}
	\frac{C_1(t)}{d(x)} \le u(x,t) \le \frac{C_2(t)}{d(x)}
	\qquad \mbox{for a.e. } x\in\Omega,
  \ee
  and if $\frac{1}{d} \not\in L^\infty(\Omega)$, then
  \be{87.6}
	\|u(\cdot,t)\|_{L^\infty(\Omega)} = \infty
	\qquad \mbox{for a.e. } t>0.
  \ee
\end{theo}
{\cred 
The above results seem to go beyond previous knowledge even in cases when haptotactic interaction is neglected
e.g.~by formally setting $w\equiv 0$ in (\ref{0}). 
In the non-degenerate version of the correspondingly obtained linear diffusion problem, that is, when $d>0$
in $\bom$, global existence of classical solutions, smoothly approaching the steady state in (\ref{7000.2}),
can readily be established by standard methods.
As for degenerate limit cases thereof, a result on global existence of certain very weak solutions, as well as on
their stabilization toward an associated singular equilibrium, can be found in \cite{hpw_EJAM}.
A very early caveat indicating criticality of the assumption (\ref{d2}) goes back to \cite{feller},
where it is shown that if the diffusion degeneracy is slightly stronger in that $d(x)=x$ in $\Omega=(0,\infty)$,
then prescribing boundary conditions at $x=0$ in the resulting simple equation $u_t=(xu)_{xx}$
is meaningless in the sense that solutions to the initial-value problem therefor are uniquely determined 
already by their prescribed (reasonably regular) initial data.\abs
}
{\bf Main ideas.} \quad
Our analysis is rooted in the observation that in the context of non-degenerate and suitably regular
diffusion, a supposedly given smooth solution to (\ref{0}) satisfies the energy inequality
\bea{energy}
	\frac{d}{dt} \bigg\{
	\io u\ln (du)
	+ \frac{1}{2} \io d\frac{w_x^2}{g(w)}
	+ \frac{\og}{\ug^2} \io \frac{d_x^2}{d} w \bigg\} 
	+ \io \frac{(du)_x^2}{du}
	+ \frac{\ug}{4\og} \io du \frac{w_x^2}{w}
	\le 0
\eea
where our hypothesis that $\io \frac{1}{d}$ be finite warrants that the Lyapunov functional therein 
is bounded from below (cf.~Lemma \ref{lem8}). 
Thus generalizing the corresponding identity for the special case $d\equiv const.$, as already observed
in \cite{corrias_perthame_zaag2004} and frequently adapted to various related cases involving
spatially homogeneous diffusion (cf.~\cite{ssw} for a recent even quite complex example), 
(\ref{energy}) contains in its dissipated part, as a main novel ingredient,
the fraction $\frac{d_x^2}{d}$ which our assumption (\ref{d2}) 
enforces to have infinite integral around each zero of $d$ (see Lemma \ref{lem41}).
Mainly due to this circumstance, considerable efforts will be undertaken in Section \ref{sect2} 
to carefully design a sequence of regularized problems, indexed by a small positive parameter $\eps$,
that will involve nondegenerate diffusion in the respective first equation 
as well as a parabolic approximation of the second equation in (\ref{0}), 
and at the core of which
the construction of suitable approximations $\deps$ and $w_{0\eps}$ to $d$ and $w_0$, respectively,
is guided by the intention to remain basically consistent with the structure expressed in (\ref{energy}).
In Section \ref{sect3} this will enable us to obtain an
approximate counterpart of (\ref{energy}) and derive correspondingly implied a priori estimates for 
the respective solutions $(\ueps,\weps)$ in the central Lemma \ref{lem8}, inter alia containing a 
regularized variant of the global dissipation property
\be{relax0}
	\int_0^\infty \io \frac{(du)_x^2}{du} <\infty
\ee
formally resulting from (\ref{energy}).
By means of standard testing procedures, in Section \ref{sect4} these will be seen to
entail further regularity properties, now possibly $\eps$-dependent, which enable us to extend each of these approximate
solutions so as to exist globally.\\
Beyond some local-in-time estimates for $u_{\eps x}$ and $w_{\eps t}$, Section \ref{sect_equi} will thereafter
reveal two key regularity features, namely firstly uniform integrability of $\ueps$ and of $\wepsx$ with respect to
both the time variable and the approximation parameter (Lemma \ref{lem826} and Lemma \ref{lem824}), 
and secondly an approximate analogue of the relaxation property
\be{relax}
	\int_0^\infty \|u_t(\cdot,t)\|_{W^{1,\infty}(\Omega))^\star}^2 dt < \infty
\ee
formally implied by (\ref{energy}) (Lemma \ref{lem27}).
Along with a crucial strong $L^2$ compactness property of the first factor $\sqrt{\deps}\ueps$ in the corresponding
cross-diffusive flux (Lemma \ref{lem823} and Lemma \ref{lem822}),
these will allow for constructing a solution to (\ref{0}) through an appropriate extraction procedure
based on straightforward compactness arguments (Section \ref{sect6}), and thus for proving Theorem \ref{theo700}
(Section \ref{sect7}).\\
Section \ref{sect8} will then be devoted to the derivation of the stabilization results in Theorem \ref{theo7000},
where first concentrating on the solution component $u$ we will make
essential use of the weak decay information implicitly contained in (\ref{relax0}) and (\ref{relax}), as well
as a now evident equi-integrability feature of $(u(\cdot,t))_{t>0}$ (Sections \ref{sect8.1}-\ref{sect8.3}).	
Thereafter, the fact that thus $u$ approaches a positive limit will be combined with the equicontinuity of
$(w(\cdot,t))_{t>0}$, as implied by the above, to verify that the 
decreasing quantity $\|w(\cdot,t)\|_{L^\infty(\Omega)}$ must actually decay (Section \ref{sect8.4}).\abs
Finally, Section \ref{sect9} provides a proof of Theorem \ref{theo87}, with a key step consisting in deriving an
estimate of the form
\be{wx}
	\int_0^T \io \wepsx^2 \le C(T),
	\qquad T>0,
\ee
(Lemma \ref{lem82}), used to control the right-hand side in the regularized analogue of
\be{ln}
	\frac{d}{dt} \io \frac{1}{d} \ln u
	\ge \frac{1}{2} \io \frac{(du)_x^2}{(du)^2} - \frac{1}{2} \io w_x^2
\ee
adequately (Lemma \ref{lem835}).
For smooth solutions, (\ref{wx}) would trivially result as a by-product of (\ref{energy}) due to the evident fact
that as a consequence of (\ref{0}) and 
the assumptions in Theorem \ref{theo87}, $\frac{d}{w}$ would have a positive lower bound,
and hence would $\frac{d}{g(w)}$ by (\ref{g2}).
Due to positivity of $\weps$ enforced by artificial diffusion, however, a corresponding upper bound
for $\frac{\weps}{\deps}$ seems available only in certain $L^p$ spaces, with the integrability power $p$ herein
fortunately increasing with decreasing $\eps$, however (Lemma \ref{lem81}). 
Therefore, (\ref{wx}) can be obtained by means of a subtle interpolation argument (Lemma \ref{lem82}) involving
an additional regularity information on $\wepsx$ which stems from the artificially introduced dissipation and is
thus of higher order, but singular with respect to $\eps$ (Lemma \ref{lem80}).\abs
Before going into details, let us remark that due to the delicate coupling of diffusion and haptotactic
cross-diffusion in (\ref{0}), in the general framework determined by our conditions and especially by (\ref{d2})
we do not expect solutions to possess spatially global
regularity properties substantially beyond those obtained by our analysis,
as already discussed above in the context of Theorem \ref{theo700}. 
An interesting question going beyond the scope of the present work consists in describing possible further
regularity aspects inside the positivity region of $d$ where in the purely diffusive case when 
$w\equiv 0$, standard parabolic theory essentially provides smoothness up to an extent determined by the smoothness
of $d$ and $u_0$.
After all, a subsequent study in this direction will inter alia show that 
imposing the slightly stronger assumption $\io \frac{1}{d^2}<\infty$ 
on the behavior of $d$ near its zeros ensures that the quantity $du$ remains bounded in $L^p(\Omega)$ for any 
$p\in (1,\infty)$, that locally in $\{d>0\}\times (0,\infty)$
the function $u$ itself is even H\"older continuous, and that the convergence in (\ref{7000.2}) 
in fact is locally uniform in $\{d>0\}$ (\cite{stiwi_hapto}).\abs

\mysection{Approximation of (\ref{0}) by a family of regularized problems}\label{sect2}		
\subsection{A weak solution concept}
To begin with, let us specify our generalized solution concept 
in order to substantiate the goal to be pursued in the context of our existence analysis.
\begin{defi}\label{defi_weak}
  A pair $(u,w)$ of nonnegative functions
  \be{w1}
	\left\{ \begin{array}{l}
	u \in L^1_{loc}(\overline{\Omega}\times [0,\infty)), \\[1mm]
	w\in L^\infty_{loc}(\overline{\Omega}\times [0,\infty)) \cap L^1_{loc}([0,\infty);W^{1,1}(\Omega))
	\end{array} \right.
  \ee
  satisfying
  \be{w2}
	duw_x \in L^1_{loc}(\overline{\Omega}\times [0,\infty))
  \ee
  will be called a {\em global weak solution} of (\ref{0}) if
  \bea{w3}
	- \int_0^\infty \io u\varphi_t - \io u_0 \varphi(\cdot,0)
	= \int_0^\infty \io du\varphi_{xx}
	+ \int_0^\infty \io duw_x \varphi_x 
  \eea
  for all $\varphi\in C_0^\infty(\overline{\Omega}\times [0,\infty))$ such that $\varphi_x=0$ on $\pO \times (0,\infty)$
  and
  \be{w4}
	\int_0^\infty \io w\varphi_t + \io w_0\varphi(\cdot,0)
	= \int_0^\infty \io ug(w) \varphi
  \ee	
  for all $\varphi\in C_0^\infty(\overline{\Omega}\times [0,\infty))$.
\end{defi}
\subsection{Construction of energy-compatible sequences approximating $d$ and $w_0$}\label{sect2.2}
A natural first step in the construction of globally defined functions solving (\ref{0}) in the above sense
consists in considering appropriately regularized {\cb problems.
In order to allow for classical solvability, the latter should in particular involve non-degenerate 
diffusion in the respective cruicial first equation;
as smooth solvability furthermore seems to require second-order spatial differentiability of the haptoattractant therein,
apart from that a certain smoothness-enforcing regularization in the second equation appears to be in order.}
In the context of the questions addressed here, however, 
nearby approaches based e.g.~on straightforward introduction of artificial
non-degenerate diffusion in both sub-problems of (\ref{0}) apparently need to face two essential challenges:
Firstly, our assumption (\ref{d2}) of suitably weak degeneracy implicitly forces $d$ to be non-smooth near possible zeros;
in particular, the function $d_x$ appearing as a coefficient in the divergence-like reformulation
of the diffusion operator $(du)_{xx}=(du_x+d_x u)_x$ need not belong to any of the spaces $L^p(\Omega)$ for $p>1$;
accordingly, for guaranteeing the existence of suitably smooth solutions to our regularized problems
it seems adequate to approximate $d$ by appropriate functions each of which, beyond being strictly positive, 
is sufficiently regular.
Secondly, and more drastically, in view of our goal to exploit the energy structure (\ref{energy}) 
formally associated with (\ref{0}), 
unlike in situations when only global solvability is strived for (\cite{surulescu_win1})
our design of regularization will be restricted to approximate
problems which are essentially consistent with this structure.
Here, in view of a considerably strong singularity of $\frac{d_x^2}{d}$ necessarily appearing 
near any zero of $d$ (Lemma \ref{lem41}),
a particularly crucial role will be played by the last integral $\io \frac{d_x^2}{d} w$ arising in the 
Lyapunov functional in (\ref{energy}), especially at the initial time where it seems far from obvious 
how far our mere assumptions in (\ref{init}) and
(\ref{w0}) may warrant boundedness of the respective expression when $d$ is replaced by 
approximate variants; accordingly, our regularization procedure will moreover include a suitable modification of
$w_0$ near zeros of $d$.\abs
In order to adequately cope with both these challenges, 
{\cred in this section we describe a possible
}
construction of a sequence of approximate versions of (\ref{0}), indexed by a small parameter $\eps\in (0,1)$ 
which will eventually be restricted so as to run along an appropriately chosen decreasing 
sequence $(\eps_j)_{j\in\N}\subset (0,1)$ (see Lemma \ref{lem45}).
{\cred In order to avoid abundant technicalities at this stage, we postpone details of the corresponding
analysis to an appendix below.}\abs
As a first {\cred step within our procedure, we will} make sure that $d$ can monotonically 
be approximated by a family of smooth positive functions $\deps$ with convenient further properties.
\begin{lem}\label{lem43}
  Suppose that $d\in C^0(\overline{\Omega})$ is such that (\ref{d1}) holds.
  Then there exists a family $(\deps)_{\eps\in (0,1)} \subset C^\infty(\overline{\Omega})$ with the properties that
  as $\eps\searrow 0$ we have
  \be{43.1}
	\deps \to d
	\qquad \mbox{in } L^\infty(\Omega)
  \ee
  and
  \be{43.02}
	\depsx \to d_x
	\qquad \mbox{in $L^\infty_{loc}(\{d>0\}\cap\Omega)$ and in $L^p_{loc}(\{d>0\})$ for all } p\in [1,\infty),
  \ee
  that
  \be{43.01}
	\deps \le d_{\eps'}
	\quad \mbox{in $\Omega$ whenever } 0<\eps\le\eps'<1,
  \ee
  that for all $\eps\in (0,1)$ we have $\deps>0$ in $\overline{\Omega}$,
  \be{43.2}
	\depsx =0
	\qquad \mbox{on } \pO
  \ee
  and
  \be{43.3}
	\deps \le \|d\|_{L^\infty(\Omega)}+1
	\qquad \mbox{in } \Omega,
  \ee
  and such that
  \be{43.5}
	\eps^2 \io \frac{\depsx^2}{\deps^3} \le 1
  \ee
  and
  \be{43.55}
	\sqrt{\eps} \io \frac{\depsx^4}{\deps^2} \le 1
  \ee
  as well as
  \be{43.6}
	\eps^\frac{1}{4} \cdot \frac{1}{\inf_{x\in\Omega} \deps(x)} \le 1
  \ee
  and
  \be{43.99}
	\eps^\frac{1}{4} \cdot \Big\| \frac{\depsx}{\deps}\Big\|_{L^\infty(\Omega)} \le 1
  \ee
  for all $\epsin$.
\end{lem}
In view of (\ref{43.1}) and (\ref{43.02}), taking $\eps\searrow 0$ in the expression $\io \frac{\depsx^2}{\deps} w_0$
will not go along with any difficulty in the special case when $w_0$ has compact support in $\{d>0\}$.
That it is reasonable to use such functions for the approximation of a general $w_0$, 
beyond the required regularity assumptions merely satisfying (\ref{w0}),
is indicated by the observation to be made in Lemma \ref{lem411}, which itself is prepared by the
following implication of our assumptions on $d$.
\begin{lem}\label{lem41}
  Let $d\in C^0(\overline{\Omega})$ be nonnegative and such that (\ref{d1}) and (\ref{d2}) are satisfied.
  Then for any set $\Omega_0\subset \overline{\Omega}$ which is relatively open in $\overline{\Omega}$ and such that
  $\Omega_0 \cap \{d=0\} \ne \emptyset$, we have
  \be{41.2}
	\int_{\Omega_0} \frac{d_x^2}{d} = \infty.
  \ee
\end{lem}
We can thereby easily assert that any $w_0$ compatible with the hypotheses of Theorem \ref{theo700} indeed must vanish
at each zero of $d$.
\begin{lem}\label{lem411}
  Let $d\in C^0(\overline{\Omega})$ be nonnegative and such that (\ref{d1}) and (\ref{d2}) are valid,
  and suppose that $w_0\in C^0(\overline{\Omega})$ is a nonnegative function fulfilling (\ref{w0}).
  Then 
  \be{411.3}
	w_0(x)=0
	\qquad \mbox{for all } x\in \{d=0\}.
  \ee
\end{lem}
We shall next use the above fact together with our overall regularity assumption that $\sqrt{w_0}$ belongs to
$W^{1,2}(\Omega)$ to construct a monotone sequence of approximations to $w_0$ which are all compactly supported
in $\{d>0\}$, and which moreover are compatible with the energy functional in (\ref{energy}) in the sense that
not only the third but also the second intergal therein remains bounded along this sequence.
\begin{lem}\label{lem42}
  Assume that the nonnegative function $d\in C^0(\overline{\Omega})$ satisfies (\ref{d1}) and (\ref{d2}),
  and that $w_0$ complies with (\ref{init}) and (\ref{w0}).
  Then there exists $(w_{0j})_{j\in\N} \subset L^\infty(\Omega)$ such that for all $j\in\N$ we have
  $w_{0j}\ge 0$ in $\Omega$ and $\sqrt{w_{0j}} \in W^{1,2}(\Omega)$ as well as
  \be{42.4}
	\supp w_{0j} \subset \{d>0\},
  \ee
  and such that
  \be{42.5}
	w_{0j} \nearrow w_0
	\quad \mbox{in } \Omega 
	\qquad \mbox{as } j\to\infty
  \ee
  and
  \be{42.6}
	\sup_{j\in\N} \io d\frac{w_{0jx}^2}{w_{0j}} < \infty.
  \ee
\end{lem}
We finally combine the outcomes of Lemma \ref{lem43} and Lemma \ref{lem42} to select a suitable
decreasing sequence $(\eps_j)_{j\in\N}\subset (0,1)$ along which the interplay of the correspondingly defined
function $d_{\eps_j}$ with a slightly shifted variant of $w_{0j}$ is favorable with 
regard to both relevant integrals appearing in the Lyapunov functional in (\ref{energy}).
\begin{lem}\label{lem45}
  Let $d\in C^0(\overline{\Omega})$ be nonnegative and such that (\ref{d1}) and (\ref{d2}) hold,
  and let $w_0$ satisfy (\ref{init}) and (\ref{w0}).
  Then there exists $(\eps_j)_{j\in\N}\subset (0,1)$ such that $\eps_j\searrow 0$ as $j\to\infty$, and such that
  for $(d_{\eps_j})_{j\in\N}$ as determined by Lemma \ref{lem43}, and for
  \be{45.3}
	w_{0\eps}(x):=w_{0j}(x) + \eps^\frac{1}{4},
	\qquad x\in\overline{\Omega}, \quad \eps=\eps_j, \quad j\in\N,
  \ee
  with $(w_{0j})_{j\in\N}$ taken from Lemma \ref{lem42}, we can find $C>0$ such that
  \be{45.4}
	\io \deps \frac{w_{0\eps x}^2}{w_{0\eps}} \le C
	\qquad \mbox{for all } \epsin
  \ee
  and
  \be{45.5}
	\io \frac{\depsx^2}{\deps} w_{0\eps} \le C
	\qquad \mbox{for all } \epsin.
  \ee
\end{lem}
\subsection{Regularized problems: local existence and extensibility}
Upon the choices specified in Lemma \ref{lem45}, for $\eps\in (\eps_j)_{j\in\N}$ we henceforth consider the
approximate variants of (\ref{0}) given by
\be{0eps}
	\left\{ \begin{array}{ll}
	u_{\eps t} = (\deps\ueps)_{xx} - (\deps \ueps w_{\eps x})_x,
	\quad & x\in\Omega, \ t>0, \\[1mm]
	w_{\eps t} = \eps \Big( \deps \frac{w_{\eps x}}{\sqrt{g(\weps)}} \Big)_x - \ueps g(\weps),
	\quad & x\in\Omega, \ t>0, \\[1mm]
	u_{\eps x}=w_{\eps x}=0,
	\quad & x\in\pO, \ t>0, \\[1mm]
	\ueps(x,0)=u_0(x), \qquad \weps(x,0)=w_{0\eps}(x),
	\quad & x\in\Omega,
	\end{array} \right.
\ee
which are all solvable at least locally in time, and for which a convenient criterion for extensibility can be obtained:
\begin{lem}\label{lem_loc}
  For each $\eps\in (\eps_j)_{j\in\N}$, there exist $\tme\in (0,\infty]$ and functions
  \be{l1}
	\left\{ \begin{array}{l}
	\ueps \in C^0(\overline{\Omega}\times [0,\tme)) \cap C^{2,1}(\overline{\Omega}\times (0,\tme)), \\[1mm]
	\weps \in C^0([0,\tme);W^{1,2}(\Omega)) \cap C^{2,1}(\overline{\Omega}\times (0,\tme)),
	\end{array} \right.
  \ee
  for which we have {\cred $\ueps> 0$ in $\overline{\Omega}\times (0,\tme)$} 
  and $\weps>0$ in $\overline{\Omega}\times [0,\tme)$,
  which solve (\ref{0eps}) in the classical sense in $\Omega\times (0,\tme)$, and which are such that
  \be{ext}
	\mbox{if $\tme<\infty$, then }
	\limsup_{t\nearrow \tme} \bigg\{ \|\ueps(\cdot,t)\|_{L^\infty(\Omega)} + \|\weps(\cdot,t)\|_{W^{1,2}(\Omega)} 
	+ \Big\|\frac{1}{\weps(\cdot,t)}\Big\|_{L^\infty(\Omega)} \bigg\} = \infty.
  \ee
\end{lem}
\proof
  In light of the positivity of both $\deps$ and $w_{0\eps}$ in $\overline{\Omega}$, as asserted by Lemma \ref{lem43}
  and Lemma \ref{lem45}, 
  this can be seen on adapting well-established arguments from the analysis of chemotaxis problems and of parabolic problems
  involving nonlinear degenerate diffusion
  (\cite{taowin_SIMA2011}, \cite{amann}, \cite{wiegner_NA}) to the present context.
\qed
The following two properties of these solutions are almost trivial but important.
\begin{lem}\label{lem1}
  Let $\eps\in (\eps_j)_{j\in\N}$. Then
  \be{mass}
	\io \ueps(\cdot,t)=\io u_0
	\qquad \mbox{for all } t\in (0,\tme)
  \ee
  and
  \be{M}
	\weps(x,t) \le M := \|w_0\|_{L^\infty(\Omega)}+1
	\qquad \mbox{for all $x\in\Omega$ and } t\in (0,\tme),
  \ee
  and furthermore we have
  \be{massw}
	\io \weps(\cdot,t) \le \io \weps(\cdot,t_0) \le \io w_{0\eps}
	\qquad \mbox{whenever } 0 < t_0<t<\tme
  \ee
  as well as
  \be{int_uw}
	\int_0^t \io \ueps\weps \le \frac{1}{\ug} \io w_{0\eps}
	\qquad \mbox{for all } t\in (0,\tme).
  \ee
\end{lem}
\proof
  The identity (\ref{mass}) immediately results on integration of the first equation in (\ref{0eps}) over
  $\Omega\times (0,t)$.
  For the derivation of (\ref{M}), we only need to observe that by the maximum principle,
  \bas
	\weps \le \|w_{0\eps}\|_{L^\infty(\Omega)}
	\qquad \mbox{in } \Omega\times (0,\tme),
  \eas
  and that herein by definition (\ref{45.3}) of $w_{0\eps}$, due to the fact that $\eps_j\le 1$ for all $j\in\N$ we have
  \bas
	w_{0\eps_j} = w_{0j} + \eps_j^\frac{1}{4} \le w_0 + 1 \le M
	\quad \mbox{in } \Omega
	\qquad \mbox{for all } j\in\N,
  \eas
  because $w_{0j} \le w_0$ in $\Omega$ for all $j\in\N$ by Lemma \ref{lem42}.\\
  Finally, since from the second equation in (\ref{0eps}) we obtain
  \bas
	\frac{d}{dt} \io \weps = - \io \ueps g(\weps) \le - \ug \io \ueps\weps
	\qquad \mbox{for all } t\in (0,\tme),
  \eas
  after an integration in time we readily infer that also (\ref{massw}) and (\ref{int_uw}) hold.
\qed
\mysection{An approximate energy inequality}\label{sect3}
In order to derive some fundamental a priori information beyond that from Lemma \ref{lem1}, we shall next make use of
our particular construction of the functions $\deps$ and $w_{0\eps}$ to establish an approximate version
of the energy inequality (\ref{energy}). 
This will be achieved in Lemma \ref{lem7}, and thereafter further exploited in Lemma \ref{lem8}, 
on the basis of three testing procedures performed in Lemma \ref{lem3},
Lemma \ref{lem4} and Lemma \ref{lem5}.\abs
We first consider the part containing the logarithmic entropy functional.
\begin{lem}\label{lem3}
  For all $\epsin$ and arbitrary $\delta>0$,
  \bea{3.1}
	\frac{d}{dt} \io \ueps \ln (\deps\ueps)
	+ \io \frac{(\deps\ueps)_x^2}{\deps\ueps}
	\le \io \deps \uepsx \wepsx
	+ \delta \io \deps\ueps \frac{\wepsx^2}{\weps} 
	+ \frac{1}{4\delta} \io \frac{\depsx^2}{\deps} \ueps \weps
  \eea
  for all $t\in (0,\tme)$.
\end{lem}
\proof
  We multiply the first equation in (\ref{0eps}) by the function $\ln (\deps\ueps)$ which by Lemma \ref{lem43} and
  the strong maximum principle is positive in $\overline{\Omega}\times (0,\tme)$.
  On integrating by parts and using (\ref{mass}) we thereby obtain the identity
  \bas
	\frac{d}{dt} \io \ueps \ln (\deps\ueps)
	+ \io \frac{(\deps\ueps)_x^2}{\deps\ueps}
	&=& \io (\deps\ueps)_x \wepsx \\
	&=& \io \deps \uepsx \wepsx 
	+ \io \depsx \ueps \wepsx
	\qquad \mbox{for all } t\in (0,\tme),
  \eas
 in which by Young's inequality, for each $\delta>0$ we have
  \bas
	\io \depsx \ueps \wepsx
	\le \delta \io \deps \ueps \frac{\wepsx^2}{\weps}
	+ \frac{1}{4\delta} \io \frac{\depsx^2}{\deps} \ueps\weps
	\qquad \mbox{for all } t\in (0,\tme),
  \eas
  so that (\ref{3.1}) directly follows.
\qed
As already observed in \cite{corrias_perthame_zaag2004} and essentially used
in numerous further precedent works on haptotaxis systems (see e.g.~\cite{litcanu_moralesrodrigo_M3AS}, 
\cite{taowin_JDE2014}), 
the interaction term in (\ref{3.1}) containing the gradients of both the population density and the attractant,
precisely appears during an appropriate testing process applied to the second equation in (\ref{0eps}).
Thanks to the dissipative character of the signal consumption mechanism in (\ref{0eps}), this furthermore
provides an absorptive term that can be used to compensate the second summand on the right of (\ref{3.1}).
The next lemma will moreover reveal the fortunate circumstance that 
the particular diffusive regularization chosen in the second equation in (\ref{0eps}) is in favorable accordance
with these stuctural properties.
\begin{lem}\label{lem4}
  Let $\epsin$. Then with $\ug$ and $\og$ taken from (\ref{g1}), we have
  \bea{4.1}
	\frac{1}{2} \frac{d}{dt} \io \deps \frac{\wepsx^2}{g(\weps)}
	+ \frac{\ug}{2\og} \io \deps \ueps \frac{\wepsx^2}{\weps}
	+ \eps \io \frac{1}{\sqrt{g(\weps)}} \Big(\deps\frac{\wepsx}{\sqrt{g(\weps)}}\Big)_x^2
	\le - \io \deps \uepsx \wepsx
  \eea
  for all $t\in (0,\tme)$.
\end{lem}
\proof
  Using that $\weps>0$ in $\overline{\Omega}\times [0,\tme)$ by Lemma \ref{lem_loc}, and that hence (\ref{g2}) warrants that
  also $g(\weps)$ is positive in $\overline{\Omega}\times [0,\tme)$, on the basis of the second equation in (\ref{0eps})
  and an integration by parts we compute
  \bea{4.2}
	\frac{d}{dt} \io \deps \frac{\wepsx^2}{g(\weps)}
	&=& 2\io \deps \frac{\wepsx}{g(\weps)} \cdot
	\bigg\{ \eps\Big(\deps \frac{\wepsx}{\sqrt{g(\weps)}}\Big)_{xx} - \Big(\ueps g(\weps)\Big)_x \bigg\} \nn\\
	& & - \io \deps \frac{\wepsx^2}{g^2(\weps)} g'(\weps) \cdot
	\bigg\{ \eps\Big(\deps \frac{\wepsx}{\sqrt{g(\weps)}}\Big)_x - \ueps g(\weps) \bigg\} \nn\\[1mm]
	&=& - 2\eps \io \Big(\deps \frac{\wepsx}{g(\weps)} \Big)_x \cdot \Big(\deps \frac{\wepsx}{\sqrt{g(\weps)}}\Big)_x 
		\nn\\
	& & - 2 \io \deps \uepsx \wepsx 
	- 2\io \deps \ueps g'(\weps) \frac{\wepsx^2}{g(\weps)} \nn\\
	& & - \eps \io \deps g'(\weps) \frac{\wepsx^2}{g^2(\weps)} \cdot \Big(\deps \frac{\wepsx}{\sqrt{g(\weps)}}\Big)_x 
		\nn\\
	& & + \io \deps\ueps g'(\weps) \frac{\wepsx^2}{g(\weps)} \nn\\[1mm]
	&=& -2\eps \io \frac{1}{\sqrt{g(\weps)}} \cdot \Big(\deps \frac{\wepsx}{\sqrt{g(\weps)}}\Big)_x^2 \nn\\
	& & -2 \io \deps\uepsx\wepsx \nn\\
	& & - \io \deps \ueps g'(\weps) \frac{\wepsx^2}{g(\weps)}
	\qquad \mbox{for all } t\in (0,\tme),
  \eea
  where we have used the pointwise identity
  \bas
	& & \hspace*{-20mm}
	2\Big(\deps\frac{\wepsx}{g(\weps)}\Big)_x
	+ \deps g'(\weps) \frac{\wepsx^2}{g^2(\weps)} \\
	&=& 2\Big(\frac{1}{\sqrt{g(\weps)}} \cdot \deps \frac{\wepsx}{\sqrt{g(\weps)}}\Big)_x
	+ \deps g'(\weps) \frac{\wepsx^2}{g^2(\weps)} \\
	&=& - \frac{g'(\weps)\wepsx}{\sqrt{g(\weps)}^3} \cdot \deps \frac{\wepsx}{\sqrt{g(\weps)}}
	+ 2\cdot \frac{1}{\sqrt{g(\weps)}} \cdot \Big(\deps \frac{\wepsx}{\sqrt{g(\weps)}}\Big)_x
	+ \deps g'(\weps) \frac{\wepsx^2}{g^2(\weps)} \\
	&=& 2\frac{1}{\sqrt{g(\weps)}} \cdot \Big(\deps \frac{\wepsx}{\sqrt{g(\weps)}}\Big)_x
	\qquad \mbox{in } \Omega\times (0,\tme).
  \eas
  Since (\ref{g1}) and (\ref{g2}) entail that
  \bas
	g'(\weps) \ge \ug
	\quad \mbox{and} \quad
	g(\weps) \le \og \weps
	\qquad \mbox{in } \Omega\times (0,\tme),
  \eas
  from (\ref{4.2}) we obtain (\ref{4.1}).
\qed
Finally, in order to absorb the rightmost summand in (\ref{3.1}) appropriately, we shall add a suitable multiple
of the inequality contained in the following.
\begin{lem}\label{lem5}
  Let $\epsin$ and $\delta>0$. Then for all $t\in (0,\tme)$,
  \be{5.1}
	\frac{d}{dt} \io \frac{\depsx^2}{\deps} \weps
	+ \ug\io \frac{\depsx^2}{\deps} \ueps\weps
	\le \delta\eps \io \frac{1}{\sqrt{g(\weps)}} \cdot \Big(\deps\frac{\wepsx}{\sqrt{g(\weps)}}\Big)_x^2
	+ \frac{\sqrt{\og M\eps}}{4\delta},
  \ee
  where $\ug, \og$ and $M$ are as in (\ref{g1}) and (\ref{M}), respectively.
\end{lem}
\proof
  By means of (\ref{0eps}), for $\epsin$ we calculate
  \be{5.2}
	\frac{d}{dt} \io \frac{\depsx^2}{\deps} \weps
	= \eps \io \frac{\depsx^2}{\deps} \cdot \Big(\deps \frac{\wepsx}{\sqrt{g(\weps)}}\Big)_x
	- \io \frac{\depsx^2}{\deps} \ueps g(\weps)
	\qquad \mbox{for all } t\in (0,\tme),
  \ee
  where thanks to (\ref{g2}),
  \be{5.3}
	- \io \frac{\depsx^2}{\deps} \ueps g(\weps)
	\le - \ug \io \frac{\depsx^2}{\deps} \ueps\weps
	\qquad \mbox{for all } t\in (0,\tme).
  \ee
  In order to estimate the first term on the right of (\ref{5.2}), we first invoke Young's inequality to see that
  for each $\delta>0$ we have
  \be{5.4}
	\eps\io \frac{\depsx^2}{\deps} \cdot \Big(\deps\frac{\wepsx}{\sqrt{g(\weps)}}\Big)_x
	\le \delta \eps \io \frac{1}{\sqrt{g(\weps)}} \cdot \Big(\deps\frac{\wepsx}{\sqrt{g(\weps)}}\Big)_x^2
	+ \frac{\eps}{4\delta} \io \frac{\depsx^4}{\deps^2} \sqrt{g(\weps)}
	\qquad \mbox{for all } t\in (0,\tme),
  \ee
  and here in the rightmost summand we recall (\ref{g2}) and (\ref{M}) to find that
  \bas
	\sqrt{g(\weps)} \le \sqrt{\og M}
	\qquad \mbox{in } \Omega\times (0,\tme).
  \eas
  Since in Lemma \ref{lem43} we have asserted that
  \bas
	\io \frac{\depsx^4}{\deps^2} \le \frac{1}{\sqrt{\eps}}
	\qquad \mbox{for all } \epsin,
  \eas
  this entails that 
  \bas
	\frac{\eps}{4\delta} \io \frac{\depsx^4}{\deps^2} \sqrt{g(\weps)}
	\le \frac{\sqrt{\og M\eps}}{4\delta}
	\qquad \mbox{for all } t\in (0,\tme),
  \eas
  so that combining (\ref{5.3}) and (\ref{5.4}) with (\ref{5.2}) yields (\ref{5.1}).
\qed
In summary, on adequately joining the above three lemmata we obtain the desired approximate analogue of the energy
inequality (\ref{energy}).
\begin{lem}\label{lem7}
  Let $\ug,\og$ and $M$ denote the constants from (\ref{g1}) and (\ref{M}).
  Then whenever $\epsin$,
  \bea{7.1}
	& & \hspace*{-30mm}
	\frac{d}{dt} \bigg\{
	\io \ueps \ln (\deps\ueps)
	+ \frac{1}{2} \io \deps \frac{\wepsx^2}{g(\weps)}
	+ \frac{\og}{\ug^2} \io \frac{\depsx^2}{\deps} \weps \bigg\} \nn\\
	& & + \io \frac{(\deps\ueps)_x^2}{\deps\ueps}
	+ \frac{\ug}{4\og} \io \deps\ueps \frac{\wepsx^2}{\weps} 
	+ \frac{\eps}{2} \io \frac{1}{\sqrt{g(\weps)}} \cdot \Big(\deps\frac{\wepsx}{\sqrt{g(\weps)}}\Big)_x^2 \nn\\
	&\le& \frac{\sqrt{\og^5 M\eps}}{2\ug^4}
	\qquad \mbox{for all } t\in (0,\tme).
  \eea
\end{lem}
\proof
  We choose the free parameters $\delta$ in Lemma \ref{lem3} and Lemma \ref{lem5} to equal
  $\frac{\ug}{4\og}$ and $\frac{\ug^2}{2\og}$, respectively, to see on linearly combining (\ref{3.1}), (\ref{4.1}) and
  (\ref{5.1}) that for all $t\in (0,\tme)$,
  \bas
	& & \hspace*{-40mm}
	\frac{d}{dt} \bigg\{
	\io \ueps \ln (\deps\ueps)
	+ \frac{1}{2} \io \deps \frac{\wepsx^2}{g(\weps)}
	+ \frac{\og}{\ug^2} \io \frac{\depsx^2}{\deps} \weps \bigg\}
	+ \io \frac{(\deps\ueps)_x^2}{\deps\ueps}
	+ \frac{\ug}{2\og} \io \deps\ueps\frac{\wepsx^2}{\weps} \nn\\
	& & + \eps \io \frac{1}{\sqrt{g(\weps)}} \cdot \Big(\deps\frac{\wepsx}{\sqrt{g(\weps)}}\Big)_x^2
	+ \frac{\og}{\ug^2} \cdot \ug \io \frac{\depsx^2}{\deps} \ueps\weps \\[2mm]
	&\le& \io \deps\uepsx\wepsx
	+ \frac{\ug}{4\og} \io \deps\ueps \frac{\wepsx^2}{\weps}
	+ \frac{\og}{\ug} \io \frac{\depsx^2}{\deps} \ueps\weps \\
	& & - \io \deps\uepsx\wepsx \\
	& & + \frac{\og}{\ug^2} \cdot \frac{\ug^2}{2\og} \cdot \eps 
	\io \frac{1}{\sqrt{g(\weps)}}\cdot \Big(\deps \frac{\wepsx}{\sqrt{g(\weps)}}\Big)_x^2
	+ \frac{\og}{\ug^2} \cdot \frac{\sqrt{\og M \eps}}{4\frac{\ug^2}{2\og}},
  \eas
  which can readily be simplified so as to yield (\ref{7.1}).
\qed
We now use Lemma \ref{lem45} to make sure that the respective energy values at the initial time are bounded from above
uniformly with respect to $\eps\in (\eps_j)_{j\in\N}$.
Therefore, an integration of (\ref{7.1}) yields the following.
\begin{lem}\label{lem8}
  There exists $C>0$ with the property that whenever $\epsin$, we have
  \be{8.1}
	\io \ueps(\cdot,t) \Big|\ln (\deps\ueps(\cdot,t))\Big|
	\le C\cdot (1+\sqrt{\eps} t)
	\qquad \mbox{for all } t\in (0,\tme)
  \ee
  and
  \be{8.2}
	\io \deps \frac{\wepsx^2(\cdot,t)}{\weps(\cdot,t)}
	\le C\cdot (1+\sqrt{\eps} t)
	\qquad \mbox{for all } t\in (0,\tme)
  \ee
  as well as
  \be{8.3}
	\int_0^t \io \frac{(\deps\ueps)_x^2}{\deps\ueps} 
	\le C\cdot (1+\sqrt{\eps} t)
	\qquad \mbox{for all } t\in (0,\tme)
  \ee
  and
  \be{8.33}
	\int_0^t \io \deps\ueps \frac{\wepsx^2}{\weps}
	\le C\cdot (1+\sqrt{\eps} t)
	\qquad \mbox{for all } t\in (0,\tme)
  \ee
  and
  \be{8.34}
	\eps \int_0^t \io \frac{1}{\sqrt{g(\weps)}} \cdot \Big(\deps\frac{\wepsx}{\sqrt{g(\weps)}}\Big)_x^2
	\le C\cdot (1+\sqrt{\eps} t)
	\qquad \mbox{for all } t\in (0,\tme).
  \ee
\end{lem}
\proof
  For $\epsin$ we obtain from Lemma \ref{lem7} that if we let $\ug>0, \og>0$ and $M>0$ be as specified in (\ref{g1})
  and (\ref{M}), then
  \bas
	\yeps(t):=\io \ueps(\cdot,t) \ln (\deps\ueps(\cdot,t))
	+ \frac{1}{2} \io \deps \frac{\wepsx^2(\cdot,t)}{g(\weps(\cdot,t))}
	+ \frac{\og}{\ug^2} \io \frac{\depsx^2}{\deps} \weps(\cdot,t),
	\qquad t\in [0,\tme),
  \eas
  and
  \bas
	\heps(t)
	&:=& \io \frac{(\deps\ueps(\cdot,t))_x^2}{\deps\ueps(\cdot,t)}
	+ \frac{\ug}{4\og} \io \deps \ueps(\cdot,t)\frac{\wepsx^2(\cdot,t)}{\weps(\cdot,t)} \\
	& & + \frac{\eps}{2} \io \frac{1}{\sqrt{g(\weps(\cdot,t))}} \cdot 
	\Big(\deps \frac{\wepsx(\cdot,t)}{\sqrt{g(\weps(\cdot,t))}} \Big)_x^2,
	\qquad t\in (0,\tme),
  \eas
  satisfy
  \be{8.5}
	\yeps'(t) + \heps(t)
	\le c_1\sqrt{\eps}
	\qquad \mbox{for all } t\in (0,\tme),
  \ee
  where $c_1:=\frac{\sqrt{\og^5 M}}{2\ug^4}$.
  To conclude (\ref{8.1})-(\ref{8.33}) from this, we observe that at the initial time we can use (\ref{43.3}) to estimate
  \be{8.6}
	\io u_0\ln (\deps u_0)
	\le c_2:=\|u_0\|_{L^\infty(\Omega)} \ln \Big\{ (\|d\|_{L^\infty(\Omega)}+1) \|u_0\|_{L^\infty(\Omega)} \Big\}
	\qquad \mbox{for all } \epsin,
  \ee
  whereas Lemma \ref{lem45} ensures the existence of $c_3>0$ and $c_4>0$ such that
  \be{8.7}
	\io \deps \frac{w_{0\eps x}^2}{w_{0\eps}} \le c_3
	\qquad \mbox{for all } \epsin
  \ee
  and
  \be{8.8}
	\io \frac{\depsx^2}{\deps} w_{0\eps} \le c_4
	\qquad \mbox{for all } \epsin.
  \ee
  Since $\frac{w_{0\eps x}^2}{g(w_{0\eps})} \le \frac{1}{\ug} \frac{w_{0\eps x}^2}{w_{0\eps}}$ in $\Omega$ by (\ref{g2}),
  (\ref{8.6})-(\ref{8.8}) show that
  \bas
	\yeps(0) \le c_5:=c_2 + \frac{c_3}{2\ug} + \frac{c_4\og}{\ug^2}
	\qquad \mbox{for all } \epsin,
  \eas
  by (\ref{8.5}) implying that for all $\epsin$,	
  \be{8.10}
	\yeps(t) + \int_0^t \heps(s)ds
	\le c_5 + c_1\sqrt{\eps} t
	\qquad \mbox{for all } t\in (0,\tme).
  \ee
  Here we note that since $\xi\ln \xi \ge -\frac{1}{e}$ for all $\xi>0$ and $\deps\ge d$ in $\Omega$ for all
  $\epsin$ thanks to Lemma \ref{lem43}, 
{\cred with $c_6:=\frac{1}{e} \io \frac{1}{d}$ being finite according to (\ref{d2}), we have}
  \bas
	- \int_{\{\deps\ueps<1\}} \ueps \ln (\deps\ueps)
	&=& - \int_{\{\deps\ueps<1\}} \frac{1}{\deps} \cdot \Big\{ \deps\ueps \cdot \ln (\deps\ueps)\Big\} \\
	&\le& \frac{1}{e} \int_{\{\deps\ueps<1\}} \frac{1}{\deps} \\
	&\le& {\cred c_6}		
	\qquad \mbox{for all $\epsin$ and } t\in (0,\tme),
  \eas
  which along with (\ref{g2}) in particular entails that
  \bas
	\yeps(t)
	&\ge& \io \ueps\ln (\deps\ueps) + \frac{1}{2\og} \io \deps \frac{\wepsx^2}{\weps} \\
	&=& \io \ueps \Big|\ln (\deps\ueps)\Big|
	+ 2\int_{\{\deps\ueps<1\}} \ueps \ln (\deps\ueps)
	+ \frac{1}{2\og} \io \deps\frac{\wepsx^2}{\weps} \\
	&\ge& \io \ueps \Big|\ln (\deps\ueps)\Big|
	-2c_6
	+ \frac{1}{2\og} \io \deps\frac{\wepsx^2}{\weps}
	\qquad \mbox{for all $\epsin$ and } t\in (0,\tme).
  \eas
  Therefore, (\ref{8.10}) implies that for all $\epsin$ we have
  \bas
	\io \ueps\Big|\ln (\deps\ueps)\Big|	
	+ \frac{1}{2\og} \io \deps \frac{\wepsx^2}{\weps} 
	+ \int_0^t \heps(s)ds
	\le 2c_6+c_5+c_1\sqrt{\eps} t
	\qquad \mbox{for all } t\in (0,\tme),
  \eas
  which in view of the definition of $\heps$ yields all claimed inequalities.
\qed
\mysection{Global existence in the approximate problems}\label{sect4}
With the above information at hand, we can now make sure that in fact all our approximate solutions are global in time.
To achieve this in Lemma \ref{lem93}
on the basis of the extensibility criterion in Lemma \ref{lem_loc}, for each individual $\epsin$
we will derive further estimates which may depend on $\eps$.
We begin with a pointwise lower estimate for $\weps$ that we obtain by a comparison argument combined with
Lemma \ref{lem8}, and that will be used in Lemma \ref{lem90}.
\begin{lem}\label{lem900}
  Assume that $\tme<\infty$ for some $\epsin$. Then there exists $C(\eps)>0$ such that
  \be{900.1}
	\weps(x,t) \ge C(\eps)
	\qquad \mbox{for all $x\in\Omega$ and } t\in (0,\tme).
  \ee
\end{lem}
\proof
  We first observe that under the current hypothesis, Lemma \ref{lem8} says that 
  \be{900.2}
{\cred	\int_0^{\tme} \io \Big| \Big(\sqrt{\deps\ueps}\Big)_x\Big|^2 <\infty,}
  \ee
  and we claim that along with (\ref{mass}) this provides sufficient regularity information on
  the absorption coefficient function 
  $\ueps$ in the second equation in (\ref{0eps}) to rule out finite-time formation of zeros of $\weps$
  in the sense of (\ref{900.1}).
  To verify this, we use the continuity of the embedding 
  $W^{1,2}(\Omega)\hra L^\infty(\Omega)$ as well as (\ref{mass}) and (\ref{43.3}) to fix $c_1>0$ and $c_2>0$
  fulfilling
  \bas
	\int_0^{\tme} \Big\|\sqrt{\deps\ueps(\cdot,t)}\Big\|_{L^\infty(\Omega)}^2 dt
	&\le& c_1 \int_0^{\tme} \bigg\{ \Big\|\Big(\sqrt{\deps\ueps(\cdot,t)}\Big)_x\Big\|_{L^2(\Omega)}^2 
	+ \Big\|\sqrt{\deps\ueps(\cdot,t)}\Big\|_{L^2(\Omega)}^2 \bigg\} dt \\
	&\le& c_2 \int_0^{\tme} \bigg\{ \Big\|\Big(\sqrt{\deps\ueps(\cdot,t)}\Big)_x\Big\|_{L^2(\Omega)}^2 +1 \bigg\} dt,
  \eas
  and to see that thus (\ref{900.2}) entails that 
  \bas
	\int_0^{\tme} \Big\|\sqrt{\deps\ueps(\cdot,t)}\Big\|_{L^\infty(\Omega)}^2 dt < \infty.
  \eas
  hence, by positivity of $\deps$ in $\bom$, also the number
  \bas
	c_3:=\int_0^{\tme} \|\ueps(\cdot,t)\|_{L^\infty(\Omega)} dt
  \eas
  is finite, which 
  in particular implies that the solution
  $y\in C^1([0,\tme))$ of the initial-value problem
  \be{900.9}
	\left\{ \begin{array}{l}
	y'(t) = - \og \|\ueps(\cdot,t)\|_{L^\infty(\Omega)} \cdot y(t), \qquad t\in (0,\tme), \\[1mm]
	y(0)=\eps^\frac{1}{4},
  	\end{array} \right.
  \ee
  satisfies
  \be{900.10}
	y(t) = \eps^\frac{1}{4} e^{-\og \int_0^t \|\ueps(\cdot,s)\|_{L^\infty(\Omega)} ds}
	\ge c_4:=\eps^\frac{1}{4} e^{-\og c_3}
	\qquad \mbox{for all } t\in (0,\tme).
  \ee
{\cred It can now readily be verified that $\overline{\Omega}\times [0,\tme) \ni (x,t) \mapsto y(t)$ 
  is a classical subsolution to the 
  initial-boundary problem solved by $\weps$ in $\Omega\times (0,\tme)$, so that by (\ref{900.10}), 
  $\weps(x,t) \ge c_4$ for all $x\in\Omega$ and $t\in (0,\tme)$, which yields (\ref{900.1}).
}
\qed
This lower bound enables us to suitably estimate singular denominators appearing in the following lemma which, 
apart from that and the positivity of $\deps$, again only relies on Lemma \ref{lem8} only.
\begin{lem}\label{lem90}
  Let $\epsin$, and suppose that $\tme<\infty$. Then there exists $C(\eps)>0$ such that
  \be{90.01}
	\io \wepsx^4(\cdot,t) \le C(\eps)
	\qquad \mbox{for all } t\in (0,\tme).
  \ee
\end{lem}
\proof
  According to Lemma \ref{lem8} and Lemma \ref{lem900}, our hypothesis that $\tme<\infty$ again warrants that
  \be{90.1}
	\int_0^{\tme} \io \frac{(\deps\ueps)_x^2}{\deps\ueps} <\infty,
  \ee
  and that moreover with some $c_1>0$ and $c_2>0$ we have
  \be{90.2}
	\io \deps \frac{\wepsx^2}{\weps} \le c_1
	\qquad \mbox{for all } t\in (0,\tme)
  \ee
  as well as
  \be{90.11}
	\weps\ge c_2
	\qquad \mbox{in } \Omega\times (0,\tme).
  \ee
  Since combining the Gagliardo-Nirenberg inequality with 
  (\ref{mass}) yields positive constants $c_3$ and $c_4$ such that
  \bas
	\int_0^{\tme} \io (\deps\ueps)^3
	&=& \int_0^{\tme} \Big\|\sqrt{\deps\ueps(\cdot,t)}\Big\|_{L^6(\Omega)}^6 dt \\
	&\le& c_3 \int_0^{\tme} \bigg\{ \Big\|\Big(\sqrt{\deps\ueps(\cdot,t)}\Big)_x\Big\|_{L^2(\Omega)}^2 
	\Big\|\sqrt{\deps\ueps(\cdot,t)}\Big\|_{L^2(\Omega)}^4 
	+ \Big\|\sqrt{\deps\ueps(\cdot,t)}\Big\|_{L^2(\Omega)}^6 \bigg\} dt \\
	&\le& c_4 \int_0^{\tme} \bigg\{ \Big\|\Big(\sqrt{\deps\ueps(\cdot,t)}\Big)_x\Big\|_{L^2(\Omega)}^2 +1 \bigg\} dt,
  \eas
  it follows from (\ref{90.1}) that also
  \bas
	\int_0^{\tme} \io (\deps\ueps)^3 < \infty.
  \eas
  As
  \be{90.6}
	c_5:=\inf_{x\in\Omega} \deps(x)
  \ee
  is positive thanks to Lemma \ref{lem43}, this implies that
  \be{90.8}
	c_6:=\int_0^{\tme} \io \ueps^3
  \ee
  is finite, whereas (\ref{90.2}) and (\ref{M}) show that with some $c_7>0$ we have
  \be{90.66}
	\io \wepsx^2 \le c_7
	\qquad \mbox{for all } t\in (0,\tme).
  \ee
  We now use the second equation in (\ref{0eps}) to compute
  \bea{90.12}
	\hspace*{-5mm}
	\frac{1}{4} \frac{d}{dt} \io \wepsx^4
	&=& -3\io \wepsx^2 \wepsxx w_{\eps t} \nn\\
	&=& - 3\eps \io \wepsx^2 \wepsxx \cdot \Big(\deps \frac{\wepsx}{\sqrt{g(\weps)}}\Big)_x
	+ 3\io \io \ueps g(\weps) \wepsx^2 \wepsxx \nn\\
	&=& -3\eps \io \deps \frac{1}{\sqrt{g(\weps)}} \wepsx^2 \wepsxx^2
	+ \frac{3}{2} \eps \io \deps \frac{g'(\weps)}{\sqrt{g(\weps)}^3} \wepsx^4 \wepsxx \nn\\
	& & -3\eps \io \depsx \frac{1}{\sqrt{g(\weps)}} \wepsx^3 \wepsxx
	+ 3\io \ueps g(\weps) \wepsx^2 \wepsxx
	\qquad \mbox{for all } t\in (0,\tme),
  \eea
  where by Young's inequality, (\ref{90.11}), (\ref{90.6}), (\ref{g1}), (\ref{g2}) and (\ref{M}),
  for all $t\in (0,\tme)$ we can estimate
  \bea{90.13}
	\frac{3}{2} \eps \io \deps \frac{g'(\weps)}{\sqrt{g(\weps)}^3} \wepsx^4 \wepsxx
	&\le& \eps \io \deps \frac{1}{\sqrt{g(\weps)}} \wepsx^2 \wepsxx^2
	+ \frac{9}{16}\eps \io \deps \frac{g'^2(\weps)}{\sqrt{g(\weps)}^5} \wepsx^6 \nn\\
	&\le& \eps \io \deps \frac{1}{\sqrt{g(\weps)}} \wepsx^2 \wepsxx^2
	+ c_8 \io \wepsx^6
  \eea	
  and 
  \bea{90.14}
	-3\eps \io \depsx \frac{1}{\sqrt{g(\weps)}} \wepsx^3 \wepsxx
	&\le& \eps \io \deps \frac{1}{\sqrt{g(\weps)}} \wepsx^2 \wepsxx^2 
	+ \frac{9}{4}\eps \io \frac{\depsx^2}{\deps} \frac{1}{\sqrt{g(\weps)}} \wepsx^4 \nn\\
	&\le& \eps \io \deps \frac{1}{\sqrt{g(\weps)}} \wepsx^2 \wepsxx^2 
	+ c_9 \io \wepsx^4 \nn\\
	&\le& \eps \io \deps \frac{1}{\sqrt{g(\weps)}} \wepsx^2 \wepsxx^2 
	+ c_9 \io \wepsx^6 + c_9 |\Omega|
  \eea
  as well as
  \bea{90.15}
	3\io \ueps g(\weps) \wepsx^2 \wepsxx
	&\le& \frac{\eps}{2} \io \deps \frac{1}{\sqrt{g(\weps)}} \wepsx^2 \wepsxx^2 
{\cred 	+ \frac{9}{2\eps} \io \frac{1}{\deps} \ueps^2 \sqrt{g(\weps)}^5 \wepsx^2 
} \nn\\
	&\le& \frac{\eps}{2} \io \deps \frac{1}{\sqrt{g(\weps)}} \wepsx^2 \wepsxx^2 
	+ c_{10} \io \ueps^2 \wepsx^2 \nn\\
	&\le& \frac{\eps}{2} \io \deps \frac{1}{\sqrt{g(\weps)}} \wepsx^2 \wepsxx^2 
	+ c_{10} \io \ueps^3 + c_{10} \io \wepsx^6
  \eea
  with $c_8:=\frac{9\eps \og^2 \|\deps\|_{L^\infty(\Omega)}}{16\sqrt{\ug c_2}^5}$,
  $c_9:=\frac{9\eps\|\depsx\|_{L^\infty(\Omega)}^2}{4c_5\sqrt{\ug c_2}}$
  and {\cred $c_{10}:=\frac{9\sqrt{\og M}^5}{2\eps c_5}$.}
  Since (\ref{90.6}), (\ref{g2}) and (\ref{M}) moreover entail that writing $c_{11}:=\frac{\eps c_5}{2\sqrt{\og M}}$ we obtain
  \bas
	\frac{\eps}{2} \io \deps \frac{1}{\sqrt{g(\weps)}} \wepsx^2 \wepsxx^2
	\ge c_{11} \io \wepsx^2 \wepsxx^2
	= \frac{c_{11}}{4} \io (\wepsx^2)_x^2
	\qquad \mbox{for all } t\in (0,\tme),
  \eas
  and since the Gagliardo-Nirenberg inequality and Young's inequality together with (\ref{90.66})
  show that with some positive constants $c_{12}, c_{13}$ and $c_{14}$ we have
  \bas
	(c_8+c_9+c_{10}) \io \wepsx^6
	&=& (c_8+c_9+c_{10}) \|\wepsx^2\|_{L^3(\Omega)}^3 \\
	&\le& c_{12} \Big\|(\wepsx^2)_x\Big\|_{L^2(\Omega)}^\frac{4}{3} \|\wepsx^2\|_{L^1(\Omega)}^\frac{5}{3}
	+ c_{12} \|\wepsx^2\|_{L^1(\Omega)}^3 \\
	&\le& c_{13} \Big\|(\wepsx^2)_x\Big\|_{L^2(\Omega)}^\frac{4}{3} + c_{13} \\
	&\le& \frac{c_{11}}{4} \io (\wepsx^2)_x^2
	+ c_{14}
	\qquad \mbox{for all } t\in (0,\tme),
  \eas
  on combining (\ref{90.12})-(\ref{90.15}) we therefore see that
  \bas
	\frac{1}{4} \frac{d}{dt} \io \wepsx^4
	\le c_9|\Omega| + c_{10} \io \ueps^3 + c_{14}
	\qquad \mbox{for all } t\in (0,\tme),
  \eas
  and that hence
  \bas
	\io \wepsx^4(\cdot,t)
	\le \io w_{0\eps x}^4 + 4c_{10} c_6 + 4(c_9|\Omega|+c_{14}) \cdot t
	\qquad \mbox{for all } t\in (0,\tme)
  \eas
  because of (\ref{90.8}). Again since $\tme$ was assumed to be finite, this entails (\ref{90.01}).
\qed
The above regularity information on the haptotactic gradient is now sufficient to warrant an $\eps$-dependent bound
for $\ueps$ in $L^p(\Omega)$ for arbitrarily large $p$.
\begin{lem}\label{lem91}
  Assume that $\tme<\infty$ for some $\epsin$. Then for all $p\in (1,\infty)$ there exists $C(\eps,p)>0$ such that
  \be{91.1}
	\io \ueps^p(\cdot,t) \le C(\eps,p)
	\qquad \mbox{for all } t\in (0,\tme).
  \ee
\end{lem}
\proof
  We test the first equation in (\ref{0eps}) against $\ueps^{p-1}$ and use Young's inequality to see that
  \bas
	\frac{1}{p} \frac{d}{dt} \io \ueps^p
	&=& \io \ueps^{p-1} \cdot \Big\{ (\deps\ueps)_x - \deps\ueps \wepsx \Big\}_x \\
	&=& - (p-1) \io \deps \ueps^{p-2} \uepsx^2
	-(p-1) \io \depsx \ueps^{p-1} \uepsx
	+ (p-1) \io \deps \ueps^{p-1} \uepsx \wepsx \\
	&\le& - \frac{p-1}{2} \io \deps \ueps^{p-2} \uepsx^2
	+ (p-1) \io \frac{\depsx^2}{\deps} \ueps^p
	+ (p-1) \io \deps \ueps^p \wepsx^2
	\qquad \mbox{for all } t\in (0,\tme),
  \eas
  so that since $\deps$ is smooth and positive throughout $\overline{\Omega}$, we can find $c_1>0$ and $c_2>0$ fulfilling
  \be{91.2}
	\frac{d}{dt} \io \ueps^p + c_1 \io (\ueps^\frac{p}{2})_x^2
	\le c_2 \io \ueps^p + c_2 \io \ueps^p \wepsx^2
	\qquad \mbox{for all } t\in (0,\tme).
  \ee
  Here using the Cauchy-Schwarz inequality, the Gagliardo-Nirenberg inequality and Young's inequality along with 
  (\ref{mass}) and the estimate from Lemma \ref{lem90}, we obtain positive constants $c_3, c_4, c_5$ and $c_6$ such that
  \bas
	c_2 \io \ueps^p \wepsx^2
	&\le& c_2 \cdot \bigg\{ \io \ueps^{2p} \bigg\}^\frac{1}{2} \cdot
	\bigg\{ \io \wepsx^4 \bigg\}^\frac{1}{2} \\
	&\le& c_3 \|\ueps^\frac{p}{2}\|_{L^4(\Omega)}^2 \\
	&\le& c_4 \|(\ueps^\frac{p}{2})_x\|_{L^2(\Omega)}^\frac{2p-1}{p+1}
	\|\ueps^\frac{p}{2}\|_{L^\frac{2}{p}(\Omega)}^\frac{3}{p+1}
	+ c_4\|\ueps^\frac{p}{2}\|_{L^\frac{2}{p}(\Omega)}^2 \\
	&\le& c_5 \| (\ueps^\frac{p}{2})_x\|_{L^2(\Omega)}^\frac{2p-1}{p+1} + c_5 \\
	&\le& c_1 \io (\ueps^\frac{p}{2})_x^2 + c_6
	\qquad \mbox{for all } t\in (0,\tme),
  \eas
  because $\frac{2p-1}{p+1}<2$. Therefore, (\ref{91.2}) entails that
  \bas
	\frac{d}{dt} \io \ueps^p \le c_2 \io \ueps^p + c_6
	\qquad \mbox{for all } t\in (0,\tme)
  \eas
  and thus, upon integration, that
  \bas
	\io \ueps^p \le \bigg\{ \io u_0^p \bigg\} \cdot e^{c_2 t} + \frac{c_6}{c_2} (e^{c_2 t}-1)
	\qquad \mbox{for all } t\in (0,\tme),
  \eas
  which implies (\ref{91.1}).
\qed
By means of a standard result based on a Moser-type iteration, along with Lemma \ref{lem90} this readily
yields boundedness of $\ueps$ whenever $\tme<\infty$.
\begin{lem}\label{lem92}
  If $\tme<\infty$ for some $\epsin$, then there exists $C(\eps)>0$ such that
  \be{92.1}
	\|\ueps(\cdot,t)\|_{L^\infty(\Omega)} \le C(\eps)
	\qquad \mbox{for all } t\in (0,\tme).
  \ee
\end{lem}
\proof
  We rewrite the first equation in (\ref{0eps}) in the form
  \bas
	u_{\eps t} = (\deps \uepsx)_x + h_{\eps x},
	\qquad x\in\Omega, \ t\in (0,\tme),
  \eas
  with
  \bas
	\heps(x,t):=\depsx(x)\ueps(x,t) - \deps(x)\ueps(x,t) \wepsx(x,t),
	\qquad x\in\Omega, \ t\in (0,\tme),
  \eas
  and note that for any fixed $q\in (3,4)$, by the H\"older inequality we have	
  \bas
	\io |\heps(\cdot,t)|^q
	&\le& 2^{q-1} \io |\deps|^q \ueps^q 
	+ 2^{q-1} \io \deps^q \ueps^q |\wepsx|^q \\
	&\le& 2^{q-1} |\depsx\|_{L^\infty(\Omega)}^q \io \ueps^q
	+ 2^{q-1}\|\deps\|_{L^\infty(\Omega)}^q \cdot
	\bigg\{ \io \ueps^\frac{4q}{4-q}\bigg\}^\frac{4-q}{4} \cdot \bigg\{ \io \wepsx^4 \bigg\}^\frac{q}{4} 
  \eas
  for all $t\in (0,\tme)$.
  As Lemma \ref{lem90} and Lemma \ref{lem91} guarantee that
  \be{92.3}
	\sup_{t\in (0,\tme)} \io \wepsx^4(\cdot,t)<\infty
	\qquad \mbox{and} \qquad
	\sup_{t\in (0,\tme)} \io \ueps^p(\cdot,t)<\infty
	\quad \mbox{for all } p\in (1,\infty),
  \ee
  this implies that $\heps$ belongs to $L^\infty((0,\tme);L^q(\Omega))$, 
  so that using that $q>3$ we may 
  apply a known Moser-type result on boundedness in scalar parabolic 
  equations (\cite[Lemma A.1]{taowin_JDE_subcrit}) to see that 	
  along with the identities
  \bas
	\uepsx=0
	\quad \mbox{and} \quad
	\heps=0
	\qquad \mbox{on } \pO\times (0,\tme),
  \eas
  the latter being asserted by the fact that $\depsx=0$ on $\pO$ by (\ref{43.2}),
  the second property in (\ref{92.3}) is sufficient to warrant (\ref{92.1}).
\qed
In conclusion, finite-time blow-up cannot occur in any of the approximate problems.
\begin{lem}\label{lem93}
  For all $\epsin$, the solution of (\ref{0eps}) is global in time.
\end{lem}
\proof
  In view of the extensibility criterion (\ref{ext}), we only need to collect the outcomes of Lemma \ref{lem900},
  Lemma \ref{lem90} and Lemma \ref{lem92}.
\qed
\mysection{Further $\eps$-independent regularity properties of (\ref{0eps})}
\subsection{Equi-integrability properties}\label{sect_equi}
Now a key to both our existence proof and our stabilization results consists in the observation that due to
Lemma \ref{lem8}, and again due to the assumed
integrability of $\frac{1}{d}$, the solution component $\ueps$ enjoys {\cred a certain} doubly uniform integrability
property.
{\cred In order to prepare this and also our subsequent analysis, let us introduce
\be{omega}
	\od(\delta):=\sup \bigg\{ \int_E \frac{1}{d} \ \bigg| \ E\subset \Omega
	\mbox{ is measurable with } |E| \le \delta \bigg\}
	\qquad \mbox{for } \delta>0,
\ee
and observe that then our integrability assumption (\ref{d2}) warrants that
\be{omega1}
	\od(\delta)\to 0
	\qquad \mbox{as } \delta\searrow 0.
\ee
Along with Lemma \ref{lem8}, this will entail the following.
}
\begin{lem}\label{lem826}
  For all $\eta>0$ there exists $\delta>0$ such that whenever $\epsin$,
  \be{826.1}
	\int_E \ueps(\cdot,t) \le \eta
	\qquad \mbox{for all } t\in \Big(0,\frac{1}{\sqrt{\eps}}\Big)
	\quad \mbox{and any measurable $E\subset\Omega$ such that } |E|\le\delta.
  \ee
\end{lem}
\proof
  According to Lemma \ref{lem8} we can fix $c_1>0$ such that for all $\epsin$ we have
  \bas
	\io \ueps \Big| \ln (\deps\ueps)\Big| \le c_1 \cdot (1+\sqrt{\eps}t)
	\qquad \mbox{for all } t>0,,
  \eas
  whence in particular
  \be{826.2}
	\io \ueps \Big| \ln (\deps\ueps) \Big| \le 2c_1
	\qquad \mbox{for all } t\in \Big(0,\frac{1}{\sqrt{\eps}}\Big).
  \ee
  We then let $\eta>0$ be given and pick $L>1$ large enough fulfilling
  {\cred $\frac{2c_1}{\ln L} \le \frac{\eta}{2}$,
  and thereafter make use of (\ref{omega1}) in choosing some $\delta>0$ such that
  $\od(\delta) \le \frac{\eta}{2L}$.
  Then for any measurable $E\subset\Omega$ with $|E|\le\delta$
  and each $\epsin$, again using 
  the fact that
} 
  $\deps\ge d$ we find that
  \bas
	\int_E \ueps
	&=& \int_{E\cap\{\deps\ueps>L\}} \ueps
	+ \int_{E\cap\{\deps\ueps\le L\}} \ueps \\
	&\le& \int_{E\cap\{\deps\ueps>L\}} \ueps \cdot \frac{\ln (\deps\ueps)}{\ln L}
	+ \int_{E\cap\{\deps\ueps\le L\}} \frac{L}{\deps} \\
	&\le& \frac{1}{\ln L} \cdot \int_{E\cap\{\deps\ueps>L\}} \ueps \Big| \ln (\deps\ueps)\Big|
	+ L \int_{E\cap\{\deps\ueps\le L\}} \frac{1}{\deps} \\
	&\le& \frac{1}{\ln L} \cdot \io \ueps \Big|\ln (\deps\ueps)\Big| + L\int_E \frac{1}{\deps} \\
	&\le& \frac{1}{\ln L} \cdot 2c_1 
	+ L {\cred \od(\delta)} \\
	&\le& \frac{\eta}{2}+\frac{\eta}{2}=\eta
	\qquad \mbox{for all } t\in\Big(0,\frac{1}{\sqrt{\eps}}\Big),
  \eas
  as claimed.
\qed
{\cred Likewise, the weighted $L^2$ estimate for $\wepsx$ in Lemma \ref{lem8} can be turned into a corresponding 
equi-integrability statement for $\wepsx$, 
and apart from that it implies an additional boundedness property of $\weps$ in a space compactly 
embedded into $C^0(\overline{\Omega})$.}
\begin{lem}\label{lem824}
  For all $\eta>0$ there exists $\delta>0$ with the property that for arbitrary $\epsin$,
  \be{824.1}
	\int_E |\wepsx(\cdot,t)| \le \eta
	\qquad \mbox{for all } t\in \Big(0,\frac{1}{\sqrt{\eps}}\Big)
	\quad \mbox{whenever $E\subset\Omega$ is measurable with } |E|\le\delta.
  \ee
{\cred
  Moreover, there exists $C>0$ such that for arbitrary $\epsin$,
  \be{824.11}
	\|\weps(\cdot,t)\|_Y \le C
	\qquad \mbox{for all } t\in \Big(0,\frac{1}{\sqrt{\eps}}\Big),
  \ee	
  where the Banach space $Y$ is defined by
  \be{Y}
	Y:=\Bigg\{ \varphi\in C^0(\bom) \ \Bigg| \ 
	\|\varphi\|_Y:=\|\varphi\|_{L^\infty(\Omega)}
	+ \sup_{x,y\in\Omega, \ x\ne y} \frac{|\varphi(x)-\varphi(y)|}{\sqrt{\od(|x-y|)}} <\infty \Bigg\}
  \ee
  with $\od$ as in (\ref{omega}).
}
\end{lem}
\proof
  From Lemma \ref{lem8} and (\ref{M}) we obtain $c_1>0$ such that for all $\epsin$,
  \bas
	\io \deps \wepsx^2 \le c_1\cdot (1+\sqrt{\eps} t)
	\qquad \mbox{for all } t>0
  \eas
  and hence
  \be{824.2}
	\io \deps \wepsx^2 \le 2c_1
	\qquad \mbox{for all } t\in\Big(0,\frac{1}{\sqrt{\eps}}\Big).
  \ee
{\cred
  Therefore, an application of the Cauchy-Schwarz inequality shows that for arbitrary measurable $F\subset\Omega$
  we can estimate
  \bea{824.3}
	\int_F |\wepsx| 
	\le \bigg\{ \int_F\deps \wepsx^2 \bigg\}^\frac{1}{2} \cdot \bigg\{ \int_F \frac{1}{\deps} \bigg\}^\frac{1}{2}
	\le \sqrt{2c_1} \sqrt{\od(|F|)}
	\qquad \mbox{for all } t\in\Big(0,\frac{1}{\sqrt{\eps}}\Big).
  \eea
  In particular, if given $\eta>0$ we let $\delta>0$ be such that $\od(\delta)\le \frac{\eta^2}{2c_1}$, then for each
  measurable $E\subset\Omega$ fulfilling $|E|\le \delta$ we conclude from (\ref{824.3}) that
  \bas
	\int_E |\wepsx| \le \sqrt{2c_1} \sqrt{\od(\delta)} \le \eta
	\qquad \mbox{for all } t\in\Big(0,\frac{1}{\sqrt{\eps}}\Big)
  \eas
  and that thus (\ref{824.1}) holds. 
  Furthermore, for $x\in\Omega$ and $y\in\Omega$ with $y<x$, a second application of (\ref{824.3}), now to 
  $F:=(y,x)$, shows that 
  \bas
	|\weps(x,t)-\weps(y,t)|
	= \bigg| \int_y^x \wepsx(z,t)dz \bigg|
	\le \int_y^x |\wepsx(z,t)|dz
	\le \sqrt{2c_1} \sqrt{\od(|x-y|)}
	\quad \mbox{for all } t\in\Big(0,\frac{1}{\sqrt{\eps}}\Big),
  \eas
  which together with a similar estimate in the case $y>x$ establishes (\ref{824.11}).
}
\qed
\subsection{A local estimate for $\uepsx$}
In order to {\cred ultimately} achieve pointwise convergence of $\ueps$ along a subsequence of $(\eps_j)_{j\in\N}$
through a compactness argument based on the Aubin-Lions lemma in Lemma \ref{lem25}, let us combine the 
weighted estimate for $(\sqrt{\deps\ueps})_x$ from Lemma \ref{lem8} with (\ref{mass}) and the boundedness properties
of $\deps$ inside $\{d>0\}$ to derive the following local but unweighted integral estimate for $\uepsx$ itself.
\begin{lem}\label{lem825}
  Let $K\subset \{d>0\}\cap\Omega$ be compact. Then there exists $C(K)>0$ such that for all $\epsin$,
  \be{825.1}
	\int_0^T \bigg\{ \int_K |\uepsx(x,t)| dx \bigg\}^2 dt \le C(K)\cdot (1+T)
	\qquad \mbox{for all } T>0.
  \ee
\end{lem}
\proof
  According to Lemma \ref{lem43}, our assumption on $K$ ensures that with some $c_1>0$ we have
  \be{825.2}
	\deps \ge c_1
	\quad \mbox{in $K$ \quad for all } \epsin,
  \ee
  and that moreover $\depsx\to d_x$ in $L^\infty(K)$, whence there exists $c_2>0$ fulfilling
  \be{825.3}
	|\depsx| \le c_2
	\quad \mbox{in $K$ \quad for all } \epsin.
  \ee
  We now make use of the fact that Lemma \ref{lem8} yields $c_3>0$ satisfying
  \bas
	\int_0^T \io \frac{(\deps\ueps)_x^2}{\deps\ueps} \le c_3\cdot (1+T)
	\qquad \mbox{for all $T>0$ and } \epsin,
  \eas
  which namely implies that for any such $\eps$ we have
  \bas
	c_3\cdot (1+T)
	&\ge& \int_0^T \int_K \frac{(\deps\uepsx + \depsx\ueps)^2}{\deps\ueps} \\
	&\ge& \int_0^T \int_K \frac{\frac{1}{2} \deps^2 \uepsx^2 - \depsx^2 \ueps^2}{\deps\ueps} \\
	&=& \frac{1}{2} \int_0^T \int_K \deps \frac{\uepsx^2}{\ueps}
	- \int_0^T \int_K \frac{\depsx^2}{\deps} \ueps
	\qquad \mbox{for all } T>0,
  \eas
  because $(\xi+\eta)^2 \ge \frac{1}{2}\xi^2  -\eta^2$ for all $\xi\in\R$ and $\eta\in\R$.
  In view of (\ref{825.2}), (\ref{825.3}) and (\ref{mass}), this shows that
  \bas
	\frac{c_1}{2} \int_0^T \int_K \frac{\uepsx^2}{\ueps}
	&\le& \frac{1}{2} \int_0^T \int_K \deps \frac{\uepsx^2}{\ueps} \\
	&\le& \int_0^T \int_K \frac{\depsx^2}{\deps} \ueps + c_3\cdot (1+T) \\
	&\le& \frac{c_2^2}{c_1} \int_0^T \io \ueps + c_3\cdot (1+T) \\
	&=& \frac{c_2^2}{c_1} T \io u_0 + c_3\cdot (1+T)
	\qquad \mbox{for all } T>0,
  \eas
  which readily implies (\ref{825.1}) upon the observation that
  \bas
	\int_0^t \bigg\{ \int_K |\uepsx| \bigg\}^2
	&\le& \int_0^T \bigg\{ \int_K \frac{\uepsx^2}{\ueps} \bigg\} \cdot \bigg\{ \int_K \ueps \bigg\} \\
	&\le& \bigg\{ \int_0^T \int_K \frac{\uepsx^2}{\ueps} \bigg\} \cdot \bigg\{ \io u_0 \bigg\}
	\qquad \mbox{for all } T>0
  \eas
  according to the Cauchy-Schwarz inequality and (\ref{mass}).
\qed
\subsection{Time regularity}
As a final preparation for our subsequence extraction, let us derive some regularity features of the respective
time derivatives. The first of these, again resulting from Lemma \ref{lem8},
is actually asymptotically independent of the length of the time interval
appearing therein, and hence can serve below as a first information on decay of temporal oscillations.
\begin{lem}\label{lem27}
  There exists $C>0$ such that for all $\epsin$,
  \be{27.1}
	\int_0^T \|u_{\eps t}(\cdot,t)\|_{(W^{1,\infty}(\Omega))^\star}^2 dt \le C\cdot (1+\sqrt{\eps}T)
	\qquad \mbox{for all } T>0.
  \ee
\end{lem}
\proof
  We fix $t>0$ and $\psi\in W^{1,\infty}(\Omega)$ such that $\|\psi\|_{W^{1,\infty}(\Omega)} \le 1$, and then obtain on
  testing the first equation in (\ref{0eps}) by $\psi$ and using the Cauchy-Schwarz inequality and (\ref{43.3})
  as well as (\ref{mass}) and (\ref{M}) that
  \bas
	\bigg| \io u_{\eps t}(\cdot,t)\psi\bigg|
	&=& \bigg| - \io (\deps\ueps)_x \psi_x + \io \deps\ueps\wepsx \psi_x \bigg| \\
	&\le& \io \Big|(\deps\ueps)_x \Big| + \io \deps\ueps |\wepsx| \\
	&\le& \bigg\{ \io \frac{(\deps\ueps)_x^2}{\deps\ueps} \bigg\}^\frac{1}{2} \cdot
	\bigg\{ \io \deps\ueps \bigg\}^\frac{1}{2}
	+ \bigg\{ \io \deps\ueps \frac{\wepsx^2}{\weps} \bigg\}^\frac{1}{2} \cdot
	\bigg\{ \io \deps\ueps\weps \bigg\}^\frac{1}{2} \\
	&\le& \bigg\{ \io \frac{(\deps\ueps)_x^2}{\deps\ueps} \bigg\}^\frac{1}{2} \cdot
	(\|d\|_{L^\infty(\Omega)}+1)^\frac{1}{2} \cdot \bigg\{ \io u_0\bigg\}^\frac{1}{2} \\
	& & + \bigg\{ \io \deps\ueps \frac{\wepsx^2}{\weps} \bigg\}^\frac{1}{2} \cdot
	(\|d\|_{L^\infty(\Omega)}+1)^\frac{1}{2} \cdot \bigg\{ \io u_0\bigg\}^\frac{1}{2} \cdot \sqrt{M}
  \eas
  for all $\epsin$.
  Writing $c_1:=(\|d\|_{L^\infty(\Omega)}+1) \cdot\io u_0$, we thus infer that for any such $\eps$,
  \bas
	\|u_{\eps t}(\cdot,t)\|_{(W^{1,\infty}(\Omega))^\star}^2
	\le 2c_1 \io \frac{(\deps\ueps)_x^2}{\deps\ueps}
	+ 2c_1 M \io \deps\ueps \frac{\wepsx^2}{\weps}
	\qquad \mbox{for all } t>0,
  \eas
  which in view of Lemma \ref{lem8} implies (\ref{27.1}) on integration in time.
\qed
Next, the estimates from Lemma \ref{lem8} imply the following temporally local estimate for $w_{\eps t}$ in a
straightforward manner.
\begin{lem}\label{lem827}
  Let $T>0$. Then there exists $C(T)>0$ such that 
  \be{827.1}
	\int_0^T \|w_{\eps t}(\cdot,t)\|_{L^1(\Omega)}^2 dt \le C(T)
	\qquad \mbox{for all } \epsin.
  \ee
\end{lem}
\proof
  By directly using the second equation in (\ref{0eps}) we can estimate
  \be{827.2}
	\bigg\{ \io |w_{\eps t}| \bigg\}^2
	\le 2\cdot \bigg\{ \eps\cdot \io \Big| \Big(\deps \frac{\wepsx}{\sqrt{g(\weps)}}\Big)_x \Big| \bigg\}^2
	+ 2\cdot\bigg\{ \io \ueps g(\weps) \bigg\}^2
	\qquad \mbox{for all } t>0,
  \ee
  where due to the Cauchy-Schwarz inequality, (\ref{g2}) and (\ref{M}),
  \bas
	\bigg\{ \eps\cdot \io \Big| \Big(\deps \frac{\wepsx}{\sqrt{g(\weps)}}\Big)_x \Big| \bigg\}^2
	&\le& \eps^2 \cdot 
	\bigg\{ \io \frac{1}{\sqrt{g(\weps)}} \cdot \Big(\deps\frac{\wepsx}{\sqrt{g(\weps)}}\Big)_x^2 \bigg\} \cdot
	\bigg\{ \io \sqrt{g(\weps)} \bigg\} \\
	&\le& \sqrt{\og M} |\Omega| \cdot \eps
	\io \frac{1}{\sqrt{g(\weps)}} \cdot \Big(\deps\frac{\wepsx}{\sqrt{g(\weps)}}\Big)_x^2
	\qquad \mbox{for all } t>0,
  \eas
  because $(\eps_j)_{j\in\N}\subset (0,1)$. 
  As (\ref{g2}) and (\ref{M}) together with (\ref{mass}) assert that
  \bas
	\bigg\{ \io \ueps g(\weps)\bigg\}^2
	\le \og^2 M^2 \cdot \bigg\{ \io \ueps \bigg\}^2
	= \og^2 M^2 \cdot \bigg\{ \io u_0 \bigg\}^2
	 \qquad \mbox{for all } t>0,
  \eas
  in view of Lemma \ref{lem8} we thus obtain (\ref{827.1}) from (\ref{827.2}).
\qed
\mysection{Global existence in the degenerate problem}\label{sect6}
\subsection{Construction of limit functions}
By means of a straightforward extraction procedure based on our estimates collected so far as well as
standard compactness arguments, we can now construct a limit object that will finally turn out to solve
(\ref{0}) in the considered generalized sense.
\begin{lem}\label{lem25}
  There exist a subsequence $(\eps_{j_k})_{k\in\N}$ of $(\eps_j)_{j\in\N}$  and nonnegative functions
  \be{25.1}
	\left\{ \begin{array}{l}
	u\in L^1_{loc}(\bom\times [0,\infty)) \cap C^0\big([0,\infty);(W^{1,\infty}(\Omega))^\star\big)
	\qquad \mbox{and} \\[1mm]
	w \in L^\infty(\Omega\times (0,\infty)) \cap C^0([0,\infty);L^1(\Omega)) \cap L^1_{loc}([0,\infty);W^{1,1}(\Omega))
	\end{array} \right.
  \ee
  such that 
  \begin{eqnarray}
	\ueps \to u
	\quad & & \mbox{a.e.~in } \Omega\times (0,\infty),
	\label{25.2} \\
	\ueps \to u
	\quad & & \mbox{in } L^1_{loc}(\overline{\Omega}\times [0,\infty)),
	\label{25.3} \\ 
	\ueps \to u
	\quad & & \mbox{in } C^0_{loc}([0,\infty);(W^{1,\infty}(\Omega))^\star),
	\label{25.33} \\ 
	\weps \to w
	\quad & & {\cred \mbox{in } C^0_{loc}(\bom\times [0,\infty)),}
	\label{25.44} \\
	\wepsx \wto w_x
	\quad & & \mbox{in } L^1_{loc}(\bom\times [0,\infty))
	\qquad \mbox{and}
	\label{25.5} \\
	\sqrt{\deps} \wepsx \wto \sqrt{d}w_x
	\quad & & \mbox{in } L^2_{loc}(\bom\times [0,\infty))
	\label{25.6}
  \eea
  as $\eps=\eps_{j_k} \searrow 0$.
\end{lem}
\proof
  We first combine Lemma \ref{lem825} with (\ref{mass}) to see that for any open $\Omega_0\subset \Omega$
  satisfying $\bom_0 \subset \{d>0\}\cap\Omega$,
  \bas
	(u_{\eps_j})_{j\in\N} 
	\quad \mbox{is bounded in $L^2((0,T);W^{1,1}(\Omega_0))$ for all } T>0,
  \eas
  whereas Lemma \ref{lem27}, asserting that
  \be{25.100}
	(u_{\eps_j t})_{j\in\N}
	\quad \mbox{is bounded in } L^2\Big( (0,T);(W^{1,\infty}(\Omega))^\star\Big) \mbox{ for all } T>0,
  \ee	
  entails that
  \bas
	(u_{\eps_j t})_{j\in\N}
	\quad \mbox{is bounded in } L^2\Big( (0,T);(W_0^{1,\infty}(\Omega_0))^\star\Big) \mbox{ for all } T>0
  \eas
  due to the observation that the trivial extension $\psi$ of any $\psi_0\in W_0^{1,\infty}(\Omega_0)$ to all of
  $\Omega$ satisfies $\psi\in W^{1,\infty}(\Omega)$ with 
  $\|\psi\|_{W^{1,\infty}(\Omega)} \le \|\psi_0\|_{W_0^{1,\infty}(\Omega_0)}$.
  For any such $\Omega_0$, in view of the compactness of the first of the embeddings 
  $W^{1,1}(\Omega_0) \hra L^2(\Omega) \hra (W_0^{1,\infty}(\Omega_0))^\star$ the Aubin-Lions lemma (\cite{temam})
  thus guarantees that
  \bas
	(u_{\eps_j})_{j\in\N} 
	\quad \mbox{is relatively compact in } L^2(\Omega_0\times (0,T)),
  \eas
  so that since $d$ is continuous in $\bom$, and since our assumption that $\frac{1}{d} \in L^1(\Omega)$ especially ensures
  that $d>0$ a.e.~in $\Omega$, by means of a straightforward successive extraction procedure we obtain a decreasing
  subsequence $(\eps_{j_k})_{k\in\N}$ of $(\eps_j)_{j\in\N}$ and a nonnegative measurable function
  $u: \Omega\times (0,\infty)\to\R$ such that (\ref{25.2}) holds.
  As from Lemma \ref{lem826} we particularly know that
  \bas
	(u_{\eps_j})_{j\in\N}
	\quad \mbox{is equi-integrable in $\Omega\times (0,T)$ for all } T>0,
  \eas
  due to (\ref{25.2}) we may invoke the Vitali convergence theorem to see that also (\ref{25.3}) holds 
  along this sequence.
  Moreover, combining (\ref{25.100}) with the fact that
  \bas
	(u_{\eps_j})_{j\in\N}
	\quad \mbox{is bounded in } L^\infty((0,T);L^1(\Omega)) \mbox{ for all } T>0
  \eas
  due to (\ref{mass}), we may make use of the compactness of the embedding 
  $L^1(\Omega)\hra (W^{1,\infty}(\Omega))^\star$ in employing the Arzel\`a-Ascoli theorem to conclude that 
  \bas
	(u_{\eps_j})_{j\in\N}
	\quad \mbox{is relatively compact in } C^0\Big([0,T];(W^{1,\infty}(\Omega))^\star\Big) \mbox{ for all } T>0,
  \eas
  and that hence on modification of $u$ on a null set of times we can also achieve (\ref{25.33}).\abs
  As for the second solution component, we first note that 
{\cred as a consequence of Lemma \ref{lem824},
  with $Y$ as introduced in (\ref{Y}) we have that
  \bas
	(w_{\eps_j})_{j\in\N}
	\quad \mbox{is bounded in } L^\infty((0,T);Y)
	\quad \mbox{for all } T>0,
  \eas
}
  so that since due to Lemma \ref{lem827},
  \bas
	(w_{\eps_j t})_{j\in\N}
	\quad \mbox{is bounded in } L^2((0,T);L^1(\Omega))
	\quad \mbox{for all } T>0,
  \eas
{\cred
  and since $Y$ is compactly emdedded into $C^0(\bom)$ according to the Arzel\`a-Acsoli theorem,
  another application of an Aubin-Lions lemma shows that
  \be{25.7}
	(w_{\eps_j})_{j\in\N}
	\quad \mbox{is relatively compact in } C^0([0,T];C^0(\bom))
	\quad \mbox{for all } T>0.
  \ee
  As combining Lemma \ref{lem824} with the Dunford-Pettis theorem apart from that warrants that
  \bas
	(w_{\eps_j x})_{j\in\N}
	\quad \mbox{is relatively compact with respect to the weak topology in } L^1(\Omega\times (0,T))
  \eas
  for all $T>0$, we may assume on passing to a further subsequence if necessary
  that also (\ref{25.44}) and (\ref{25.5}) hold,
  and since furthermore Lemma \ref{lem8} implies that
  \bas
	\Big(\sqrt{d_{\eps_j}} w_{\eps_j x}\Big)_{j\in\N}
	\quad \mbox{is bounded in } L^\infty((0,T);L^2(\Omega))
	\quad \mbox{for all } T>0,
  \eas
  upon a final extraction process we can also achieve (\ref{25.6}).
}
\qed
\subsection{Strong convergence of $\sqrt{\deps}\ueps$ in $L^2$}\label{sect6.2}
In view of (\ref{25.6}), for appropriate passing to the limit in the regularized counterpart of the haptotactic
integral in (\ref{w3}) it seems in order to assert strong convergence of the expression
$\sqrt{\deps}\ueps$ with respect to the norm in $L^2(\Omega\times (0,T))$ for fixed $T>0$.
In achieving this on the basis of the Vitali convergence theorem, we will make use of the
following generalization of the Gagliardo-Nirenberg inequality that can be obtained by straighforward
adaptation of the argument in \cite{biler_hebisch_nadzieja} (cf.~also \cite[Lemma A.5]{taowin_JDE2014}).
\begin{lem}\label{lem21}
  There exists $C>0$ such that for any choice of $\eta\in (0,1)$ one can find $C(\eta)>0$ with the property that
  \be{21.1}
	\|\varphi\|_{L^\infty(\Omega)}^{4}
	\le \eta \|\varphi_x\|_{L^2(\Omega)}^2 \cdot
	\Big\| \varphi \big|\ln |\varphi| \big|^\frac{1}{2} \Big\|_{L^2(\Omega)}^2
	+ C \|\varphi\|_{L^2(\Omega)}^{4}
	+ C(\eta)
	\qquad \mbox{for all } \varphi\in W^{1,2}(\Omega).
  \ee
\end{lem}
We can thereby once more exploit the estimates for $\ueps$ from Lemma \ref{lem8} to infer 
the following spatio-temporal equi-integrability property of $\deps\ueps^2$.
\begin{lem}\label{lem823}
  Let $T>0$. Then for all $\eta>0$ one can find $\delta>0$ such that for any choice of $\epsin$,
  \be{823.1}
	\int\int_E \deps\ueps^2 \le \eta
	\qquad \mbox{for all measurable $E\subset \Omega\times (0,T)$ fulfilling } |E|\le\delta.
  \ee
\end{lem}
\proof
  In conclusion of Lemma \ref{lem8}, we can fix $c_1>0$ and $c_2>0$ such that for all $\epsin$ we have
  \be{823.2}
	\io \ueps \Big|\ln (\deps\ueps)\Big| \le c_1
	\qquad \mbox{for all } t\in (0,T)
  \ee
  and
  \be{823.3}
	\int_0^T \io \frac{(\deps\ueps)_x^2}{\deps\ueps} \le c_2.
  \ee
  Then for arbitrary $\eta>0$, applying Lemma \ref{lem21} and using that 
  $c_3:=\io \frac{1}{d}$ is finite,
  we may pick $c_4>0$ such that
  \be{823.4}
	\|\varphi\|_{L^\infty(\Omega)}^4
	\le \frac{4\eta}{c_1 c_2 c_3(\|d\|_{L^\infty(\Omega)}+1)} \|\varphi_x\|_{L^2(\Omega)}^2 
	\Big\| \varphi \ln^\frac{1}{2} |\varphi| \Big\|_{L^2(\Omega)}^2
	+ c_4 \|\varphi\|_{L^2(\Omega)}^4 + c_4
	\qquad \mbox{for all } \varphi\in W^{1,2}(\Omega)
  \ee
  and abbreviate $c_5:=c_4(\|d\|_{L^\infty\Omega)}+1)^2 \Big\{ \io u_0\Big\}^2 +c_4$.
  Now once more since $\frac{1}{d}\in L^1(\Omega)$ {\cred and hence also $\frac{1}{d} \in L^1(\Omega\times (0,T))$,
  we can find $\delta>0$ such that
  \be{823.111}
	\int\int_E \frac{1}{d} \le \frac{\eta}{2c_5}
	\quad \mbox{for each measurable $E\subset\Omega\times (0,T)$ satisfying } |E|\le \delta.
  \ee
  In order to derive (\ref{823.1}) from this, we observe that by (\ref{823.4}),
  \bea{823.112}
	\deps\ueps^2
	&=& \frac{\sqrt{\deps\ueps}^4}{\deps} \nn\\
	&\le& \frac{4\eta}{c_1 c_2 c_3 (\|d\|_{L^\infty(\Omega)}+1) \deps}
	\Big\| (\sqrt{\deps\ueps})_x \Big\|_{L^2(\Omega)}^2 
	\Big\| \sqrt{\deps\ueps} \cdot \sqrt{|\ln \sqrt{\deps\ueps}|} \Big\|_{L^2(\Omega)}^2
	+ \frac{c_4}{\deps} \| \sqrt{\deps\ueps}\|_{L^2(\Omega)}^4
	+ \frac{c_4}{\deps} \nn\\
	&=& \frac{\eta}{2c_1 c_2 c_3(\|d\|_{L^\infty(\Omega)}+1) \deps}
	\cdot \bigg\{ \io \frac{(\deps\ueps)_x^2}{\deps\ueps} \bigg\}
	\cdot \bigg\{ \io \deps\ueps \Big|\ln (\deps\ueps)\Big| \bigg\} \nn\\
	& & + \frac{c_4}{\deps} \cdot \bigg\{ \io \deps\ueps \bigg\}^2
	+ \frac{c_4}{\deps}
	\qquad \mbox{for all $x\in\Omega$ and } t>0,
  \eea
  where according to (\ref{823.2}), (\ref{43.3}) and (\ref{mass}),
  \bas
	\io \deps\ueps \Big|\ln (\deps\ueps)\Big|
	\le (\|d\|_{L^\infty(\Omega)} +1) \cdot \io \ueps \Big|\ln (\deps\ueps)\Big|
	\le (\|d\|_{L^\infty(\Omega)} +1) \cdot c_1
	\qquad \mbox{for all } \in (0,T)
  \eas
  and
  \bas
	\frac{c_4}{\deps} \cdot \bigg\{ \io \deps\ueps\bigg\}^2 + \frac{c_4}{\deps}
	\le \frac{c_4}{\deps} \cdot (\|d\|_{L^\infty(\Omega)}+1)^2 \cdot \bigg\{ \io u_0\bigg\}^2
	+ \frac{c_4}{\deps}
	= \frac{c_5}{\deps}
	\qquad \mbox{for all $x\in\Omega$ and } t>0.
  \eas
  Therefore, given any measurable $E\subset (0,T)$ with $|E|\le\delta$, we infer on integrating (\ref{823.112}) that due
  to (\ref{823.3}) and (\ref{823.111}), indeed we have
  \bas
	\int\int_E \deps\ueps^2
	&\le& \frac{\eta}{2c_2 c_3} \cdot \int\int_E \frac{1}{\deps} 
	\cdot \bigg\{\io \frac{(\deps\ueps)_x^2}{\deps\ueps} \bigg\}
	+ c_5 \int\int_E \frac{1}{\deps} \\
	&\le& \frac{\eta}{2c_2 c_3} \cdot \int_0^T \bigg\{ \io \frac{1}{d} \bigg\} \cdot
	\bigg\{ \frac{(\deps\ueps)_x^2}{\deps\ueps} \bigg\}
	+ c_5 \int\int_E \frac{1}{d} \\
	&\le& \frac{\eta}{2} + \frac{\eta}{2}=\eta
	\qquad \mbox{for all } \epsin,
  \eas
  again because $\deps\ge d$.
}
\qed
In consequence, the Vitali convergence theorem entails the desired strong convergence feature.
\begin{lem}\label{lem822}
  With $(\eps_{j_k})_{k\in\N}$ taken from Lemma \ref{lem25}, we have
  \bas
	\sqrt{\deps} \ueps \to \sqrt{d}u
	\quad \mbox{in } L^2_{loc}(\bom\times [0,\infty)) 
	\qquad \mbox{as } \eps=\eps_{j_k}\searrow 0.
  \eas
\end{lem}
\proof
  In view of the Vitali convergence theorem, this is a direct consequence of Lemma \ref{lem823}
  when combined with the fact that due to Lemma \ref{lem43} and Lemma \ref{lem25} we have
  $\sqrt{\deps}\ueps\to \sqrt{d}u$ a.e.~in $\Omega\times (0,\infty)$ as
  $\eps=\eps_{j_k}\searrow 0$.
\qed
\subsection{Global existence in (\ref{0})}
We are now prepared for appropriate limit procedures in each of the integrals related to (\ref{w3}) and (\ref{w4}).
\begin{lem}\label{lem255}
  The pair $(u,w)$ obtained in Lemma \ref{lem25} is a global generalized solution of (\ref{0}) in the sense
  of Definition \ref{defi_weak}.
\end{lem}
\proof
  The regularity properties in (\ref{w1}) are implied by (\ref{25.1}), whereas if we take $(\eps_{j_k})_{k\in\N}$
  as provided by Lemma \ref{lem25}, then the strong $L^2$ convergence property of
  $(\sqrt{d_{\eps_{j_k}}} u_{\eps_{j_k}})_{k\in\N}$ asserted by Lemma \ref{lem822} along with the weak $L^2$
  approximation feature of $(\sqrt{d_{\eps_{j_k}}} w_{\eps_{j_k} x})_{k\in\N}$ gained in Lemma \ref{lem25} warrants that
  \be{255.1}
	\deps\ueps\wepsx
	= (\sqrt{\deps}\ueps)\cdot(\sqrt{\deps}\wepsx)
	\wto duw_x
	\qquad \mbox{in } L^1_{loc}(\bom\times [0,\infty))
  \ee
  as $\eps=\eps_{j_k}\searrow 0$, and that hence also (\ref{w2}) holds.\\
  The verification of the integral identity (\ref{w3}) is now straightofrward:
  Fixing an arbitrary
  $\varphi\in C_0^\infty(\overline{\Omega}\times [0,\infty))$ such that $\varphi_x=0$ on $\pO\times (0,\infty)$,
  we obtain from the first equation in (\ref{0eps}) that for each $\epsin$,
  \be{255.2}
	- \int_0^\infty \io \ueps \varphi_t
	- \io u_0\varphi(\cdot,0)
	= \int_0^\infty \io \deps\ueps \varphi_{xx}
	+ \int_0^\infty \io \deps \ueps \wepsx \varphi_x,
  \ee
  where (\ref{255.1}) ensures that
  \bas
	\int_0^\infty \io \deps\ueps \wepsx \varphi_x
	\to \int_0^\infty \io duw_x \varphi_x
	\qquad \mbox{as } \eps=\eps_{j_k}\searrow 0,
  \eas
  whereas using that $\deps\to d$ in $L^\infty(\Omega)$ as $\eps=\eps_j\searrow 0$ we infer from (\ref{25.3}) that
  \bas
	-\int_0^\infty \io \ueps \varphi_t \to - \int_0^\infty \io u\varphi_t
	\quad \mbox{and} \quad
	- \int_0^\infty \io \deps\ueps \varphi_{xx} 
	\to \int_0^\infty \io du\varphi_{xx}
	\qquad \mbox{as } \eps=\eps_{j_k}\searrow 0,
  \eas
  so that (\ref{255.2}) entails (\ref{w3}).\\
  Likewise, for fixed
  $\varphi\in C_0^\infty(\overline{\Omega}\times [0,\infty))$ the second equation in
  (\ref{0eps}) yields
  \be{255.3}
	- \int_0^\infty \io \weps \varphi_t
	- \io w_{0\eps} \varphi(\cdot,0)
	= - \eps \int_0^\infty \io \deps \frac{\wepsx}{\sqrt{g(\weps)}} \varphi_x
	- \int_0^\infty \io \ueps g(\weps)\varphi
  \ee
  for all $\epsin$, 
  where according to our construction of $(w_{0\eps_j})_{j\in\N}$ in Lemma \ref{lem45} and Lemma \ref{lem42} we know that
  \bas
	- \io w_{0\eps_{j_k}} \varphi(\cdot,0)
	&=& - \io w_{0 j_k} \varphi(\cdot,0)
	- \eps_{j_k}^\frac{1}{4} \io \varphi(\cdot,0) \\
	&\to& -\io w_0 \varphi(\cdot,0)
	\qquad \mbox{as } k\to\infty,
  \eas
{\cred 
  and where according to the uniform convergence statement in {\cred (\ref{25.44})}, the $L^1$ approximation property
  (\ref{25.3})
}
  and the continuity of $g$ on $[0,\infty)$,
  \bas
	-\int_0^\infty \io \weps\varphi_t 
	\to - \int_0^\infty \io w\varphi_t
	\qquad \mbox{as } \eps=\eps_{j_k}\searrow 0
  \eas
  and
  \bas
	-\int_0^\infty \io \ueps g(\weps)\varphi
	\to - \int_0^\infty \io ug(w)\varphi
	\qquad \mbox{as } \eps=\eps_{j_k}\searrow 0.
  \eas
  Since the Cauchy-Schwarz inequality together with (\ref{g2}) implies that
  \bas
	\bigg| - \eps \int_0^\infty \io \deps \frac{\wepsx}{\sqrt{g(\weps)}} \varphi_x \bigg|
	&\le& \eps \cdot \bigg\{ \int_0^\infty \io \frac{1}{\sqrt{g(\weps)}} 
		\cdot \Big(\deps \frac{\wepsx}{\sqrt{g(\weps)}}\Big)_x^2 \bigg\}^\frac{1}{2} \cdot 
	\bigg\{ \int_0^\infty \io \sqrt{g(\weps)} \varphi^2 \bigg\}^\frac{1}{2} \\
	&\le& \sqrt{\eps} \cdot \bigg\{ \sqrt{\eps} \int_0^\infty \io \frac{1}{\sqrt{g(\weps)}}
		\cdot \Big(\deps \frac{\wepsx}{\sqrt{g(\weps)}}\Big)_x^2 \bigg\}^\frac{1}{2} \cdot 
	(\og M)^\frac{1}{4} \cdot \bigg\{ \int_0^\infty \io \varphi^2 \bigg\}^\frac{1}{2}
  \eas
  for all $\epsin$, and that hence
  \bas
	- \eps \int_0^\infty \io \deps \frac{\wepsx}{\sqrt{g(\weps)}} \varphi_x
	\to 0
	\qquad \mbox{as } \eps=\eps_j\searrow 0
  \eas
  thanks to Lemma \ref{lem8}, on taking $\eps=\eps_{j_k}\searrow 0$ in (\ref{255.3}) we also obtain (\ref{w4}).
\qed
\mysection{Further regularity properties of $(u,w)$. Proof of Theorem \ref{theo700}}\label{sect7}
{\cred
In order to complete the proof of Theorem \ref{theo700}, but also to further prepare our subsequent asymptotic analysis,
let us use the equi-integrability and equicontinuity properties contained in Section \ref{sect_equi} to firstly 
derive corresponding conclusions
for the respective limit functions, and to secondly assert the continuity and mass conservation
properties claimed in Theorem \ref{theo700}.
}
{\cred 
\begin{lem}\label{lem333}
  Let $(\eps_{j_k})_{k\in\N}$ be as in Lemma \ref{lem25}. 
  The solution component $u$ belongs to $C^0_w([0,\infty);L^1(\Omega))$,
  and with $(\eps_{j_k})_{k\in\N}$ taken from Lemma \ref{lem25}, we have
  \be{333.1}
	\ueps(\cdot,t) \wto u(\cdot,t)
	\quad \mbox{in } L^1(\Omega)
	\qquad \mbox{as } \eps=\eps_{j_k}\searrow 0.
  \ee
  Moreover,
  \be{333.2}
	(u(\cdot,t))_{t>0}
	\quad \mbox{is equi-integrable over } \Omega
  \ee
  and
  \be{333.3}
	(w(\cdot,t))_{t>0}
	\quad \mbox{is equicontinuous in } \bom.
  \ee
\end{lem}
\proof
  Once more using that with $Y$ as in (\ref{Y}), the family $(w_{\eps_j})_{j\in\N}$ is bounded in 
  $L^\infty((0,\infty);Y)$, we directly see from (\ref{25.44}) that $(w(\cdot,t))_{t>0}$ is bounded in $Y$
  and hence equicontinuous in $\bom$ according to (\ref{Y}).\abs
  Next, fixing an arbitrary 
  $t>0$ we know from (\ref{25.33}) that $\ueps(\cdot,t)\to u(\cdot,t)$ in $(W^{1,\infty}(\Omega))^\star$
  as $\eps=\eps_{j_k}\searrow 0$, whereas Lemma \ref{lem826} shows that $(u_{\eps_j}(\cdot,t))_{j\in\N}$
  is relatively compact with respect to the weak topology in $L^1(\Omega)$ due to the Dunford-Pettis theorem.
  Combining these two properties implies that any accumulation point of $(u_{\eps_{j_k}}(\cdot,t))_{k\in\N}$
  in the weak topology of $L^1(\Omega)$ must coincide with $u(\cdot,t)$, hence implying that $u(\cdot,t)\in L^1(\Omega)$
  and $\ueps(\cdot,t)\wto u(\cdot,t)$ in $L^1(\Omega)$ along the entire sequence $\eps=\eps_{j_k}\searrow 0$.
  Having thus verified (\ref{333.1}), in view of the fact that this entails $\int_E \ueps(\cdot,t) \to \int_E u(\cdot,t)$
  as $\eps=\eps_{j_k}\searrow 0$ for each measurable $E\subset\Omega$, we immediately also obtain
  (\ref{333.2}) as a consequence of Lemma \ref{lem826}.
  Finally, the inclusion $u\in C^0_w([0,\infty);L^1(\Omega))$ can be seen by quite a similar argument:
  Given $t_0\ge 0$ and
  $(t_k)_{k\in\N}\subset (0,\infty)$ such that $t_k\to t_0$ as $k\to\infty$, again relying on (\ref{25.33})
  we note that $u(\cdot,t_k)\to u(\cdot,t)$ in $(W^{1,\infty}(\Omega))^\star$ as $k\to\infty$, whereas
  (\ref{333.2}) in conjunction with the Dunford-Pettis theorem warrants that $(u(\cdot,t_k))_{k\in\N}$ is
  relatively compact with respect to the weak topology in $L^1(\Omega)$. As thus $u(\cdot,t)$ is the only
  cluster point of $(u(\cdot,t_k))_{k\in\N}$ in the latter space, we infer that indeed $u(\cdot,t_k)\wto u(\cdot,t_0)$
  in $L^1(\Omega)$ as $k\to\infty$.
\qed
}
Thus particularly knowing that not only $w(\cdot,t)$ but also $u(\cdot,t)$ is a well-defined element of $L^1(\Omega)$ for all $t>0$, we can proceed to formulate corresponding dissipation and conservation properties in this space,
both of which being of great importance for our stabilization proof below.
\begin{lem}\label{lem256}
  We have
  \be{massu}
	\io u(\cdot,t) = \io u_0
	\qquad \mbox{for all } t>0
  \ee
  and
  \be{256.2}
	\|w(\cdot,t)\|_{L^1(\Omega)} \le \|w(\cdot,t_0)\|_{L^1(\Omega)} \le \|w_0\|_{L^1(\Omega)}
	\qquad \mbox{for all $t_0>0$ and any } t>t_0
  \ee
  as well as
  \be{256.3}
	\int_0^\infty \io uw \le \frac{1}{\ug} \io w_0.
  \ee
\end{lem}
\proof
  The conservation property (\ref{massu}) is an immediate consequence of (\ref{mass}) and {\cred Lemma \ref{lem333}.}
  For the derivation of (\ref{256.2}) and (\ref{256.3}) we integrate the second equation in (\ref{0eps})
  to see that
  \be{256.4}
	\frac{d}{dt} \io \weps = - \io \ueps g(\weps)
	\qquad \mbox{for all } t>0,
  \ee
  whence in particular
  \be{256.5}
	\io \weps(\cdot,t) \le \io \weps(\cdot,t_0) \le \io w_{0\eps}
	\qquad \mbox{for all $t_0>0$ and any } t\in (t_0,\infty).
  \ee
  Recalling that by Lemma \ref{lem45} and Lemma \ref{lem42},
  \be{256.6}
	\io w_{0\eps_j} = \io w_{0j} + \eps_j^\frac{1}{4} |\Omega| \to \io w_0
	\qquad \mbox{as } j\to\infty,
  \ee
  in view of (\ref{25.44}) we thus obtain (\ref{256.2}) from (\ref{256.5}).\\
  Finally, further integration of (\ref{256.4}) shows that due to (\ref{g2}),

  \bas
	\ug \int_0^t \io \ueps\weps
	\le \int_0^t \io \ueps g(\weps)
	= \io w_{0\eps} - \io \weps(\cdot,t) \le \io w_{0\eps}
	\qquad \mbox{for all } t>0,
  \eas
  so that (\ref{256.6}) along with (\ref{25.2}) and {\cred (\ref{25.44})} 
  establishes (\ref{256.3}) by means of Fatou's lemma.
\qed
The proof of our main result on global existence, regularity and mass conservation is thereby complete:\abs
\proofc of Theorem \ref{theo700}.\quad
  In Lemma \ref{lem255} we have seen that $(u,w)$ is a global generalized solution of (\ref{0}) in the desired sense.
  The additional boundedness and continuity properties in (\ref{700.01}) as well as the mass conservation law
  (\ref{700.9}) readily result from {\cred Lemma \ref{lem333}, Lemma \ref{lem25} and Lemma \ref{lem256}.}
\qed
\mysection{Stabilization. Proof of Theorem \ref{theo7000}}\label{sect8}
We next intend to properly exploit the global dissipative properties expressed in Lemma \ref{lem8},
Lemma \ref{lem27} and Lemma \ref{lem256} so as to derive the convergence results claimed in Theorem \ref{theo7000}.
We will first concentrate on the respective statement concerning $u$ and thereafter consider the decay of the
component $w$.
\subsection{{\cred An averaged stabilization property of $u$}}\label{sect8.1}
Let us first state a conseqence of Lemma \ref{lem27} for the limit $u$ in a form which does no longer involve 
time derivatives but rather 
concentrates on the quantity $u$ itself, but which still reflects an appropriate relaxation property in the large time limit.
{\cred
The argument underlying the following lemma was kindly pointed out to us by one of the reviewers.
}
{\cred
\begin{lem}\label{lem100}
  For each $\varphi\in L^\infty(\Omega)$, we have
  \be{100.1}
	\sup_{\tau\in [0,1]} \bigg| \io \Big(u(\cdot,t+\tau)-u(\cdot,t)\Big) \cdot\varphi \ \bigg| \to 0
	\qquad \mbox{as } t\to\infty.
  \ee
\end{lem}
\proof
  Given $\eta>0$, thanks to the equi-integrability property (\ref{333.2}) we can fix $\delta>0$ such that whenever
  $E\subset\Omega$ is measurable with $|E|\le \delta$, we have
  \be{100.2}
	2 \cdot \Big( 2\|\varphi\|_{L^\infty(\Omega)}+1\Big) \int_E u(\cdot,t) \le \frac{\eta}{4}
	\qquad \mbox{for all } t>0.
  \ee
  Next, employing a standard regularization procedure we can find $(\varphi_k)_{k\in\N}\subset X:=W^{1,\infty}(\Omega)$
  such that 
  \be{100.22}
	\|\varphi_k\|_{L^\infty(\Omega)} \le \|\varphi\|_{L^\infty(\Omega)} +1
	\quad \mbox{for all } k\in\N
	\qquad \mbox{and} \qquad
	\varphi_k\to\varphi
	\quad \mbox{a.e.~in $\Omega$ as } k\to\infty.
  \ee
  Due to Egorov's theorem, the latter approximation property in particular enables us to pick $k_0\in\N$ and a
  measurable $E\subset\Omega$ such that $|E|\le\delta$ and
  \be{100.3}
	2\cdot \bigg\{\io u_0 \bigg\} \cdot \|\varphi-\varphi_{k_0}\|_{L^\infty(\Omega\setminus E)} \le \frac{\eta}{4}.
  \ee
  Finally, Lemma \ref{lem27} asserts the existence of $c_1>0$ such that
  \bas
	\int_0^T \|u_{\eps t}(\cdot,t)\|_{X^\star}^2 dt \le c_1\cdot (1+\sqrt{\eps}T),
	\qquad \mbox{for all } T>0
	\mbox{ and any } \epsin,
  \eas
  from which it readily follows by means of Lemma \ref{lem25} and a lower semicontinuity argument that
  $\int_0^\infty \|u_t(\cdot,t)\|_{X^\star}^2 dt<\infty$ and that hence we can choose
  $t_0>0$ large enough fulfilling
  \be{100.4}
	\|\varphi_{k_0}\|_X \cdot \int_{t_0}^\infty \|u_t(\cdot,t)\|_{X^\star}^2 dt \le \frac{\eta}{2}.
  \ee
  Now decomposing the expression under consideration according to
  \bas	
	\io \Big(u(\cdot,t+\tau)-u(\cdot,t)\Big) \cdot\varphi
	= \io \Big(u(\cdot,t+\tau)-u(\cdot,t)\Big) \cdot\varphi_{k_0}
	+ \io \Big(u(\cdot,t+\tau)-u(\cdot,t)\Big) \cdot (\varphi-\varphi_{k_0})
  \eas
  for $t>0$ and $\tau\in [0,1]$, by using the Cauchy-Schwarz inequality and (\ref{100.4}) we may estimate
  \bas
	\bigg| 	\io \Big(u(\cdot,t+\tau)-u(\cdot,t)\Big) \cdot\varphi_{k_0} \bigg|
	&=& \bigg| \int_t^{t+\tau} \langle u_t(\cdot,s),\varphi_{k_0} \rangle ds \bigg| \\
	&\le& \|\varphi_{k_0}\|_{X} \int_t^{t+\tau} \|u_t(\cdot,s)\|_{X^\star} ds \\
	&\le& \|\varphi_{k_0}\|_{X} \int_t^\infty \|u_t(\cdot,s)\|_{X^\star}^2 ds \\
	&\le& \frac{\eta}{2}
	\qquad \mbox{for all $t\ge t_0$ and any } \tau\in [0,1],
  \eas
  where $\langle \cdot , \cdot, \rangle$ denotes the duality pairing between $X^\star$ and $X$.
  Since furthermore from (\ref{massu}), (\ref{100.22}), (\ref{100.3}) and (\ref{100.2}) we know that  
  \bas
	\bigg| 	\io \Big(u(\cdot,t+\tau)-u(\cdot,t)\Big) \cdot (\varphi-\varphi_{k_0})\bigg|
	&\le& \bigg\{ \int_{\Omega\setminus E} \Big(u(\cdot,t+\tau)-u(\cdot,t)\Big) \bigg\}
	\cdot \|\varphi-\varphi_{k_0}\|_{L^\infty(\Omega\setminus E)} \\
	& & + \|\varphi-\varphi_{k_0}\|_{L^\infty(E)} \cdot \bigg\{ \int_E u(\cdot,t+\tau) + \int_E u(\cdot,t) \bigg\} \\
	&\le& 2\cdot\bigg\{ \io u_0 \bigg\} \cdot \|\varphi-\varphi_{k_0}\|_{L^\infty(\Omega\setminus E)} \\
	& & + \Big(2\|\varphi\|_{L^\infty(\Omega)} +1\Big)
	\cdot \bigg\{ \int_E u(\cdot,t+\tau) + \int_E u(\cdot,t) \bigg\} \\
	&\le& \frac{\eta}{4} + \frac{\eta}{4} = \frac{\eta}{2}
	\qquad \mbox{for all $t>0$ and each } \tau\in  [0,1],
  \eas
  we thus infer that
  \bas
	\bigg| 	\io \Big(u(\cdot,t+\tau)-u(\cdot,t)\Big) \cdot\varphi \bigg|
	\le \eta 
	\qquad \mbox{for all $t\ge t_0$ and } \tau\in [0,1],
  \eas
  as intended.
\qed
}

\subsection{Decaying deviation of $du$ from its spatial average}
Next aiming at a direct exploitation of (\ref{8.3}), in view of the fact that through a Poincar\'e inequality
the spatial gradients appearing therein control appropriate $L^p$ norms of deviations from respective spatial means,
let us briefly address the spatial averages relevant to our approach in the following.
\begin{lem}\label{lem922}
  The function $\mu$ defined on $[0,\infty)$ by letting
  \be{mu}
	\mu(t):=\frac{1}{|\Omega|} \io d(x)u(x,t)dx,
	\qquad t>0,
  \ee
  is bounded and continuous on $[0,\infty)$, and with $(\eps_{j_k})_{k\in\N}$ as provided by Lemma \ref{lem25}
  we have
  \be{922.1}
	\mu_\eps(t) \to \mu(t)
{\cred	\quad \mbox{for all $t>0$}}	
	\qquad \mbox{as } \eps=\eps_{j_k}\searrow 0,
  \ee
  where we have set 
  \bas
	\mu_\eps(t):=\frac{1}{|\Omega|} \io \deps(x)\ueps(x,t)dx, \qquad t\ge 0, \ \epsin.
  \eas
{\cred
  Moreover, 
  \be{103.1}
	\sup_{\tau\in [0,1]} \Big|\mu(t+\tau)-\mu(t)\Big| \to 0
	\qquad \mbox{as } t\to\infty.
  \ee
}
\end{lem}
\proof
  As $d$ is bounded, the continuity of $\mu$ is an immediate consequence of {\cred Lemma \ref{lem333}}, 
  whereas its boundedness is evident from (\ref{massu}).
  The approximation property (\ref{922.1}) results upon observing that 
  {\cred Lemma \ref{lem333} asserts that as $\eps=\eps_{j_k}\searrow 0$, for all $t>0$ we have
  $\ueps(\cdot,t)\wto u(\cdot,t)$ in $L^1(\Omega)$} and hence also $\deps\ueps(\cdot,t)\to du(\cdot,t)$
  in $L^1(\Omega)$ due to the fact that $\deps\to d$ in $L^\infty(\Omega)$ by Lemma \ref{lem43}.
  {\cred Finally, (\ref{103.1}) directly results on applying Lemma \ref{lem100} to $\varphi:=d$.}
\qed
In terms of the function $\mu$ thus defined, (\ref{8.3}) implies the following.
\begin{lem}\label{lem102}
  With $\mu$ as defined in (\ref{mu}), we have
  \be{102.1}
	\int_0^\infty \|du(\cdot,t)-\mu(t)\|_{L^1(\Omega)}^2 dt < \infty.
  \ee
\end{lem}
\proof
  According to a Poincar\'e inequality we can find $c_1>0$ such that
  \bas
	\|\varphi-\overline{\varphi}\|_{L^1(\Omega)} \le c_1\|\varphi_x\|_{L^1(\Omega)}
	\qquad \mbox{for all } \varphi\in W^{1,1}(\Omega),
  \eas
  so that for arbitrary $T>0$
  we may once more combine the Cauchy-Schwarz inequality with (\ref{43.3}) and (\ref{mass})
  to see that with $\mu_\eps$ as introduced in Lemma \ref{lem922} we have
  \bas
	\int_0^T \|\deps\ueps(\cdot,t)-\mu_\eps(t)\|_{L^1(\Omega)}^2 dt
	&\le& c_1^2 \int_0^T \bigg\{ \io \Big| (\deps\ueps(\cdot,t))_x \Big| \bigg\}^2 dt \\
	&\le& c_1^2 \int_0^T \bigg\{ \io \frac{(\deps\ueps)_x^2}{\deps\ueps} \bigg\} \cdot
	\bigg\{ \io \deps\ueps \bigg\} \\
	&\le& c_2 \int_0^T \io \frac{(\deps\ueps)_x^2}{\deps\ueps}
	\qquad \mbox{for all } \epsin
  \eas
  with $c_2:=c_1^2(\|d\|_{L^\infty(\Omega)}+1)\io u_0$. Since Lemma \ref{lem8} provides $c_3>0$ such that
  \bas
	\int_0^T \io \frac{(\deps\ueps)_x^2}{\deps\ueps} \le c_3\cdot (1+\sqrt{\eps}T)
	\qquad \mbox{for all $T>0$ and } \epsin,
  \eas
  from this we infer that for all $T>0$,
  \be{102.3}
	\int_0^T \|\deps\ueps(\cdot,t)-\mu_\eps(t)\|_{L^1(\Omega)}^2 dt
	\le 2c_2 c_3
	\qquad \mbox{whenever $\epsin$ is such that } \eps\le \frac{1}{T^2}.
  \ee
  We now use that as a particular consequence of Lemma \ref{lem25} we have $\deps\ueps \to du$ in 
  $L^2_{loc}([0,\infty);L^1(\Omega))$ as $\eps=\eps_{j_k}\searrow 0$, which together with Lemma \ref{lem922} guarantees that
  for all $T>0$ and any $\tau\in (0,T)$,
  \bas
	\int_\tau^T \|\deps\ueps(\cdot,t)-\mu_\eps(t)\|_{L^1(\Omega)}^2 dt	
	\to \int_\tau^T \|du(\cdot,t)-\mu(t)\|_{L^1(\Omega)}^2 dt
	\qquad \mbox{as } \eps=\eps_{j_k}\searrow 0.
  \eas
  Therefore, (\ref{102.3}) implies that
  \bas
	\int_\tau^T \|du(\cdot,t)-\mu(t)\|_{L^1(\Omega)}^2 dt 
	\le 2c_2 c_3
	\qquad \mbox{for all $T>0$ and $\tau\in (0,T)$},
  \eas
  from which (\ref{102.1}) results on taking $\tau\searrow 0$ and $T\to\infty$.
\qed
{\cred
Once again relying on Lemma \ref{lem100}, we thereby indeed arrive at the main result of this section.
}
\begin{lem}\label{lem104}
  With $\mu$ as defined in (\ref{mu}), we have
  \be{104.1}
	du(\cdot,t)-\mu(t) \wto 0
	\quad \mbox{in } L^1(\Omega)
	\qquad \mbox{as } t\to\infty.
  \ee
\end{lem}
\proof
{\cred
  We fix $\varphi\in L^\infty(\Omega)$ and $\eta>0$ and then obtain from Lemma \ref{lem100} that there exists
  $t_1>0$ such that
  \be{104.2}
	\sup_{\tau\in [0,1]} \bigg| \io d\varphi \cdot \Big(u(\cdot,t+\tau)-u(\cdot,t)\Big) \bigg| 
	\le \frac{\eta}{3}
	\qquad \mbox{for all } t\ge t_1,
  \ee
  whereas (\ref{103.1}) says that with some $t_2\ge t_1$ we have
  \be{104.3}
	\|\varphi\|_{L^1(\Omega)} \cdot \sup_{\tau\in [0,1]} \Big|\mu(t+\tau)-\mu(t)\Big| \le \frac{\eta}{3}
	\qquad \mbox{for all } t\ge t_2,
  \ee
  and finally invoking Lemma \ref{lem102} we can pick $t_0\ge t_2$ satisfying
  \be{104.4}
	\|\varphi\|_{L^\infty(\Omega)} \cdot
	\int_{t_0}^\infty \|du(\cdot,t)-\mu(t)\|_{L^1(\Omega)}^2 dt \le \frac{\eta}{3}.
  \ee
  We now write
  \bea{104.5}
	\io \Big(du(\cdot,t)-\mu(t)\Big)\cdot\varphi
	&=& \int_0^1 \io \Big( d(x) u(x,t) - \mu(t)\Big)\cdot\varphi(x) dxd\tau \nn\\
	&=& \int_0^1 \io d(x) \Big( u(x,t)-u(x,t+\tau)\Big)\cdot\varphi(x) dxd\tau \nn\\
	& & + \int_0^1 \io \Big(d(x) u(x,t+\tau)-\mu(t+\tau)\Big) \cdot\varphi(x) dxd\tau \nn\\
	& & + \int_0^t \Big(\mu(t+\tau)-\mu(t)\Big) \cdot \io \varphi(x) dx d\tau
	\qquad \mbox{for } t>0,
  \eea
  and use (\ref{104.2}) to see that herein for all $t\ge t_0\ge t_1$,
  \bas
	\bigg| \int_0^1 \io d(x) \Big( u(x,t)-u(x,t+\tau)\Big)\cdot\varphi(x) dxd\tau \bigg|
	\le \sup_{\tau\in [0,1]}  \bigg| \io d\varphi \cdot \Big(u(\cdot,t+\tau)-u(\cdot,t)\Big) \bigg| 
	\le \frac{\eta}{3}.
  \eas
  Moreover, (\ref{104.3}) entails that
  \bas
	\bigg| \int_0^t \Big(\mu(t+\tau)-\mu(t)\Big) \cdot \io \varphi(x) dx d\tau \bigg|
	\le \|\varphi\|_{L^1(\Omega)} \cdot \sup_{\tau\in [0,1]} \Big|\mu(t+\tau)-\mu(t)\Big| 
	\le \frac{\eta}{3}
	\qquad \mbox{for all } t\ge t_0\ge t_2,
  \eas
  while combining the Cauchy-Schwarz inequality with (\ref{104.4}) shows that
  \bas
	\bigg|\int_0^1 \io \Big(d(x) u(x,t+\tau)-\mu(t+\tau)\Big) \cdot\varphi(x) dxd\tau \bigg|
	&\le& \|\varphi\|_{L^\infty(\Omega)} \cdot
	\int_0^1 \big\| d u(\cdot,t+\tau) - \mu(t+\tau)\big\|_{L^1(\Omega)} d\tau \\
	&=& \|\varphi\|_{L^\infty(\Omega)} \cdot \int_t^{t+1} \|du(\cdot,s)-\mu(s)\|_{L^1(\Omega)} ds \\
	&\le& \|\varphi\|_{L^\infty(\Omega)} \cdot \int_t^\infty \|du(\cdot,s)-\mu(s)\|_{L^1(\Omega)}^2 ds \\
	&\le& \frac{\eta}{3}
	\qquad \mbox{for all } t\ge t_0.
  \eas
  In summary, (\ref{104.5}) implies that
  \bas
	\bigg| \io \Big(du(\cdot,t)-\mu(t)\Big)\cdot\varphi \bigg| \le \eta
	\qquad \mbox{for all } t\ge t_0
  \eas
  and thereby yields (\ref{104.1}).
\qed
}
\subsection{Weak $L^1$ convergence of $u$}\label{sect8.3}
{\cred
We are now in the position to address the claimed convergence statement concerning the quantity $u$ itself.
As a last preparation, let us use Lemma \ref{lem104} and again the uniform integrability of
$(u(\cdot,t))_{t>0}$ to derive the following.
\begin{lem}\label{lemc}
  Let $\mu$ be as in (\ref{mu}). Then for each $\varphi\in L^\infty(\Omega)$,
  \be{c1}
	\io u(\cdot,t)\varphi - \mu(t) \io \frac{\varphi}{d} \to 0
	\qquad \mbox{as } t\to\infty.
  \ee
\end{lem}
\proof
  Observing that $|\{d\le \nu\}|\to 0$ as $\nu\searrow 0$, 
  for fixed $\varphi\in L^\infty(\Omega)$ and $\eta>0$ we first employ (\ref{333.2}) to pick $\nu>0$ small enough such
  that
  \be{c2}
	\|\varphi\|_{L^\infty(\Omega)} \cdot \int_{\{d\le \nu\}} u(\cdot,t) \le \frac{\eta}{3}
	\qquad \mbox{for all } t>0
  \ee
  and such that moreover
  \be{c3}
	c_1 \|\varphi\|_{L^\infty(\Omega)} \cdot \int_{\{d\le\nu\}} \frac{1}{d} \le \frac{\eta}{3},
  \ee
  with $c_1:=\sup_{t>0} \mu(t)$ being finite according to Lemma \ref{lem922}.
  As $\frac{\varphi}{d} \cdot\chi_{\{d>\nu\}}$ belongs to $L^\infty(\Omega)$, we may now rely on Lemma \ref{lem104}
  in choosing $t_0>0$ suitably large such that
  \be{c4}
	\bigg| \int_{\{d>\nu\}} \Big(du(\cdot,t)-\mu(t)\Big)\cdot \frac{\varphi}{d} \bigg|
	\le \frac{\eta}{3}
	\qquad \mbox{for all } t\ge t_0.
  \ee
  Then in the identity
  \be{c5}
	\io u(\cdot,t)\varphi - \mu(t) \io \frac{\varphi}{d}
	= \int_{\{d\le \nu\}} u(\cdot,t)\varphi
	+ \int_{\{d>\nu\}} \Big( du(\cdot,t)-\mu(t)\Big) \cdot \frac{\varphi}{d}
	- \mu(t) \int_{\{d\le\nu\}} \frac{\varphi}{d},
	\quad t>0,
  \ee
  we may use (\ref{c2}) to estimate
  \bas
	\bigg| \int_{\{d\le\nu\}} u(\cdot,t)\varphi \bigg|
	\le \|\varphi\|_{L^\infty(\Omega)} \cdot \int_{\{d\le\nu\}} u(\cdot,t) 
	\le \frac{\eta}{3}
	\qquad \mbox{for all } t>0,
  \eas
  and apply (\ref{c3}) to see that
  \bas
	\bigg| -\mu(t)\cdot\int_{\{d\le\nu\}} \frac{\varphi}{d} \bigg|
	\le \mu(t) \|\varphi\|_{L^\infty(\Omega)} \cdot \int_{\{d\le\nu\}} \frac{1}{d}
	\le \frac{\eta}{3}
	\qquad \mbox{for all } t>0.
  \eas
  In view of (\ref{c4}), from (\ref{c5}) we thus infer that
  \bas
	\bigg| \io u(\cdot,t)\varphi - \mu(t)\cdot\io \frac{\varphi}{d} \bigg|
	\le\eta
	\qquad \mbox{for all } t\ge t_0
  \eas
  and conclude.
\qed
}
{\cred Two applications thereof now yield the claimed stabilization property of $u$.
}
\begin{lem}\label{lem815}
  Let $\mu_\infty$ be as specified in Theorem \ref{theo7000}.
  Then 
  \be{815.1}
	u(\cdot,t) \wto \frac{\mu_\infty}{d}
	\quad \mbox{in } L^1(\Omega)
	\qquad \mbox{as } t\to\infty.
  \ee
\end{lem}
{\cred
\proof
  A first application of Lemma \ref{lemc} shows that due to (\ref{massu}), with $\mu$ as in (\ref{mu}) we have
  \bas
	\io u_0 - \mu(t) \cdot \io \frac{1}{d}
	= \io u(\cdot,t) - \mu(t) \cdot \io \frac{1}{d}
	\to 0
	\qquad \mbox{as } t\to\infty
  \eas
  and that hence by definition of $\mu_\infty$,
  \bas
	\mu(t)\to \mu_\infty
	\qquad \mbox{as } t\to\infty.
  \eas
  Therefore, once more employing Lemma \ref{lemc}, now with arbitrary $\varphi\in L^\infty(\Omega)$, shows that
  for any such $\varphi$,
  \bas
	\io u(\cdot,t)\varphi
	= \bigg\{ \io u(\cdot,t)\varphi - \mu(t)\cdot\io \frac{\varphi}{d} \bigg\} 
	+ \mu(t)\cdot \io \frac{\varphi}{d}
	\to \mu_\infty \io \frac{\varphi}{d}
	\qquad \mbox{as } t\to\infty,
  \eas
  and that thus (\ref{815.1}) holds.
\qed
}

\subsection{Uniform decay of $w$. Proof of Theorem \ref{theo7000}}\label{sect8.4}
Now since we already know that $u(\cdot,t)$ stabilizes 
with respect to the weak topology in $L^1(\Omega)$ to a positive limit function as $t\to\infty$,
thanks to the equicontinuity feature of $w$ expressed {\cred in Lemma \ref{lem333}} 
the integral decay property (\ref{256.3}) can be used to derive the following.
\begin{lem}\label{lem67}
  We have
  \be{67.1}
	\int_t^{t+1} \io w \to 0
	\qquad \mbox{as } t\to\infty.
  \ee
\end{lem}
\proof
  Given $\eta>0$, relying on the fact that $u_0\not\equiv 0$ and that hence the number $\mu_\infty$
  in Theorem \ref{theo7000} is positive, we can fix $\delta>0$ small enough such that with
  $c_1:=\frac{\mu_\infty}{\|d\|_{L^\infty(\Omega)}}$ we have
  \be{67.2}
	\delta\cdot\io u_0 \le \frac{c_1 \eta}{6}.
  \ee
  We thereafter once again make use of {\cred Lemma \ref{lem333} which in conjunction with the Arzel\`a-Ascoli theorem   
  ensures that the set $(w(\cdot,t))_{t>0}$}
  is relatively compact in $C^0(\bar\Omega)$, implying that there exist $k_0\in\N$ and $(w_k)_{k\in\{1,...,k_0\}}
  \subset C^0(\bar\Omega)$ with the property that for all $t>0$ one can choose $k(t)\in \{1,...,k_0\}$
  fulfilling
  \be{67.3}
	\|w(\cdot,t)-w_{k(t)}\|_{L^\infty(\Omega)} \le \delta.
  \ee
  Since $\{1,...,k_0\}$ is finite, thanks to the fact that $u(\cdot,t)\wto \frac{\mu_\infty}{d}$ in $L^1(\Omega)$
  as $t \to \infty$, as asserted by Lemma \ref{lem815}, it is then possible to pick
  $t_1>0$ such that
  \be{67.4}
	\bigg| \io u(\cdot,t) w_k - \io \frac{\mu_\infty}{d}w_k \bigg|
	\le \frac{c_1 \eta}{6}
	\qquad \mbox{for all $k\in \{1,...,k_0\}$ and each } t>t_1.
  \ee
  Finally, the integrability property (\ref{256.3}) enables us to find $t_0>t_1$ such that
  \be{67.5}
	\int_{t_0}^\infty \io uw \le \frac{c_1 \eta}{6},
  \ee
  and we claim that these choices guarantee that
  \be{67.6}
	\int_t^{t+1} \io w \le \eta
	\qquad \mbox{for all } t> t_0.
  \ee
  To verify this, we split
  \bea{67.7}
	\int_t^{t+1} \io u(x,s)w(x,s)dxds	
	&=& \int_t^{t+1} \io \frac{\mu_\infty}{d(x)} w(x,s) dxds \nn\\
	& & + \int_t^{t+1} \io u(x,s)\cdot \Big\{ w(x,s)-w_{k(s)}(x)\Big\} dxds \nn\\
	& & + \int_t^{t+1} \bigg\{ \io u(x,s) w_{k(s)}(x) dx - \io \frac{\mu_\infty}{d(x)} w_{k(s)}(x) dx \bigg\} ds \nn\\
	& & + \int_t^{t+1} \io \frac{\mu_\infty}{d(x)} \cdot \Big\{ w_{k(s)}(x)-w(x,s)\Big\} dxds
	\qquad \mbox{for } t>0
  \eea
  and use (\ref{massu}) together with (\ref{67.3}) and (\ref{67.2}) to see that
  \bas
	\bigg| \int_t^{t+1} \io u(x,s)\cdot \Big\{ w(x,s)-w_{k(s)}(x)\Big\} dxds \bigg|
	&\le& \int_t^{t+1} \bigg\{ \io u(\cdot,s) \bigg\} \cdot \|w(\cdot,s)-w_{k(s)}\|_{L^\infty(\Omega)} ds \\
	&=& \bigg\{ \io u_0 \bigg\} \cdot \int_t^{t+1} \|w(\cdot,s)-w_{k(s)}\|_{L^\infty(\Omega)} ds \\
	&\le& \bigg\{ \io u_0 \bigg\}\cdot\delta \\
	&\le& \frac{c_1\eta}{6}
	\qquad \mbox{for all } t>0,
  \eas
  and that, similarly, by definition of $\mu_\infty$ we have
  \bas
	\bigg| \int_t^{t+1} \io \frac{\mu_\infty}{d(x)} \cdot \Big\{ w_{k(s)}(x)-w(x,s)\Big\} dxds \bigg|
	&\le& \mu_\infty \cdot \bigg\{ \io \frac{1}{d} \bigg\} \cdot 
		\int_t^{t+1} \|w_{k(s)}-w(\cdot,s)\|_{L^\infty(\Omega)} ds \\
	&\le& \mu_\infty \cdot \bigg\{ \io \frac{1}{d} \bigg\} \cdot \delta \\
	&\le& \frac{c_1\eta}{6}
	\qquad \mbox{for all } t>0.
  \eas
  As moreover (\ref{67.4}) along with our restriction $t_0>t_1$ ensures that
  \bas
	\bigg| \int_t^{t+1} \bigg\{ \io u(x,s) w_{k(s)}(x) dx - \io \frac{\mu_\infty}{d(x)} w_{k(s)}(x) dx \bigg\} ds \bigg|
	&\le& \int_t^{t+1} \max_{k\in\{1,...,k_0\}} \bigg| \io u(\cdot,s)w_k - \io \frac{\mu_\infty}{d} w_k\bigg| ds \\
	&\le& \frac{c_1\eta}{6}
	\qquad \mbox{for all } t>t_0,
  \eas
  from (\ref{67.7}) we altogether obtain that
  \bas
	\int_t^{t+1} \io uw
	\ge \int_t^{t+1} \io \frac{\mu_\infty}{d} w
	- \frac{c_1\eta}{6}
	- \frac{c_1\eta}{6}
	- \frac{c_1\eta}{6}
	= \int_t^{t+1} \io \frac{\mu_\infty}{d} w
	- \frac{c_1\eta}{2}
	\qquad \mbox{for all } t>t_0.
  \eas
  Since apart from that
  \bas
	\io \frac{\mu_\infty}{d} w \ge \frac{\mu_\infty}{\|d\|_{L^\infty(\Omega)}} \io w = c_1 \io w
	\qquad \mbox{for all } t>0,
  \eas
  combined with (\ref{67.5}) this shows that
  \bas
	c_1 \int_t^{t+1} \io w
	\le \int_t^{t+1} \io uw + \frac{c_1\eta}{2} 
	\le \frac{c_1\eta}{2} + \frac{c_1\eta}{2}
	= c_1\eta
	\qquad \mbox{for all } t>t_0
  \eas
  and thereby establishes (\ref{67.6}).
\qed
Together with the monotonicity information (\ref{256.2}), this entails decay of $w(\cdot,t)$ with respect to the norm
in $L^1(\Omega)$.
\begin{lem}\label{lem68}
  We have
  \be{68.1}
	\|w(\cdot,t)\|_{L^1(\Omega)} \to 0
	\qquad \mbox{as } t\to\infty.
  \ee
\end{lem}
\proof
  {\cred Since from (\ref{256.2}) we know that
  \bas
	\|w(\cdot,t)\|_{L^1(\Omega)} \le \int_{t-1}^t \|w(\cdot,s)\|_{L^1(\Omega)} ds
	\qquad \mbox{for all } t\ge 1,
  \eas
  this is a direct consequence of Lemma \ref{lem67}.
}
\qed
Again by Lemma \ref{lem333}, the topological information herein can be improved.
\begin{lem}\label{lem69}
  We have
  \be{69.1}
	\|w(\cdot,t)\|_{L^\infty(\Omega)} \to 0
	\qquad \mbox{as } t\to\infty.
  \ee
\end{lem}
\proof
  If this was false, there would exist $(t_k)_{k\in\N}\subset (0,\infty)$,
  $(x_k)_{k\in\N}\subset\bar\Omega$ and $c_1>0$ such that
  \bas
	w(x_k,t_k)\ge c_1
	\qquad \mbox{for all } k\in\N,
  \eas
  which due to the equicontinuity of $(w(\cdot,t_k))_{k\in\N}$ asserted by Lemma \ref{lem333} would entail that with some
  $\delta>0$ we would have
  \bas
	w(x,t_k) \ge \frac{c_1}{2}
	\qquad \mbox{for all } x\in (x_k-\delta,x_k+\delta)\cap \Omega
	\ \mbox{and each } k\in\N.
  \eas
  This, however, would be incompatible with the outcome of Lemma \ref{lem68}.
\qed
We have thereby actually already completed the derivation of our main results concerning the 
large time behavior in (\ref{0}).\abs
\proofc of Theorem \ref{theo7000}. \quad
  We only need to combine Lemma \ref{lem815} with Lemma \ref{lem69}.
\qed

\mysection{Instantaneous blow-up. Proof of Theorem \ref{theo87}}\label{sect9}
Finally concerned with the verification of Theorem \ref{theo87}, we will pursue a strategy based on the additional
dissipative structure expressed in the identity
\be{ln1}
	\frac{d}{dt} \io \frac{1}{d} \ln u
	= \io \frac{(du)_x^2}{(du)^2} 
	- \io \frac{(du)_x}{du} w_x
\ee
formally associated with (\ref{0}).
In order to appropriately cope with the latter summand herein, 
even at the level of approximate solutions the preparation of a spatio-temporal estimate for $w_x$
seems in order.
In the limit problem (\ref{0}), this could formally be obtained in a trivial manner under our assumption
that $\frac{w_0}{d}$ and hence $\frac{w}{d}$ be bounded, together with the boundedness of $\io d\frac{w_x^2}{w}$ 
implied by (\ref{energy}). 
At the level of approximate solutions, however, in view of diffusion-induced positivity of $\weps$
considerable additional efforts seem necessary to guarantee appropriate boundedness properties of $\frac{\weps}{\deps}$. 
Our approach toward this will therefore be restricted to the derivation of corresponding $L^p$ bounds for large but
finite $p$ only (Lemma \ref{lem81}), 
thereby requiring to involve additional higher-order regularity features of $\wepsx$, possibly
depending on $\eps$ in a singular manner (Lemma \ref{lem80}), 
to achieve the desired $L^2$ estimate through an interpolation argument (Lemma \ref{lem82}).
Thereafter, on the basis of a regularized counterpart of (\ref{ln1}) we will see in Section \ref{sect9.2}
that our hypothesis that $\io \frac{1}{d} \ln \frac{1}{d}$ be finite, by guaranteeing boundedness of 
the functional $\io \frac{1}{d}\ln u$ in (\ref{ln1}) from above
(Lemma \ref{lem84}), allows for deducing space-time $L^2$ bounds on $(\ln (\deps\ueps))_x$ (Lemma \ref{lem835})
and hence for deriving Theorem \ref{theo87}.
\subsection{An $L^2$ estimate for $\wepsx$ implied by boundedness of $\frac{w_0}{d}$}
Let us first interpolate between two regularity estimates for $\wepsx$ from Lemma \ref{lem8} to achieve
the following bound involving a high integrability power but a singular dependence on $\eps$.
\begin{lem}\label{lem80}
  There exists $C>0$ with the property that for any choice of $\epsin$,
  \be{80.1}
	\int_0^T \bigg\| \deps \frac{\wepsx(\cdot,t)}{\sqrt{\weps(\cdot,t)}} \bigg\|_{L^\infty(\Omega)}^4 dt
	\le \frac{C}{\eps} \cdot (1+T) \cdot (1+\sqrt{\eps}T)^2
	\qquad \mbox{for all } T>0.
  \ee
\end{lem}
\proof
  From Lemma \ref{lem8} we obtain $c_1>0$ and $c_2>0$ such that
  \bas
	\eps \int_0^T \io \frac{1}{\sqrt{g(\weps)}} \cdot \Big(\deps \frac{\wepsx}{\sqrt{g(\weps)}}\Big)_x^2 
	\le c_1 \cdot (1+\sqrt{\eps} T) 
	\qquad \mbox{for all } T>0
  \eas
  and
  \bas
	\io \deps \frac{\wepsx^2}{\weps} \le c_2 \cdot (1+\sqrt {\eps} t)
	\qquad \mbox{for all } t>0,
  \eas
  which in view of (\ref{g2}), (\ref{M}) and (\ref{43.3}) entails that
  \bea{80.2}
	\eps \int_0^T \io \Big(\deps \frac{\wepsx}{\sqrt{g(\weps)}} \Big)_x^2
	&\le& \sqrt{\og M} \cdot \eps \int_0^T \io 
		\frac{1}{\sqrt{g(\weps)}} \cdot \Big(\deps\frac{\wepsx}{\sqrt{g(\weps)}}\Big)_x^2 \nn\\
	&\le& c_3 \cdot (1+\sqrt {\eps} T) 
	\qquad \mbox{for all } T>0
  \eea
  and
  \bea{80.3}
	\io \Big(\deps \frac{\wepsx}{\sqrt{g(\weps)}}\Big)^2
	&=& \io \frac{\deps\weps}{g(\weps)} \cdot \deps \frac{\wepsx^2}{\weps} \nn\\
	&\le& \frac{\|d\|_{L^\infty(\Omega)}+1}{\ug} \io \deps \frac{\wepsx^2}{\weps} \nn\\
	&\le& c_4 \cdot (1+\sqrt{\eps}t)
	\qquad \mbox{for all } t>0
  \eea
  with obvious choices of $c_3>0$ and $c_4>0$.
  Now since the Gagliardo-Nirenberg inequality provides $c_5>0$ fulfilling
  \bas
	\|\varphi\|_{L^\infty(\Omega)}^4 \le c_5 \|\varphi_x\|_{L^2(\Omega)}^2 \|\varphi\|_{L^2(\Omega)}^2
	+ c_5 \|\varphi\|_{L^2(\Omega)}^4
	\qquad \mbox{for all } \varphi\in W^{1,2}(\Omega),
  \eas
  combining (\ref{80.2}) with (\ref{80.3}) we infer that
  \bas
	\int_0^T \bigg\|\deps \frac{\wepsx(\cdot,t)}{\sqrt{g(\weps(\cdot,t))}} \bigg\|_{L^\infty(\Omega)}^4 dt
	&\le& c_5 \int_0^T \bigg\| \Big(\deps \frac{\wepsx(\cdot,t))}{\sqrt{g(\weps(\cdot,t))}}\Big)_x \bigg\|_{L^2(\Omega)}^2
	\bigg\| \deps \frac{\wepsx(\cdot,t)}{\sqrt{g(\weps(\cdot,t))}} \bigg\|_{L^2(\Omega)}^2 dt \\
	& & + c_5\int_0^T \bigg\| \deps \frac{\wepsx(\cdot,t)}{\sqrt{g(\weps(\cdot,t))}}\bigg\|_{L^2(\Omega)}^4 dt \\
	&\le& c_5 \cdot \frac{c_3}{\eps} (1+\sqrt{\eps} T) \cdot c_4 (1+\sqrt{\eps}T)
	+ c_5 \cdot c_4^2 (1+\sqrt{\eps}T)^2 T
	\qquad \mbox{for all } T>0,
  \eas
  which readily implies (\ref{80.1}) due to the fact that
  \bas
	\frac{1}{\sqrt{g(\weps)}} \ge \frac{1}{\sqrt{\og}\sqrt{\weps}}
	\qquad \mbox{in } \Omega\times (0,\infty)
  \eas
  by (\ref{g2}).
\qed
Next, and independently from essentially all our previous analysis,
a testing procedure applied to the second equation in (\ref{0eps}) yields the following weighted $L^p$ estimate
for $\frac{\weps}{\deps}$ for asymptotically large but yet finite $p$.
\begin{lem}\label{lem81}
  Assume that $\frac{w_0}{d}\in L^\infty(\Omega)$.
  Then there exists $C>0$ such that whenever $\epsin$,
  \be{81.1}
	\io \frac{\weps^p(\cdot,t)}{\deps^{p+1}} \le \Big(C\cdot (1+t)\Big)^{2p}
	\qquad \mbox{for all $t>0$ and any } p\in \Big[2,\frac{2}{\sqrt{\eps}}\Big].
  \ee
\end{lem}
\proof
  We integrate by parts in the second equation in (\ref{0eps}) and use the nonnegativity of $\ueps g(\weps)$ as well
  as Young's inequality to obtain
  \bea{81.2}
	\frac{1}{p} \frac{d}{dt} \io \frac{\weps^p}{\deps^{p+1}}
	&=& \io \frac{\weps^{p-1}}{\deps^{p+1}} w_{\eps t} \nn\\
	&\le& \eps \io \frac{\weps^{p-1}}{\deps^{p+1}} \cdot \Big(\deps\frac{\wepsx}{\sqrt{g(\weps)}}\Big)_x \nn\\
	&=& -(p-1)\eps \io \frac{1}{\deps^p} \cdot \frac{\weps^{p-2}}{\sqrt{g(\weps)}} \wepsx^2 
	+ (p+1)\eps \io \frac{\depsx}{\deps^{p+1}} \cdot \frac{\weps^{p-1}}{\sqrt{g(\weps)}} \wepsx \nn\\
	&\le& \frac{(p+1)^2 \eps}{4(p-1)} \io \frac{\depsx^2}{\deps^{p+2}} \cdot \frac{\weps^p}{\sqrt{g(\weps)}}
	\qquad \mbox{for all } t>0.
  \eea
  Here since according to (\ref{43.99}) we have
  \bas
	\frac{\depsx^2}{\deps^2} \le \frac{1}{\sqrt{\eps}}
	\qquad \mbox{in } \Omega,
  \eas
  due to (\ref{g2}) and our restriction $p\ge 2$ we can estimate
  \bea{81.3}
	\frac{(p+1)^2 \eps}{4(p-1)} \io \frac{\depsx^2}{\deps^{p+2}} \cdot \frac{\weps^p}{\sqrt{g(\weps)}}
	&\le& \frac{(2p)^2\eps}{4\cdot\frac{p}{2} \cdot \sqrt{\ug}} \cdot 
		\io \frac{\depsx^2}{\deps^{p+2}} \cdot \weps^{p-\frac{1}{2}} \nn\\
	&=& \frac{2p\eps}{\sqrt{\ug}} \cdot 
		\io \frac{\depsx^2}{\deps^2} \cdot \frac{\weps^\frac{2p-1}{2}}{\deps^p} \nn\\
	&\le& \frac{2p\sqrt{\eps}}{\sqrt{\ug}} \cdot \io \frac{\weps^\frac{2p-1}{2}}{\deps^p}
	\qquad \mbox{for all } t>0.
  \eea
  Since from (\ref{43.3}) we know that
  \bas
	\frac{\weps^\frac{2p-1}{2}}{\deps^p}
	&=& \deps^\frac{p-1}{2p} \cdot \Big(\frac{\weps^p}{\deps^{p+1}}\Big)^\frac{2p-1}{2p} \\
	&\le& (\|d\|_{L^\infty(\Omega)}+1)^\frac{p-1}{2p} \cdot \Big(\frac{\weps^p}{\deps^{p+1}}\Big)^\frac{2p-1}{2p} \\
	&\le& c_1 \Big(\frac{\weps^p}{\deps^{p+1}}\Big)^\frac{2p-1}{2p}
	\qquad \mbox{in } \Omega\times (0,\infty)
  \eas
  with $c_1:=(\|d\|_{L^\infty(\Omega)}+1)^\frac{1}{2}$, in view of the H\"older inequality we see that (\ref{81.3})
  entails the inequality
  \bas
	\frac{(p+1)^2\eps}{4(p-1)} \io \frac{\depsx^2}{\deps^{p+2}} \cdot \frac{\weps^p}{\sqrt{g(\weps)}}
	&\le& \frac{2c_1 p\sqrt{\eps}}{\sqrt{\ug}} \io \Big(\frac{\weps^p}{\deps^{p+1}}\Big)^\frac{2p-1}{2p} \\
	&\le& \frac{2c_1 p\sqrt{\eps}}{\sqrt{\ug}} |\Omega|^\frac{1}{2p} \cdot
	\bigg\{ \io \frac{\weps^p}{\deps^{p+1}} \bigg\}^\frac{2p-1}{2p} \\
	&\le& c_2 p\sqrt{\eps} \cdot 
	\bigg\{ \io \frac{\weps^p}{\deps^{p+1}} \bigg\}^\frac{2p-1}{2p}
	\qquad \mbox{for all } t>0,
  \eas
  where $c_2:=\frac{2c_1}{\sqrt{\ug}} \cdot \max \{|\Omega|^\frac{1}{4},1\}$.
  Therefore, (\ref{81.2}) shows that
  \bas
	\yeps(t):=\io \frac{\weps^p(\cdot,t)}{\deps^{p+1}}, \qquad t\ge 0,
  \eas
  satisfies
  \bas
	\yeps'(t) \le c_2 p^2 \sqrt{\eps} \yeps^\frac{2p-1}{2p}(t)
	\qquad \mbox{for all } t>0,
  \eas
  which on integration implies that
  \bas
	\yeps(t) \le \Big\{ \yeps^\frac{1}{2p}(0) + \frac{c_2}{2} p\sqrt{\eps}t \Big\}^{2p}
	\qquad \mbox{for all } t>0,
  \eas
  that is,
  \be{81.4}
	\io \frac{\weps^p}{\deps^{p+1}}
	\le \Bigg\{ \bigg\{ \io \frac{w_{0\eps}^p}{\deps^{p+1}} \bigg\}^\frac{1}{2p}
	+ \frac{c_2}{2} p\sqrt{\eps} t\Bigg\}^{2p}
	\qquad \mbox{for all } t>0.
  \ee
  Here thanks to the fact that from (\ref{43.6}) we know that
  \bas
	\frac{\eps^\frac{1}{4}}{\deps} \le 1
	\qquad \mbox{in } \Omega,
  \eas
  according to (\ref{43.01}) and the definition of $w_{0\eps}$ in Lemma \ref{lem45} we can use
  {\cred Lemma \ref{lem42} as well as}
  our assumption that $\frac{w_0}{d}$ be bounded in $\Omega$ to see that writing 
  $c_3:=\max\Big\{ 1 \, , \, (\io \frac{1}{d})^\frac{1}{4}\Big\}$ we have
{\cred
  \bas
	\bigg\{ \io \frac{w_{0\eps}^p}{d_{\eps}^{p+1}} \bigg\}^\frac{1}{2p}
	&\le& \Big\| \frac{w_{0\eps}}{d_{\eps}} \Big\|_{L^\infty(\Omega)}^\frac{1}{2} \cdot
	\bigg\{ \io \frac{1}{d_{\eps}} \bigg\}^\frac{1}{2p} \\
	&\le& c_3 \Big\| \frac{w_{0\eps}}{d_{\eps}} \Big\|_{L^\infty(\Omega)}^\frac{1}{2} \\
	&\le& c_3 \Big\| \frac{w_{0} + \eps^\frac{1}{4}}{d_{\eps}} \Big\|_{L^\infty(\Omega)}^\frac{1}{2} \\
	&\le& c_3 \Big\| \frac{w_{0}}{d_{\eps}}\Big\|_{L^\infty(\Omega)}^\frac{1}{2}
	+ c_3 \Big\| \frac{\eps^\frac{1}{4}}{d_{\eps}}\Big\|_{L^\infty(\Omega)}^\frac{1}{2} \\
	&\le& c_4:=c_3 \Big\| \frac{w_0}{d}\Big\|_{L^\infty(\Omega)}^\frac{1}{2} + c_3
	\qquad \mbox{for all } \cred \epsin,
  \eas
  because $w_{0\eps} \le w_0+\eps^\frac{1}{4}$ for any such $\eps$.
}
  As moreover $p\sqrt{\eps}\le 2$ due to our hypothesis, from (\ref{81.4}) we thus infer that
  \bas
	\io \frac{\weps^p}{\deps^{p+1}} \le \Big\{ c_4 + c_2 t \Big\}^{2p}
	\qquad \mbox{for all } t>0
  \eas
  and that hence (\ref{81.1}) holds with $C:=\max\{c_2,c_4\}$.
\qed
Fortunately, the largest admissible $p$ in (\ref{81.1}) is such that in the course of an interpolation argument,
an estimation of the $L^2$ norm in question only involves powers of the inequality in (\ref{80.1}) which are 
such that the singular dependence on $\eps$ therein disappears in the limit $\eps=\eps_j\searrow 0$.
\begin{lem}\label{lem82}
  Assume that $\frac{w_0}{d}\in L^\infty(\Omega)$.
  Then for all $T>0$ one can find $C(T)>0$ such that
  \bas
	\int_0^T \io \wepsx^2 \le C(T)
	\qquad \mbox{for all } \epsin.
  \eas
\end{lem}
\proof
  Let us first apply Lemma \ref{lem8}, Lemma \ref{lem80} and Lemma \ref{lem81} to fix constants $c_1\ge 1, c_2>0$ and $c_3>0$
  such that for all $\epsin \subset (0,1)$ we have
  \be{82.2}
	\io \deps \frac{\wepsx^2}{\weps} \le c_1 \cdot (1+t)
	\qquad \mbox{for all } t>0
  \ee
  and
  \be{82.3}
	\int_0^T \bigg\|\deps \frac{\wepsx(\cdot,t)}{\sqrt{\weps(\cdot,t)}} \bigg\|_{L^\infty(\Omega)}^4 dt
	\le \frac{c_2}{\eps} \cdot (1+T)^3
	\qquad \mbox{for all } T>0
  \ee
  as well as
  \be{82.4}
	\io \frac{\weps^\frac{2}{\sqrt{\eps}}}{\deps^{\frac{2}{\sqrt{\eps}}+1}} 
	\le \Big(c_3 \cdot (1+t)\Big)^\frac{4}{\sqrt{\eps}}
	\qquad \mbox{for all } t>0.
  \ee
  Then invoking the H\"older inequality we see that
  \bas
	\io \wepsx^2
	&=& \io \Big| \deps\frac{\wepsx}{\sqrt{\weps}} \Big|^{\sqrt{\eps}}
	\cdot \Big| \deps \frac{\wepsx^2}{\weps} \Big|^\frac{2-\sqrt{\eps}}{2} 
	\cdot \frac{\weps}{\deps^{1+\frac{\sqrt{\eps}}{2}}} \\
	&\le& \Big\| \deps\frac{\wepsx}{\sqrt{\weps}}\Big\|_{L^\infty(\Omega)}^{\sqrt{\eps}}
	\cdot \bigg\{ \io \deps \frac{\wepsx^2}{\weps} \bigg\}^\frac{2-\sqrt{\eps}}{2}
	\cdot \bigg\{ \io \frac{\weps^\frac{2}{\sqrt{\eps}}}{\deps^{\frac{2}{\sqrt{\eps}}+1}} \bigg\}^\frac{\sqrt{\eps}}{2} \\
	&\le& c_1^\frac{2-\sqrt{\eps}}{2} (1+t)^\frac{2-\sqrt{\eps}}{2}
	\cdot \Big(c_3\cdot (1+t)\Big)^2 \cdot 
	\bigg\|\deps\frac{\wepsx}{\sqrt{\weps}}\bigg\|_{L^\infty(\Omega)}^{\sqrt{\eps}} \\
	&\le& c_1 c_3^2 (1+t)^3 \Big\|\deps\frac{\wepsx}{\sqrt{\weps}}\Big\|_{L^\infty(\Omega)}^{\sqrt{\eps}}
	\qquad \mbox{for all } t>0,
  \eas
  because $c_1\ge 1$.
  In order to make use of (\ref{82.3}) here, we integrate with respect to the time variable and once more employ
  the H\"older inequality to find that again since $(\eps_j)_{j\in\N}\subset (0,1)$,
  \bas
	\int_0^T \io \wepsx^2
	&\le& c_1 c_3^2 \int_0^T (1+t)^3 
	\bigg\| \deps \frac{\wepsx(\cdot,t)}{\sqrt{\weps(\cdot,t)}} \bigg\|_{L^\infty(\Omega)}^{\sqrt{\eps}} dt \\
	&\le& c_1 c_3^2 (1+T)^3 \int_0^T 
	\bigg\| \deps \frac{\wepsx(\cdot,t)}{\sqrt{\weps(\cdot,t)}} \bigg\|_{L^\infty(\Omega)}^{\sqrt{\eps}} dt \\
	&\le& c_1 c_3^2 (1+T)^3 \cdot \Bigg\{ \int_0^T 
	\bigg\| \deps \frac{\wepsx(\cdot,t)}{\sqrt{\weps(\cdot,t)}} \bigg\|_{L^\infty(\Omega)}^4 dt 
		\Bigg\}^\frac{\sqrt{\eps}}{4} 
	\cdot T^\frac{4-\sqrt{\eps}}{4}\\
	&\le& c_1 c_3^2 (1+T)^3 \cdot \Big\{ \frac{c_2}{\eps} \cdot (1+T)^3 \Big\}^\frac{\sqrt{\eps}}{4} 
	\cdot T^\frac{4-\sqrt{\eps}}{4}\\
	&\le& c_1 c_2^\frac{1}{4} c_3^2 (1+T)^4 \eps^{-\frac{\sqrt{\eps}}{4}}
	\qquad \mbox{for all } T>0.
  \eas
  Consequently, the proof can be completed by the observation that
  \bas
	\eps^{-\frac{\sqrt{\eps}}{4}} = e^{-\frac{\sqrt{\eps}}{2}\ln\sqrt{\eps}} \le e^\frac{1}{2e}
	\qquad \mbox{for all } \epsin
  \eas
  due to the fact that $\xi\ln\xi>-\frac{1}{e}$ for all $\xi>0$.
\qed
\subsection{A bound for $|\ln (\deps\ueps)|$. Proof of Theorem \ref{theo87}}\label{sect9.2}
In order to prepare our estimates for the absolute value of $\ln (\deps\ueps)$, let us first make sure that
this quantity cannot attain large negative values throughout $\Omega$.
\begin{lem}\label{lem841}
  There exists $C>0$ such that for each $\epsin$,
  \be{841.1}
	\io \deps\ueps(\cdot,t) \ge C
	\qquad \mbox{for all } t\in \Big(0,\frac{1}{\sqrt{\eps}}\Big).
  \ee
\end{lem}
\proof
  From Lemma \ref{lem826} we obtain $\delta>0$ such that for all $\epsin$,
  \bas
	\int_E \ueps \le \frac{1}{2} \io u_0
	\quad \mbox{for all } t\in \Big(0,\frac{1}{\sqrt{\eps}}\Big)
	\mbox{ and each measurable $E\subset\Omega$ such that } |E|\le\delta,
  \eas
  so that since by integrability of $\frac{1}{d}$ we can find $\nu>0$ fulfilling
  $|\{d\le \nu\}| \le \delta$, we infer that
{\cred
  \bas
	\int_{\{d\le\nu\}} \ueps
	\le \frac{1}{2} \io u_0
	\qquad \mbox{for all } t\in \Big(0,\frac{1}{\sqrt{\eps}}\Big).
  \eas
  Again using that $\deps\ge d$ {\cb and (\ref{massu})}, we obtain that indeed
  \bas
	\io \deps\ueps 
	\ge \io d\ueps
	\ge \nu \cdot \int_{\{d\ge\nu\}} \ueps
	= \nu\cdot \bigg\{ \io \ueps - \int_{\{d\le\nu\}} \ueps \bigg\}
	\ge \frac{\nu}{2} \io \ueps
	{\cb \ = \ \frac{\nu}{2} \io u_0}
  \eas
  for any such $\eps$ and $t$.
}
\qed
In view of (\ref{mass}), this entails an upper bound for the spatial minimum
of $|\ln (\deps\ueps)|$.
\begin{lem}\label{lem836}
  Let $T>0$. Then there exists $C>1$ with the property that for all $\epsin$ and any $t\in (0,T)$ one can find
  $x_0(t,\eps)\in\bom$ such that
  \be{836.1}
	\frac{1}{C} \le \deps(x_0(t,\eps)) \ueps(x_0(t,\eps),t) \le C.
  \ee
\end{lem}
\proof
  Since Lemma \ref{lem841} along with (\ref{mass}) and (\ref{43.3}) ensures the existence of $c_1>0$ and $c_2>0$ such that
  for all $\epsin$ we have
  \bas
	c_1 \le \io \deps\ueps \le c_2
	\qquad \mbox{for all } t\in (0,T),
  \eas
  by a mean-value theorem we can pick some $x_0=x_0(t,\eps)\in\bom$ fulfilling
  \bas
	\deps(x_0)\ueps(x_0,t) 
	= \frac{1}{|\Omega|} \io \deps\ueps(\cdot,t)
	\in \Big[\frac{c_1}{|\Omega|}, \frac{c_2}{|\Omega|} \Big],
  \eas
  so that the claim results on taking $C>1$ suitably large.
\qed
Now a straightforward application of Young's inequality yields the following inequality which inter alia
entails an upper bound for the functional on the right of (\ref{ln1}) at the approximate level.
\begin{lem}\label{lem84}
  Suppose that
  \be{84.01}
	\io \frac{1}{d} \ln \frac{1}{d}<\infty.
  \ee
  Then there exists $C>0$ such that for each $\epsin$,
  \be{84.1}
	\int_{\{\ueps(\cdot,t)\ge 1\}} \frac{1}{\deps} \ln \ueps(\cdot,t) \le C
	\qquad \mbox{for all } t>0.
  \ee
\end{lem}
\proof
  As $\xi\eta \le \xi\ln \xi + e^{\eta-1}$ for all $\xi>0$ and $\eta\in\R$, we may use (\ref{mass}) to see that for 
  all $\epsin$,
  \bea{84.2}
	\int_{\{\ueps(\cdot,t)\ge 1\}} \frac{1}{\deps} \ln \ueps(\cdot,t)
	&\le& \io \frac{1}{\deps} \ln \frac{1}{\deps} + \frac{1}{e} \io \ueps(\cdot,t) \nn\\
	&=& \io \frac{1}{\deps} \ln \frac{1}{\deps} + \frac{1}{e} \io u_0
	\qquad \mbox{for all } t>0.
  \eea
  Since by {\cred monotonicity of $(1,\infty)\ni \xi\mapsto \xi\ln \xi$ and} Lemma \ref{lem43} we have
  \bas
	\io \frac{1}{\deps} \ln \frac{1}{\deps}
	\le \int_{\{\deps<1\}} \frac{1}{\deps} \ln \frac{1}{\deps}
	\le c_1:=\int_{\{d<1\}} \frac{1}{d} \ln \frac{1}{d}
	\qquad \mbox{for all } \epsin,
  \eas
  and since $c_1$ is finite according to our assumption (\ref{84.01}) and the boundedness of $d$, this already yields 
  (\ref{84.1}).
\qed
{\cred
In exploiting the regularized variant of (\ref{ln1}), we shall moreover make use of the following elementary
lemma concerned with an ODE comparison.
}
\begin{lem}\label{lem834}
  Let $T>0$, and suppose that $y\in C^0([0,T))\cap C^1((0,T))$ is such that
  \be{834.1}
	y'(t) + ay_+^2(t) \le h(t)
	\qquad \mbox{for all } t\in (0,T)
  \ee
  with some $a>0$ and some nonnegative $h\in L^1((0,T))\cap C^0((0,T))$. Then 
  \be{834.2}
	y(t) \le \frac{1}{at} + \int_0^t h(s)ds
	\qquad \mbox{for all } t\in (0,T).
  \ee
\end{lem}
\proof
  {\cred
  Since the expression on the right-hand side of (\ref{834.2}) defines a supersolution of the problem
  in (\ref{834.1}) which diverges to $+\infty$ as $t\searrow 0$, 
  this readily results from an ODE comparison argument.
}
\qed
We are now prepared for our analysis of the quasi-dissipative structure suggested by (\ref{ln1}),
relying on the assumption that $\frac{1}{d}$ belong to $L\log L(\Omega)$ through Lemma \ref{lem84}.
\begin{lem}\label{lem835}
  Suppose that
  \be{835.01}
	\io \frac{1}{d} \ln \frac{1}{d}< \infty
  \ee
  and that $\frac{w_0}{d}\in L^\infty(\Omega)$.
  Then for all $T>0$ and any $\tau\in (0,T)$ there exists $C(T,\tau)>0$ such that
  \be{835.02}
	\io \frac{1}{\deps} \ln \ueps(\cdot,t) \ge - C(T,\tau)
	\qquad \mbox{for all } t\in (\tau,T)
  \ee
  and
  \be{835.1}
	\int_\tau^T \io \frac{(\deps\ueps)_x^2}{(\deps\ueps)^2} \le C(T,\tau)
  \ee
  whenever $\epsin$.
\end{lem}
\proof
  We multiply the first equation in (\ref{0eps}) by $\frac{1}{\deps\ueps}$ and integrate by parts over $\Omega$ to see that
  \bas
	\frac{d}{dt} \io \frac{1}{\deps} \ln \ueps
	&=& \io \frac{1}{\deps\ueps} u_{\eps t} \\
	&=& \io \frac{(\deps\ueps)_x^2}{(\deps\ueps)^2} 
	- \io \frac{(\deps\ueps)_x}{\deps\ueps} \wepsx 
	\qquad \mbox{for all } t>0,
  \eas
  where by Young's inequality,
  \bas
{\cred
	\io \frac{(\deps\ueps)_x}{\deps\ueps} \wepsx 
	\le \frac{1}{2} \io \frac{(\deps\ueps)_x^2}{(\deps\ueps)^2} + \frac{1}{2} \io \wepsx^2
	\qquad \mbox{for all } t>0,
}
  \eas
  so that
  \bas
	\zeps(t):= - \io \frac{1}{\deps} \ln \ueps(\cdot,t), \qquad t\ge 0,
  \eas
  satisfies
  \be{835.2}
	\zeps'(t) + \frac{1}{2} {\cred \io} \frac{(\deps\ueps)_x^2}{(\deps\ueps)^2} 
	\le \frac{1}{2} \io \wepsx^2
	\qquad \mbox{for all } t>0.
  \ee
  Now given $T>0$, we apply Lemma \ref{lem836} to gain $c_1>1$ such that for all $\epsin$ and each $t\in (0,T)$
  we can pick $x_0=x_0(t,\eps)\in\bom$ such that
  \bas
	\frac{1}{c_1} \le \deps(x_0)\ueps(x_0,t) \le c_1,
  \eas
  which implies that
  \bas
	\Big| \ln \Big(\deps(x_0)\ueps(x_0,t)\Big) \Big| \le c_2:=\ln c_1.
  \eas
  Since by means of {\cred the Cauchy-Schwarz inequality} we obtain that
  \bas
	\Big| \ln \Big(\deps(x)\ueps(x,t)\Big)
	- \ln \Big(\deps(x_0)\ueps(x_0,t)\Big) \Big|
	\le \sqrt{|\Omega|} \cdot 
	\bigg\{ \io \bigg| \bigg(\ln \Big(\deps\ueps(\cdot,t)\Big)\bigg)_x \bigg|^2 \bigg\}^\frac{1}{2}
	\qquad \mbox{for all } x\in\Omega,
  \eas
  this entails that
  \bas
	\Big| \ln \Big(\deps(x) \ueps(x,t)\Big) \Big|
	\le c_2 + \sqrt{|\Omega|} \cdot
	\bigg\{ \io \frac{(\deps\ueps)_x^2}{(\deps\ueps)^2} \bigg\}^\frac{1}{2}
	\qquad \mbox{for all {\cred $x\in\Omega$ and}  } t\in (0,T)
  \eas
  and that hence
  \bea{835.3}
	\frac{1}{4} \io \frac{(\deps\ueps)_x^2}{(\deps\ueps)^2}
	&\ge& {\cred
	\frac{1}{4|\Omega|} \cdot \bigg\{ \Big\| \ln (\deps\ueps)\Big\|_{L^\infty(\Omega)} - c_2 \bigg\}_+^2} \nn\\
	&\ge& \frac{1}{8|\Omega|} \Big\| \ln (\deps\ueps)\Big\|_{L^\infty(\Omega)}^2 - c_3
	\qquad \mbox{for all } t\in (0,T)
  \eea
  with $c_3:=\frac{c_2^2}{4|\Omega|}$, because {\cred $(\xi-\eta)_+^2 \ge \frac{1}{2}\xi^2 - \eta^2$ for all
  $\xi\ge 0$ and $\eta\ge 0$.}
  Now since again using that $\xi\ln\xi\ge-\frac{1}{e}$ for all $\xi>0$ we can estimate
  \bas
	\zeps(t)
	&=& - \io \frac{1}{\deps} \ln (\deps\ueps) 
	+ \io \frac{1}{\deps} \ln \deps \\
	&\le& - \io \frac{1}{\deps} \ln (\deps\ueps) 
	+ \frac{|\Omega|}{e}
	\qquad \mbox{for all } t>0,
  \eas
  and since Lemma \ref{lem43} warrants that
  \bas
	- \io \frac{1}{\deps} \ln (\deps\ueps)
	&\le& \Big\|\ln (\deps\ueps)\Big\|_{L^\infty(\Omega)} \cdot \io \frac{1}{\deps} \\
	&\le& c_4 \Big\|\ln (\deps\ueps)\Big\|_{L^\infty(\Omega)}
	\qquad \mbox{for all } t>0
  \eas
  with $c_4:=\io \frac{1}{d}<\infty$, from (\ref{835.3}) we thus infer that
  \bas
	\frac{1}{4} \io \frac{(\deps\ueps)_x^2}{(\deps\ueps)^2}
	\ge \frac{1}{8c_4^2|\Omega|} \cdot \Big\{ \zeps(t)-\frac{|\Omega|}{e}\Big\}_+^2 - c_3
	\qquad \mbox{for all } t\in (0,T).
  \eas
  Consequently, writing $c_5:=\frac{1}{8c_4^2 |\Omega|}$ we see that (\ref{835.2}) entails the inequality
  \be{835.4}
	\zeps'(t) + c_5 \cdot \Big\{\zeps(t)-\frac{|\Omega|}{e}\Big\}_+^2
	+ \frac{1}{4} \io \frac{(\deps\ueps)_x^2}{(\deps\ueps)^2}
	\le \frac{1}{2} \io \wepsx^2(\cdot,t) + c_3
	\qquad \mbox{for all } t\in (0,T),
  \ee
  from which in view of Lemma \ref{lem834} we firstly conclude that
  \bas
	\zeps(t) - \frac{|\Omega|}{e} \le \frac{1}{c_5 t} + \frac{1}{2} \int_0^t \io \wepsx^2 + c_3 t
	\qquad \mbox{for all } t\in (0,T).
  \eas
  Since Lemma \ref{lem82} provides $c_6>0$ such that
  \be{835.5}
	\int_0^T \io \wepsx^2 \le c_6
	\qquad \mbox{for all } \epsin,
  \ee
  for arbitrary $\tau\in (0,T)$ and each $\epsin$ this entails the one-sided inequality
  \be{835.6}
	\zeps(t) \le c_7:=\frac{|\Omega|}{e} + \frac{1}{c_5\tau} + \frac{c_6}{2} + c_3 T
	\qquad \mbox{for all } t\in [\tau,T),
  \ee
  thus particularly establishing (\ref{835.02}).\abs
  In order to achieve a corresponding upper bound, we now make use of our assumption (\ref{835.01}), which allows us
  to invoke Lemma \ref{lem84} to find $c_8>0$ fulfilling
  \be{835.7}
	-\zeps(t)
	\le \int_{\{\ueps(\cdot,t)\ge 1\}} \frac{1}{\deps} \ln \ueps(\cdot,t)
	\le c_8
	\qquad \mbox{for all } t>0.
  \ee
  Therefore, namely, on integrating (\ref{835.4}) and relying on (\ref{835.6}) and again (\ref{835.5}) we see that
  \bas
	\frac{1}{4} \int_\tau^T \io \frac{(\deps\ueps)_x^2}{(\deps\ueps)^2}
	&\le& \zeps(\tau)-\zeps(T) + \frac{1}{2} \int_\tau^T \io \wepsx^2 + c_3(T-\tau) \\
	&\le& c_7+c_8+\frac{c_6}{2} + c_3(T-\tau)
	\qquad \mbox{for all } \epsin,
  \eas
  and that thus also (\ref{835.1}) is valid.
\qed
In order to turn this into a two-sided estimate for the quantity $\ln (\deps\ueps)$ itself, we once more
rely on Lemma \ref{lem84} to assert a spatial $L^1$ bound therefor.
\begin{lem}\label{lem85}
  Assume that $\io \frac{1}{d} \ln \frac{1}{d}<\infty$ and $\frac{w_0}{d}\in L^\infty(\Omega)$.
  Then for all $T>0$ and $\tau\in (0,T)$ there exists $C(T,\tau)>0$ such that for all
  $\epsin$ we have
  \be{85.1}
	\io \Big|\ln \Big(\deps\ueps(\cdot,t)\Big) \Big| \le C(T,\tau)
	\qquad \mbox{for all } t\in (\tau,T).
  \ee
\end{lem}
\proof
  In the inequality
  \be{85.2}
	\io \Big|\ln (\deps\ueps)\Big|
	\le \io |\ln\deps| + \io |\ln \ueps|,
	\qquad t>0,
  \ee
  we may first use that the validity of $\ln\xi\le \xi$ for all $\xi>0$ entails that $|\ln\xi| \le \xi+\frac{1}{\xi}$
  for all $\xi>0$, so that according to Lemma \ref{lem43}, writing $c_1:=\|d\|_{L^\infty(\Omega)}+1$ we have
  \be{85.3}
	\io |\ln\deps| \le \io \deps + \io \frac{1}{\deps} \le c_1|\Omega| + \io \frac{1}{d}<\infty.
  \ee
  Likewise, in
  \be{85.4}
	\io |\ln\ueps| = \int_{\{\ueps\ge 1\}} \ln \ueps - \int_{\{\ueps<1\}} \ln\ueps,
	\qquad t>0,
  \ee
  we have
  \be{85.5}
	\int_{\{\ueps\ge 1\}} \ln\ueps
	\le \int_{\{\ueps\ge 1\}} \ueps \le\io \ueps=\io u_0
	\qquad \mbox{for all } t>0
  \ee
  by (\ref{mass}), whereas
  \bea{85.6}
	- \int_{\{\ueps<1\}} \ln\ueps
	&=& - \int_{\{\ueps<1\}} \deps\cdot\frac{1}{\deps}\ln\ueps \nn\\
	&\le& -c_1 \int_{\{\ueps<1\}} \frac{1}{\deps} \ln\ueps \nn\\
	&=& -c_1 \io \frac{1}{\deps}\ln \ueps
	+ c_1 \int_{\{\ueps\ge 1\}} \frac{1}{\deps} \ln \ueps
	\qquad \mbox{for all } t>0.
   \eea
  Since Lemma \ref{lem84} provides $c_2>0$ such that
  \bas
	\int_{\{\ueps\ge 1\}} \frac{1}{\deps}\ln\ueps \le c_2
	\qquad \mbox{for all } t>0,
  \eas
  and since Lemma \ref{lem835} says that given any $T>0$ and $\tau\in (0,T)$ we can find {\cred $c_3(T,\tau)>0$}
  fulfilling
  \bas
	\io \frac{1}{\deps} \ln\ueps \ge -{\cred c_3(T,\tau)}
	\qquad \mbox{for all } t\in (\tau,T),
  \eas
  from (\ref{85.4}), (\ref{85.5}) and (\ref{85.6}) we conclude that
  \bas
	\bigg| \io \ln\ueps \bigg| \le \io u_0 + c_1 {\cred c_3(T,\tau)} + c_1 c_2
	\qquad \mbox{for all } t\in (\tau,T),
  \eas
  which together with (\ref{85.2}) and (\ref{85.3}) verifies (\ref{85.1}).
\qed
Now by interpolation, the latter in conjunction with Lemma \ref{lem835} entails (\ref{87.4}).
\begin{lem}\label{lem86}
  Assume that $\io \frac{1}{d} \ln \frac{1}{d}<\infty$ and $\frac{w_0}{d}\in L^\infty(\Omega)$.
  Then
  \be{86.1}
	\int_\tau^T \Big\| \ln \Big(du(\cdot,t)\Big) \Big\|_{L^\infty(\Omega)}^3 dt < \infty
	\qquad \mbox{for all $T>0$ and } \tau\in (0,T).
  \ee
\end{lem}
\proof
  Given $T>0$ and $\tau\in (0,T)$, from Lemma \ref{lem835} and Lemma \ref{lem85} we obtain $c_1>0$ and $c_2>0$ such that
  for all $\epsin$,
  \be{86.2}
	\int_\tau^T \io \Big(\ln (\deps\ueps)\Big)_x^2 \le c_1
  \ee
  and
  \be{86.3}
	\io \Big|\ln (\deps\ueps)\Big| \le c_2
	\qquad \mbox{for all } t\in (\tau,T).
  \ee
  As a Gagliardo-Nirenberg inequality says that with some $c_3>0$ we have
  \bas
	\|\varphi\|_{L^\infty(\Omega)}^3
	\le c_3 \|\varphi_x\|_{L^2(\Omega)}^2 \|\varphi\|_{L^1(\Omega)}
	+c_3 \|\varphi\|_{L^1(\Omega)}^3
	\qquad \mbox{for all } \varphi\in W^{1,2}(\Omega),
  \eas
  from this we infer that
  \bea{86.4}
	\int_\tau^T \Big\| \ln\Big(\deps\ueps(\cdot,t)\Big) \Big\|_{L^\infty(\Omega)}^3 dt
	&\le& c_3 \int_\tau^T \Big\| \Big(\ln\Big(\deps\ueps(\cdot,t)\Big)\Big)_x \Big\|_{L^2(\Omega)}^2
	\Big\| \ln \Big(\deps\ueps(\cdot,t)\Big)\Big\|_{L^1(\Omega)} dt \nn\\
	& & + c_3 \int_\tau^T \Big\| \ln\Big(\deps\ueps(\cdot,t)\Big)\Big\|_{L^1(\Omega)}^3 dt \nn\\[2mm]
	&\le& c_4:= c_1 c_2 c_3 + c_2^3 c_3 T
  \eea
  for all $\epsin$.
  Now since {\cred Lemma \ref{lem25}} along with Lemma \ref{lem43}
  warrants that with $(\eps_{j_k})_{k\in\N}$ as introduced in Lemma \ref{lem25}, for a.e.~$t>0$ 
  we have $\deps\ueps(\cdot,t)\to du(\cdot,t)$ a.e.~in $\Omega$ and hence
  \bas
	\Big\|\ln \Big(du(\cdot,t)\Big)\Big\|_{L^\infty(\Omega)} 
	\le \liminf_{\eps=\eps_{j_k} \searrow 0} \Big\| \ln \Big(\deps\ueps(\cdot,t)\Big)\Big\|_{L^\infty(\Omega)}
	\qquad \mbox{for a.e. } t>0,
  \eas
  using Fatou's lemma we thus obtain from (\ref{86.4}) that
  \bas
	\int_\tau^T \Big\| \ln \Big(du(\cdot,t)\Big)\Big\|_{L^\infty(\Omega)}^3 dt
	&\le& \int_\tau^T \liminf_{\eps=\eps_{j_k}\searrow 0} 
		\Big\| \ln \Big(\deps\ueps(\cdot,t)\Big)\Big\|_{L^\infty(\Omega)}^3 dt \\
	&\le& \liminf_{\eps=\eps_{j_k}\searrow 0} 
		\int_\tau^T \Big\| \ln \Big(\deps\ueps(\cdot,t)\Big)\Big\|_{L^\infty(\Omega)}^3 dt \\[2mm]
	&\le& c_4
  \eas
  and conclude.
\qed
We thereby readily arrive at our main result on diffusive effects at intermediate time scales.\abs
\proofc of Theorem \ref{theo87}. \quad
  The integrability property (\ref{87.4}) has precisely been asserted by Lemma \ref{lem86}.
  As a consequence, we may choose a null set $N_0\subset (0,\infty)$ such that $\ln (du(\cdot,t))\in L^\infty(\Omega)$
  for all $t\in (0,\infty)\setminus N_0$, whence if for such $t$ we abbreviate 
  $c_1(t):=\|\ln (du(\cdot,t))\|_{L^\infty(\Omega)}$, then
  \bas
	-c_1(t) \le \ln \Big(d(x)u(x,t)\Big) \le c_1(t)
	\qquad \mbox{for a.e. } x\in\Omega,
  \eas
  that is,
  \bas
	\frac{e^{-c_1(t)}}{d(x)} \le u(x,t) \le \frac{e^{c_t(t)}}{d(x)}
	\qquad \mbox{for a.e. } x\in\Omega
  \eas
  whenever $t\in (0,\infty)\setminus N_0$. This yields (\ref{87.5}), whereupon (\ref{87.6}) becomes obvious.
\qed
{\cred
\mysection{Appendix}
This appendix is devoted to the details of the approximation procedures underlying Section \ref{sect2.2}.\\
Let us first construct a family of smooth positive approximations to $d$ with the properties listed
in Lemma \ref{lem43}.\abs
\proofc of Lemma \ref{lem43}. \quad
  Without loss of generality we may assume that $\Omega=(-R,R)$ with some $R>0$, and fix a sequence 
  $(K_j)_{j\in\N}$ of compact subsets of $\{d>0\}$ such that $K_j \subset K_{j+1}$ for all $j\in\N$ and
  $\bigcup_{j\in\N} K_j=\{d>0\}$, whence for $\widetilde K_j:=K_j \cap [-R+\frac{1}{j},R-\frac{1}{j}]$, $j\in\N$,
  we have $\widetilde K_j\subset \widetilde K_{j+1}$ for all $j\in\N$ and $\bigcup_{j\in\N} \widetilde K_j=\{d>0\}\cap\Omega$.
  We first observe that then by continuity of $d$ in $\overline{\Omega}$ and of $d_x$ in $\{d>0\}$, 
  for each $\delta\in (0,1)$,
  \be{43.60}
	\psi_\delta(x):=\left\{ \begin{array}{ll}
	d(-R), \qquad & x\le -\frac{R}{1+\delta}, \\[1mm]
	d\Big((1+\delta)x\Big), \qquad & -\frac{R}{1+\delta} \le x \le \frac{R}{1+\delta}, \\[1mm]
	d(R), \qquad & x\ge\frac{R}{1+\delta},
	\end{array} \right.
  \ee
  defines a function $\psi_\delta\in C^0(\R)$ fulfilling 
  $\psi_\delta \in  
  C^1(\{\psi_\delta>0\} \cap (-\frac{R}{1+\delta}, \frac{R}{1+\delta})) \cap
  W^{1,\infty}_{loc}(\{\psi_\delta>0\})$,
  for which 
  $\psi_\delta\to d$ in $L^\infty(\Omega)$ and $\psi_{\delta x} \to d_x$ in $L^\infty_{loc}(\{d>0\}\cap\Omega)$ 
  and in $L^p_{loc}(\{d>0\})$ for all $p\in [1,\infty)$ 
  as $\delta\searrow 0$, so that for each $j\in\N$ we can pick $\delta_j\in (0,1)$ such that
  $\tphi_j(x):=\psi_{\delta_j}(x), \ x\in\overline{\Omega}$, satisfies
  \be{43.61}
	\|\tphi_j-d\|_{L^\infty(\Omega)}
	\le \frac{1}{2\cdot 3^j},
	\qquad
	\|\widetilde{\varphi}_{jx} - d_x\|_{L^\infty(\widetilde K_j)} \le \frac{1}{2j}
	\qquad \mbox{and} \qquad
	\|\widetilde{\varphi}_{jx} - d_x\|_{L^j(K_j)} \le \frac{1}{2j}
  \ee
  Next, for $\eta\in (0,1)$ letting $\rho_\eta \in C_0^\infty(\R)$ denote an arbitrary mollifier having the properties that
  $\supp \rho_\eta \subset [-\eta,\eta]$ and $\int_{\R} \rho_\eta=1$, we immediately see that
  if $\eta<\frac{\delta_j R}{1+\delta_j}$, then
  $\rho_\eta \star \tphi_j \equiv \tphi_j \equiv d(-R)$ in $(-\infty,-\frac{R}{1+\delta_j}-\eta)$ and
  $\rho_\eta \star \tphi_j \equiv \tphi_j \equiv d(R)$ in $(\frac{R}{1+\delta_j}+\eta,\infty)$ and hence, in particular,
  $(\rho_\eta\star \tphi_j)_x=0$ on $\pO$ for any such $\eta$.
  Since standard arguments (\cite{friedman}) moreover show that $\rho_\eta\star\tphi_j \to \tphi_j$ in $L^\infty(\Omega)$
  as well as $(\rho_\eta \star \tphi_j)_x \to \widetilde{\varphi}_{jx}$ in $L^\infty(\widetilde K_j)$ 	
  and in $L^p(K_j)$ for all $p\in [1,\infty)$ 
  as $\eta\searrow 0$,
  it follows that for any $j\in\N$ we may fix $\eta_j\in (0,1)$ suitably small such that for 
  $\hphi_j:=\rho_ {\eta_j}\star \tphi_j$ we have $\hphi_{jx}=0$ on $\pO$ as well as
  \be{43.62}
	\|\hphi_j - \tphi_j\|_{L^\infty(\Omega)} \le \frac{1}{2 \cdot 3^j},
	\qquad 
	\|\widehat{\varphi}_{jx} - \widetilde{\varphi}_{jx}\|_{L^\infty(\widetilde K_j)} \le \frac{1}{2j}
	\qquad \mbox{and} \qquad
	\|\widehat{\varphi}_{jx} - \widetilde{\varphi}_{jx}\|_{L^j(K_j)} \le \frac{1}{2j}
  \ee
  Writing $\varphi_j:=\hphi_j+\frac{2}{3^j}$, $j\in\N$, we thus obtain $(\varphi_j)_{j\in\N} \subset C^\infty(\overline{\Omega})$
  such that $\varphi_{jx}\equiv \widehat{\varphi}_{jx}$ in $\overline{\Omega}$ and thus still
  \be{43.66}
	\varphi_{jx}=0
	\quad \mbox{on $\pO$ \qquad for all } j\in\N,
  \ee
  that, by (\ref{43.61}) and (\ref{43.62}),
  \be{43.63}
	\|\varphi_{jx} - d_x\|_{L^\infty(\widetilde K_j)} 
	\le \|\widehat{\varphi}_{jx} - \widetilde{\varphi}_{jx} \|_{L^\infty(\widetilde K_j)}
	+ \|\widetilde{\varphi}_{jx} - d_x\|_{L^\infty(\widetilde K_j)}
	\le \frac{1}{2j}+\frac{1}{2j}=\frac{1}{j}
	\qquad \mbox{for all } j\in\N
  \ee
  and similarly
  \be{43.633}
	\|\varphi_{jx} - d_x\|_{L^j(K_j)} 
	\le \|\widehat{\varphi}_{jx} - \widetilde{\varphi}_{jx} \|_{L^j(K_j)}
	+ \|\widetilde{\varphi}_{jx} - d_x\|_{L^j(K_j)}
	\le \frac{1}{2j}+\frac{1}{2j}=\frac{1}{j}
	\qquad \mbox{for all } j\in\N,
  \ee
  and that moreover
  \be{43.7}
	d+\frac{1}{3^j} \le \varphi_j \le d+\frac{3}{3^j}
	\quad \mbox{in } \Omega
	\qquad \mbox{for all } j\in\N,
  \ee
  which in particular ensures that  
  \be{43.8}
	\varphi_{j+1} \le \varphi_j
	\quad \mbox{in } \Omega
	\qquad \mbox{for all } j\in\N.
  \ee
  Now in order to construct $(\deps)_{\eps\in (0,1)}$, we recursively define $(\eps_j)_{j\in\N_0} \subset [0,1]$ by letting
  $\eps_0:=1$ and 
  \be{43.9}
	\eps_j:=\min \Bigg\{ \frac{\eps_{j-1}}{2} \, , \, 3^{-4j} \, , \, 
	\bigg\{ \io \frac{\varphi_{jx}^2}{\varphi_j^3}\bigg\}^{-\frac{1}{2}} \, , \, 
	\bigg\{ \io \frac{\varphi_{jx}^4}{\varphi_j^2}\bigg\}^{-2} \, , \, 
	\Big\| \frac{\varphi_{jx}}{\varphi_j}\Big\|_{L^\infty(\Omega)}^{-4}
	\Bigg\},
	\qquad j\ge 1,
  \ee
  and observe that this especially guarantees that $(\eps_j)_{j\in\N_0}$ is strictly decreasing,
  and that for each $j\in\N$ we have
  $\eps_j>0$ due to (\ref{43.7}) and the inclusion $\varphi_j\in C^1(\overline{\Omega})$.
  As a consequence, introducing 
  \bas
	\deps:=\varphi_j
	\qquad \mbox{whenever $\eps\in (\eps_{j+1},\eps_j]$ for some } j\in\N_0
  \eas
  indeed yields a well-defined family $(\deps)_{\eps\in (0,1)} \subset C^\infty(\overline{\Omega})$ 
  which thanks to (\ref{43.7}),
  (\ref{43.63}), (\ref{43.633}), (\ref{43.8}), (\ref{43.66}) and the monotonicity of $(\eps_j)_{j\in\N}$ satisfies
  (\ref{43.1}), (\ref{43.3}),
  (\ref{43.02}), (\ref{43.01}) and (\ref{43.2}), and for which due to the second restriction expressed
  in (\ref{43.9}) we know from the left inequality in (\ref{43.7}) that for all $j\in\N_0$,
  \bas
	\deps \ge d+\frac{1}{3^j} \ge \frac{1}{3^j} \ge \eps_j^\frac{1}{4} \ge \eps^\frac{1}{4}>0
	\quad \mbox{in } \overline{\Omega}
	\qquad \mbox{for all } \eps \in (\eps_{j+1},\eps_j].
  \eas
  Furthermore, the third, fourth and fifth requirements in (\ref{43.9}) warrant that for any $j\in\N_0$ and each 
  $\eps\in (\eps_{j+1},\eps_j]$ we have
  \bas
	\eps^2 \io \frac{\depsx^2}{\deps^3} 
	= \eps^2 \io \frac{\varphi_{jx}^2}{\varphi_j^3} 
	\le \eps_j^2 \io \frac{\varphi_{jx}^2}{\varphi_j^3} 
	\le 1
  \eas
  and, similarly,
  \bas
	\sqrt{\eps} \io \frac{\depsx^4}{\deps^2} \le\sqrt{\eps_j} \io \frac{\varphi_{jx}^4}{\varphi_j^2} \le 1
  \eas
  as well as
  \bas
	\eps^\frac{1}{4} \Big\| \frac{\depsx}{\deps}\Big\|_{L^\infty(\Omega)}
	\le \eps_j^\frac{1}{4} \Big\|\frac{\varphi_{jx}}{\varphi_j}\Big\|_{L^\infty(\Omega)}
	\le 1,
  \eas
  and that thus also (\ref{43.5}), (\ref{43.55}) and (\ref{43.99}) are valid.
\qed
We next verify that our assumptions on $d$ and $w_0$ indeed entail the 
consequences specified in Lemma \ref{lem41} and Lemma \ref{lem411}.\abs
\proofc of Lemma \ref{lem41}. \quad
  Assuming on the contrary that $c_1:=\int_{\Omega_0} \frac{d_x^2}{d}$ be finite, by hypothesis we can find
  $x_0\in \Omega_0$ and $\delta>0$ such that $d(x_0)=0$ and either $(x_0,x_0+\delta)\subset \Omega_0$ or
  $(x_0-\delta,x_0)\subset \Omega_0$, and concentrating on the former case we know from the continuity of $d$
  that for each $x_1\in \Omega_1:=\{x\in (x_0,x_0+\delta) \ | \ d(x)>0\}$, the point
  $\tilde x_0:=\max \{ x\in [x_0,x_1] \ | \ d(x)=0\}$ belongs to $[x_0,x_1)$.
  As $d$ is positive and hence continuously differentiable on $(\tilde x_0,x_1]$, using elementary calculus we can estimate
  \bas
	\sqrt{d(x_1)}
	&=& \sqrt{d(\tilde x_0)} + \int_{\tilde x_0}^{x_1} (\sqrt{d})_x(y) dy \\
	&=& \frac{1}{2} \int_{\tilde x_0}^{x_1} \frac{d_x(y)}{\sqrt{d(y)}} dy \\
	&\le& \frac{1}{2} \bigg\{ \int_{\tilde x_0}^{x_1} \frac{d_x^2(y)}{d(y)} dy \bigg\}^\frac{1}{2} 
	\cdot \sqrt{x_1-\tilde x_0} \\
	&\le& \frac{\sqrt{c_1}}{2} \sqrt{x_1-\tilde x_0} \\
	&\le& \frac{\sqrt{c_1}}{2} \sqrt{x_1-x_0}.
  \eas
  Since $x_1\in \Omega_1$ was arbitrary and $(x_0,x_0+\delta)\setminus \Omega_1 \subset \{d=0\}$ 
  is a null set by (\ref{d2}), this entails that
  \bas
	\io \frac{1}{d}
	\ge \int_{\Omega_1} \frac{1}{d}
	\ge \int_{\Omega_1} \frac{4}{c_1(x-x_0)} dx
	= \int_{x_0}^{x_0+\delta} \frac{4}{c_1(x-x_0)} dx=\infty,
  \eas
  which in turn is incompatible with (\ref{d2}) and thereby establishes the claim.
\qed
\proofc of Lemma \ref{lem411}. \quad
  Let us assume for contradiction that there exists $x_0\in\overline{\Omega}$ such that $d(x_0)=0$ but $w_0(x_0)>0$.
  Then by continuity of $w_0$ we can find $\delta>0$ and an interval $\Omega_0\subset\overline{\Omega}$, relatively open
  in $\overline{\Omega}$, such that $w_0\ge\delta$ throughout $\Omega_0$. As $d>0$ a.e.~in $\Omega$ as a consequence 
  of (\ref{d2}), using Lemma \ref{lem41} we therefore obtain
  \bas
	\int_{\Omega_0} \frac{d_x^2}{d} w_0
	\ge \delta \int_{\Omega_0} \frac{d_x^2}{d}
	= \infty,
  \eas
  which contradicts (\ref{w0}).
\qed
We are now in the position to provide an approximation of $w_0$ in the flavor of Lemma \ref{lem42}.\abs
\proofc of Lemma \ref{lem42}. \quad
  Without loss of generality we may assume that $\{d=0\}$ is not empty.
  Then since $d$ is continuous in $\overline{\Omega}$, there exist a countable set $I\subset \N$ 
  and a family $(J_i)_{i\in I}$ of relatively 
  open proper
  subintervals $J_i$ of $\overline{\Omega}$ such that $J_i\cap J_j =\emptyset$ if $i\in I$ and $j\in I$ are such that
  $i\ne j$, and that $\bigcup_{i\in I} J_i=\{d>0\}$.
  Accordingly, for each $i\in I$ there exist $a_i\in \overline{\Omega}$ and $b_i\in \overline{\Omega}$ such that
  $(a_i,b_i) \subset J_i\subset [a_i,b_i]$, where $a_i\in J_i$ (resp., $b_i\in J_i$) if and only if 
  $a_i\in\pO$ (resp., $b_i\in\pO$).\\
  Now for fixed $i\in I$, in the case $a_i\not\in J_i$ we know from the defining properties of $J_i$ that
  $d(a_i)=0$, whence again by continuity of $d$ we have $\|d\|_{L^\infty((a_i,a_i+\delta))} \to 0$
  as $\delta\searrow 0$; likewise, if $b_i\not\in J_i$ then
  $\|d\|_{L^\infty((b_i-\delta,b_i))} \to 0$ as $\delta\searrow 0$.
  Therefore, we can recursively define $(\delta_j^{(i)})_{j\in\N} \subset (0,1)$ such that
  \be{42.7}
	\delta_j^{(i)} < \frac{b_i-a_i}{4}
	\qquad \mbox{for all } j\in\N
  \ee
  and
  \be{42.8}
	\delta_{j+1}^{(i)} \le \min \Big\{ \delta_j^{(i)} \, , \, \frac{1}{j} \Big\}
	\qquad \mbox{for all } j\in\N,
  \ee
  and such that if $a_i\not\in J_i$, then
  \be{42.9}
	\|d\|_{L^\infty((a_i,a_i+2\delta_j^{(i)}))} \le \frac{1}{2^i}
	\qquad \mbox{for all } j\in\N,
  \ee
  and that if $b_i\not\in J_i$ then
  \be{42.10}
	\|d\|_{L^\infty((b_i-2\delta_j^{(i)},b_i))} \le \frac{1}{2^i}
	\qquad \mbox{for all } j\in\N.
  \ee
  For $i\in I$ and $j\in\N$, we then introduce the piecewise linear functions
  $\zeta_j^{(i)} \in W^{1,\infty}(\Omega)$ by letting
  \be{42.11}
	\zeta_j^{(i)}:= \left\{ \begin{array}{ll}
	0 \qquad & \mbox{if } x\le a_i+\delta_j^{(i)}, \\[1mm]
	\frac{x-a_i-\delta_j^{(i)}}{\delta_j^{(i)}}
	\qquad & \mbox{if } a_i+\delta_j^{(i)}<x<a_i+2\delta_j^{(i)}, \\[1mm]
	1 & \mbox{if } a_i+2\delta_j^{(i)} \le x \le b_i-2\delta_j^{(i)}, \\[1mm]
	\frac{b_i-\delta_j^{(i)}-x}{\delta_j^{(i)}}
	& \mbox{if } b_i-2\delta_j^{(i)} < x < b_i-\delta_j^{(i)}, \\[1mm]
	0 & \mbox{if } x\ge b_i-\delta_j^{(i)}
	\end{array} \right.
  \ee
  whenever $J_i=(a_i,b_i)$ and
  \be{42.12}
	\zeta_j^{(i)}:= \left\{ \begin{array}{ll}
	1 & \mbox{if } x \le b_i-2\delta_j^{(i)}, \\[1mm]
	\frac{b_i-\delta_j^{(i)}-x}{\delta_j^{(i)}}
	& \mbox{if } b_i-2\delta_j^{(i)} < x < b_i-\delta_j^{(i)}, \\[1mm]
	0 & \mbox{if } x\ge b_i-\delta_j^{(i)}
	\end{array} \right.
  \ee
  in the case $J_i=[a_i,b_i)$ and
  \be{42.13}
	\zeta_j^{(i)}:= \left\{ \begin{array}{ll}
	0 \qquad & \mbox{if } x\le a_i+\delta_j^{(i)}, \\[1mm]
	\frac{x-a_i-\delta_j^{(i)}}{\delta_j^{(i)}}
	\qquad & \mbox{if } a_i+\delta_j^{(i)}<x<a_i+2\delta_j^{(i)}, \\[1mm]
	1 & \mbox{if } x\ge a_i+2\delta_j^{(i)} 
	\end{array} \right.
  \ee
  when $J_i=(a_i,b_i]$, and for $j\in\N$ we let
  \be{42.133}
	\zeta_j(x):=\sum_{i\in I, i\le j} \zeta_j^{(i)}(x),
	\qquad x\in\overline{\Omega},
  \ee
  as well as
  \be{42.14}
	w_{0j}(x):=\zeta_j^2(x) w_0(x),
	\qquad x\in\overline{\Omega}.
  \ee
  Then since (\ref{42.8}) in particular asserts that $\delta_j^{(i)} \searrow 0$ as $j\to\infty$ for each $i\in I$,
  from the definition of $\zeta_j$ it follows that
  \be{42.15}
	0 \le \zeta_j(x) \le \zeta_{j+1}(x)
	\qquad \mbox{for all $x\in\Omega$ and } j\in\N
  \ee
  and
  \be{42.16}
	\zeta_j(x) \nearrow 1
	\quad \mbox{as } j\to\infty
	\qquad \mbox{for all } x\in \{d>0\},
  \ee
  implying that $0 \le w_{0j} \le w_{0,j+1}$ in $\Omega$ for all $j\in\N$, and that both (\ref{42.4}) and (\ref{42.5})
  hold. Moreover, it is clear from (\ref{42.14}) and the inclusion $w_0\in W^{1,2}(\Omega)$ implied by our assumptions
  on $w_0$ that $w_{0j} \in W^{1,2}(\Omega)$ with
  \bas
	w_{0jx}^2
	&=& (\zeta_j^2 w_{0x} + 2\zeta_j \zeta_{jx} w_0)^2 \\
	&\le& 2\zeta_j^4 w_{0x}^2
	+ 8\zeta_j^2 \zeta_{jx}^2 w_0^2
	\qquad \mbox{a.e.~in } \Omega,
  \eas
  so that
  \be{42.145}
	\frac{w_{0jx}^2}{w_{0j}}
	\le 2\zeta_j^2 \frac{w_{0x}^2}{w_0}
	+ 8\zeta_{jx}^2 w_0
	\qquad \mbox{a.e.~in } \Omega.
  \ee
  Since $\frac{w_{0x}^2}{w_0}\in L^1(\Omega)$ by hypothesis, this firstly implies that for each fixed $j\in\N$ we have
  $\frac{w_{0jx}^2}{w_{0j}}\in L^1(\Omega)$ and hence $\sqrt{w_{0j}}\in W^{1,2}(\Omega)$, and according to (\ref{42.15}) and
  (\ref{42.133}), from (\ref{42.145}) we furthermore obtain that
  \bea{42.144}
	\io d\frac{w_{0jx}^2}{w_{0j}}
	&\le& 2\io d\zeta_j^2 \frac{w_{0x}^2}{w_0}
	+ 8 \io d\zeta_{jx}^2 w_0 \nn\\
	&\le& 2\io d\frac{w_{0x}^2}{w_0}
	+ 8\sum_{i\in I} \io d (\zeta_{jx}^{(i)})^2 w_0
	\qquad \mbox{for all } j\in\N.
  \eea
  In order to estimate the rightmost summand herein, we first note that according to our choice of $(J_i)_{i\in I}$,
  for all $i\in I$ we have
  \bas
	d(a_i)=0
	\mbox{ whenever } a_i\not\in J_i
	\qquad \mbox{and} \qquad
	d(b_i)=0
	\mbox{ when } b_i\not\in J_i,
  \eas
  and that thus, as a consequence of (\ref{d2}) and (\ref{w0}) when combined with Lemma \ref{lem411},
  \bas
	w_0(a_i)=0
	\mbox{ if } a_i\not\in J_i
	\qquad \mbox{and} \qquad
	w_0(b_i)=0
	\mbox{ if } b_i\not\in J_i.
  \eas
  Again since $\sqrt{w_0}\in W^{1,2}(\Omega)$, by means of the Cauchy-Schwarz inequality this implies that writing
  $c_1:=\io (\sqrt{w_0})_x^2$ we have
  \bas
	w_0(x) \le c_1 |x-a_i|
	\quad \mbox{for all } x\in\Omega
	\qquad \mbox{if } a_i\not\in J_i
  \eas
  and
  \bas
	w_0(x) \le c_1 |x-b_i|
	\quad \mbox{for all } x\in\Omega
	\qquad \mbox{if } b_i\not\in J_i,
  \eas
  so that whenever $i\in  I$ is such that $J_i=(a_i,b_i)$, in view of (\ref{42.11}) we can use (\ref{42.9}) and (\ref{42.10})
  to estimate
  \bas
	\io d (\zeta_{jx}^{(i)})^2 w_0
	&=& \frac{1}{(\delta_j^{(i)})^2} \int_{a_i+\delta_j^{(i)}}^{a_i+2\delta_j^{(i)}} dw_0
	+ \frac{1}{(\delta_j^{(i)})^2} \int_{b_i-2\delta_j^{(i)}}^{b_i-\delta_j^{(i)}} dw_0 \\
	&\le& \frac{1}{(\delta_j^{(i)})^2} \cdot \delta_j^{(i)} 
		\|d\|_{L^\infty((a_i,a_i+2\delta_j^{(i)}))}
		\|w_0\|_{L^\infty((a_i,a_i+2\delta_j^{(i)}))} \\
	& & \frac{1}{(\delta_j^{(i)})^2} \cdot \delta_j^{(i)} 
		\|d\|_{L^\infty((b_i-2\delta_j^{(i)}, b_i))}
		\|w_0\|_{L^\infty((b_i-2\delta_j^{(i)}, b_i))} \\
	&\le& 2c_1 \|d\|_{L^\infty((a_i,a_i+2\delta_j^{(i)}))}
	+ 2c_1 \|d\|_{L^\infty((b_i-2\delta_j^{(i)}, b_i))} \\
	&\le& 2c_1 \cdot \frac{1}{2^i} + 2c_1 \cdot \frac{1}{2^i} \\
	&=& \frac{4c_1}{2^i}.
  \eas
  Along with a similar reasoning in the cases $J_i=[a_i,b_i)$ and $J_i=(a_i,b_i]$, this allows us to conclude that
  \bas
	8\sum_{i\in I} \io d(\zeta_{jx}^{(i)})^2 w_0
	\le 32c_1 \sum_{i\in I} \frac{1}{2^i}
	\le 32c_1 \sum_{i=1}^\infty \frac{1}{2^i}
	< \infty
	\qquad \mbox{for all } j\in \N,
  \eas
  because $I\subset \N$. In light of our assumption (\ref{init}), from (\ref{42.144}) we thus obtain (\ref{42.6}).
\qed
Our final selection of the sequence $(\eps_j)_{j\in\N}\subset (0,1)$, as used throughout our analysis, can be
accomplished as follows.\abs
\proofc of Lemma \ref{lem45}.\quad
  For fixed $j\in\N$ we estimate
  \be{45.6}
	\io \deps \frac{(w_{0j}+\eps^\frac{1}{4})_x^2}{w_{0j}+\eps^\frac{1}{4}}
	= \io \deps \frac{w_{0jx}^2}{w_{0j}+\eps^\frac{1}{4}}
	\le \io \deps \frac{w_{0jx}^2}{w_{0j}}
	\qquad \mbox{for all } \eps\in (0,1),
  \ee
  where using the inclusion $\sqrt{w_{0j}}\in W^{1,2}(\Omega)$, as asserted by Lemma \ref{lem42}, along with the 
  monotonicity of the convergence $\deps\to d$, as obtained in Lemma \ref{lem43}, we see that
  \bas
	\io \deps \frac{w_{0jx}^2}{w_{0j}} \to \io d\frac{w_{0jx}^2}{w_{0j}}
	\qquad \mbox{as } \eps\searrow 0.
  \eas
  As $c_1:=\sup_{j\in\N} \io d\frac{w_{0jx}^2}{w_{0j}}$ is finite thanks to Lemma \ref{lem42}, from this and (\ref{45.6})
  we infer that for all $j\in\N$ we can fix $\eps^{(1)}(j)\in (0,1)$ such that
  \be{45.7}
	\io \deps \frac{(w_{0j}+\eps^\frac{1}{4})_x^2}{w_{0j}+\eps^\frac{1}{4}} \le c_1+1
	\qquad \mbox{for all } \eps\in (0,\eps^{(1)}(j)].
  \ee
  Next, for arbitrary $j\in\N$ and $\eps\in (0,1)$ we trivially split
  \be{45.8}
	\io \frac{\depsx^2}{\deps} (w_{0j}+\eps^\frac{1}{4})
	= \io \frac{\depsx^2}{\deps} w_{0j}
	+ \eps^\frac{1}{4} \io \frac{\depsx^2}{\deps}
  \ee
  and note that here due to the Cauchy-Schwarz inequality, the boundedness property (\ref{43.55}) derived in 
  Lemma \ref{lem43} ensures that
  \be{45.9}
	\eps^\frac{1}{4} \io \frac{\depsx^2}{\deps} 
	\le \eps^\frac{1}{4} |\Omega|^\frac{1}{2} \bigg\{ \io \frac{\depsx^4}{\deps^2} \bigg\}^\frac{1}{2}
	\le |\Omega|^\frac{1}{2}
	\qquad \mbox{for all } \eps\in (0,1).
  \ee
  Now since $K:=\supp w_{0j}$ is a compact subset of $\{d>0\}$ by Lemma \ref{lem42}, and since according to
  Lemma \ref{lem43} we have $\deps\to d$ in $L^\infty(\Omega)$ and $\depsx \to d_x$ in 
  $L^2_{loc}(\{d>0\})$ and hence $\frac{\depsx^2}{\deps} \to \frac{d_x^2}{d}$ in $L^1(K)$ as $\eps\searrow 0$,
  it follows that for any individual $j\in\N$,
  \bas
	\io \frac{\depsx^2}{\deps} w_{0j} \to \io \frac{d_x^2}{d} w_{0j}
	\qquad \mbox{as } \eps\searrow 0.
  \eas
  Since $w_{0j} \le w_0$ by (\ref{42.5}) and thus
  \bas
	\io \frac{d_x^2}{d} w_{0j} \le c_2:=\io \frac{d_x^2}{d} w_0
	\qquad \mbox{for all } j\in\N
  \eas
  with $c_2$ being finite thanks to our assumptions on $w_0$, we thus conclude that for any $j\in\N$ we can pick
  $\eps^{(2)}(j)\in (0,1)$ fulfilling
  \bas
	\io \frac{\depsx^2}{\deps} w_{0j} \le c_2+1
	\qquad \mbox{for all } \eps\in (0,\eps^{(2)}(j)],
  \eas
  which together with (\ref{45.8}) and (\ref{45.9}) entails that
  \bas
	\io \frac{\depsx^2}{\deps} (w_{0j}+\eps^\frac{1}{4})
	\le c_2+1 + |\Omega|^\frac{1}{2}
	\qquad \mbox{for all } \eps\in (0,\eps^{(2)}(j)].
  \eas
  In conjunction with (\ref{45.7}), this shows that if we pick any $\eps_0\in (0,1)$ and recursively define a nonincreasing
  sequence $(\eps_j)_{j\in\N}\subset (0,1)$ by letting
  \bas
	\eps_j:=\min \Big\{ \eps_{j-1} \, , \, \frac{1}{j} \, , \, \eps^{(1)}(j) \, , \, \eps^{(2)}(j) \Big\},
	\qquad j\in \N,
  \eas
  then $(w_{0\eps_j})_{j\in\N}$ as given by (\ref{45.3}) indeed satisfies (\ref{45.4}) and (\ref{45.5}).
\qed
}
{\bf Acknowledgement.}\quad
{\cred
The author would like to thank Christina Surulescu for her crucial support with regard to the embedding of this work
into the context of glioma invasion.
Furthermore, the author is grateful to Christian Stinner
for numerous useful remarks which substantially improved this manuscript.
Apart from that, the author acknowledges support of 
{\em Deutscher Akademischer Austauschdienst}
within the project {\em Qualitative analysis of models for taxis mechanisms}.
}

\end{document}